\documentclass[12pt]{article}

\title{Intersection cohomology of character varieties for punctured Riemann surfaces}
\author{Mathieu Ballandras \\ \it{Université de Paris} \\ \it{Scuola Internazionale Superiore di Studi Avanzati} \\
\it{Instituto de Ciencias Matemáticas} \\ \small{mballandras@imj-prg.fr} }

\DeclareUnicodeCharacter{0308}{}

\usepackage[english]{babel}

\usepackage[T1]{fontenc}

\usepackage[a4paper,top=3cm,bottom=2cm,left=3cm,right=3cm,marginparwidth=1.75cm]{geometry}

\usepackage{mathtools}

\usepackage[vcentermath]{youngtab}
\usepackage{young}

\usepackage[shortlabels]{enumitem}

\usepackage{amsthm}

\usepackage{relsize,exscale}

\usepackage{amsmath}
\usepackage{amssymb}
\usepackage{graphicx}
\usepackage[colorinlistoftodos]{todonotes}
\usepackage[colorlinks=true, allcolors=blue]{hyperref}
\usepackage{cleveref}
\usepackage{pst-node}
\usepackage{tikz-cd}
\usetikzlibrary{angles, intersections, quotes}

\usepackage{upgreek}

\usepackage{graphicx}

\newcommand{\commentaires}[1]{}

\makeatletter
\newcommand*\bigcdot{\mathpalette\bigcdot@{.5}}
\newcommand*\bigcdot@[2]{\mathbin{\vcenter{\hbox{\scalebox{#2}{$\m@th#1\bullet$}}}}}
\renewenvironment{abstract}{
\null\vfil
\@beginparpenalty\@lowpenalty
\begin{center}
\bfseries \abstractname
\@endparpenalty\@M
\end{center}
}
\makeatother

\newtheorem{definition}{Definition}[section]
\newtheorem{theorem}[definition]{Theorem}
\newtheorem{corollary}[definition]{Corollary}
\newtheorem{lemma}[definition]{Lemma}
\newtheorem{proposition}[definition]{Proposition}

\newtheorem{remark}[definition]{Remark}
\newtheorem{example}[definition]{Example}
\newtheorem{notations}[definition]{Notations}
\newtheorem{conjecture}[definition]{Conjecture}

\DeclareMathOperator{\Spec}{Spec}
\DeclareMathOperator{\Hom}{Hom}
\DeclareMathOperator{\SymF}{Sym}

\DeclareMathOperator{\Id}{Id}
\DeclareMathOperator{\End}{End}
\DeclareMathOperator{\Ind}{Ind}
\DeclareMathOperator{\Res}{Res}

\DeclareMathOperator{\tr}{tr}
\DeclareMathOperator{\Aut}{Aut}

\DeclareMathOperator{\Imm}{Im}

\DeclareMathOperator{\ch}{ch}

\DeclareMathOperator{\GL}{GL}
\DeclareMathOperator{\SL}{SL}

\DeclareMathOperator{\rank}{rank}
\DeclareMathOperator{\reg}{reg}
\DeclareMathOperator{\Exp}{Exp}
\DeclareMathOperator{\Log}{Log}
\DeclareMathOperator{\PGL}{PGL}
\DeclareMathOperator{\gl}{\mathfrak{gl}}

\DeclareMathOperator{\supp}{Supp}

\DeclareMathOperator{\Ree}{Re}

\DeclareMathOperator{\pdeg}{p-deg}

\newcommand{\map}[5]{
\begin{array}{ccccc}
    #1 & : & #2 & \to & #3  \\
     & & #4 & \mapsto & #5
\end{array}
}
\newcommand\KK{\mathbb K }

\newcommand{\C}{\mathbb{C}}
\newcommand{\CC}{\mathbb{C}}
\newcommand{\qlbar}{\overline{\mathbb{Q}}_l}

\newcommand{\IH}{IH}

\newcommand{\sn}{\mathfrak{S}_n}
\newcommand{\Sym}{\mathfrak{S}}
\newcommand{\Fq}{\mathbb{F}_q}
\newcommand{\Fqbar}{\overline{\mathbb{F}}_q}

\newcommand{\NN}{\mathbb N}

\newcommand{\SymX}{\SymF \left[X\right]}
\newcommand{\Part}{\mathcal{P}}

\newcommand{\Ht}{\tilde{H}}

\newcommand{\inv}{^{-1}}

\newcommand{\dcb}[1]{\mathcal{D}_c^b\left(#1\right)}
\newcommand{\hi}{\mathcal H ^i}
\newcommand{\uic}[1]{\underline{\mathcal{I C}}_{#1}^\bullet}
\newcommand{\ic}[1]{\mathcal{I C}_{#1}^\bullet}

\newcommand{\drond}{\overline{\partial}}
\newcommand{\zbar}{\overline{z}}

\newcommand{\mcmus}{\mathcal{M}_{\overline{\mathcal{C}}_{\boldsymbol{\mu},\boldsymbol{\sigma}}}}
\newcommand{\mcmur}{\mathcal{M}_{\overline{\mathcal{C}}_{\boldsymbol{\rho},\boldsymbol{\sigma}}}}
\newcommand{\cmus}{\mathcal{C}_{\boldsymbol{\mu},\boldsymbol{\sigma}}}
\newcommand{\mlps}{\widetilde{\mathcal{M}}_{\boldsymbol{L},\boldsymbol{P},\boldsymbol{\sigma}}}

\newcommand{\mcbar}{\mathcal{M}_{ \boldsymbol{ \overline{ \mathcal{C} } } }}

\begin{document}

\maketitle

\begin{abstract}
    We study intersection cohomology of character varieties for punctured Riemann surfaces with prescribed monodromies around the punctures. Relying on previous result from Mellit \cite{mellit_p} for semisimple monodromies we compute the intersection cohomology of character varieties with monodromies of any Jordan type. This proves the Poincaré polynomial specialization of a conjecture from Letellier \cite{zclosure}.
\end{abstract}

\bigskip

\tableofcontents

\section{Introduction}
Character varieties studied in this article classify rank $n$ local systems over a genus $g$ Riemann surface with $k$-punctures $(p_j)_{1\le j \le k}$. The monodromy around the puncture $p_j$ is imposed to be in the closure $\overline{\mathcal{C}}_j$ of a conjugacy class $\mathcal{C}_j$ of $\GL_n(\C)$. The character variety is an affine variety defined as a geometric invaritant theory quotient:

\begin{multline*}
\mathcal{M}_{ \boldsymbol{ \overline{ \mathcal{C} } } }:=\left\lbrace\left(A_1,B_1,\dots,A_g,B_g,X_1,\dots,X_k\right)\in \GL_n^{2g}\times \overline{\mathcal{C}}_{1}\times\dots\times \overline{\mathcal{C}}_{k}\right| \\ \left. A_1 B_1 A_1\inv B_1\inv \dots A_g B_g A_g\inv B_g\inv X_1\dots X_k =\Id \right\rbrace // \GL_n
\end{multline*}
with $\GL_n$ acting by overall conjugation. A genericity condition is imposed on the $k$-uple of conjugacy classes so that the quotient has good properties (see \ref{def_generic_class}). We study the cohomology of those varieties. As they are not smooth, it is convenient to study their \textit{intersection cohomology}. We compute the Poincaré polynomial for compactly supported intersection cohomology of those character varieties. This Poincaré polynomial encodes the dimension of the compactly supported intersection cohomology spaces $I H_c^r\left(\mathcal{M}_{ \boldsymbol{ \overline{ \mathcal{C} } } }, \qlbar \right)$ as coefficients of a polynomial:
\[
P_c \left( \mcbar, v\right):=\sum_r \dim I H_c^r\left(\mathcal{M}_{ \boldsymbol{ \overline{ \mathcal{C} } } }, \qlbar \right) v^r
\]
When the conjugacy classes are semisimple, they are closed, and the variety $\mcbar$ is smooth. Then the intersection cohomology coincides with the usual cohomology.  Cohomology of character varieties has been extensively studied in various contexts.

\subsection{Cohomology of character varieties: state of the art}
\subsubsection{One puncture with a central monodromy}
A first interesting case is when there is only one puncture and the associated monodromy is central. The genericity condition implies that the monodromy is $e^{-\frac{2 i \pi d}{n}} \Id$ with $d,n$ coprime. Then the character variety is denoted by $\mathcal{M}_{B}^d$. The index $B$ stands for Betti moduli space. Non-Abelian Hodge theory relates this Betti moduli space to a Dolbeault moduli space $\mathcal{M}^d_{Dol}$. This can be seen as a generalization of Narasimhan-Seshadri  \cite{Narasimhan_Seshadri} result relating unitary representations and holomorphic vector bundles. $\mathcal{M}^d_{Dol}$ is is the moduli space of stable Higgs bundles of rank $n$ and degree $d$. Non-Abelian Hodge correspondence was proved in rank $n=2$ by Hitchin \cite{Hitchin87} and Donaldson \cite{Donaldson}.  It was generalized to higher ranks and higher dimensions by Corlette \cite{Corlette} and Simpson \cite{Simpson_88} see also \cite{Simpson1992HiggsBA}. The correspondence is obtained as a homeomorphism between moduli spaces by Simpson \cite{Simpson_94_I,Simpson_94_II}.

Many computations of the cohomology are performed from the Dolbeault side. First Hitchin \cite{Hitchin87} computed the Poincaré polynomial in rank $n=2$. Gothen \cite{Gothen} extended the computation to rank $n=3$. Hausel-Thaddeus \cite{Hausel_Thaddeus_03,hausel_thaddeus_2004} computed the cohomology ring in rank $n=2$. Garc\'ia-Prada, Heinloth, Schmitt \cite{Garcia_Heinloth_Schmitt} gave a recursive algorithm to compute the motive of the Dolbeault moduli space. They computed an explicit expression in rank $n=4$. Garc\'ia-Prada, Heinloth \cite{Garcia_Prada_Heinloth_y_genus} obtained an explicit formula for $y$-genus in any rank.

As in the last examples, there exists more precise cohomological information than the Poincaré polynomial. The character varieties are affine, by Deligne \cite{Deligne_Hodge}, their cohomology carries a mixed-Hodge structure. The non-Abelian Hodge theory does not preserve this mixed-Hodge structure. Indeed the cohomology of the Dolbeault moduli space is pure contrarily to the cohomology of the affine character variety. De Cataldo-Hausel-Migliorini \cite{de_Cataldo_Hausel_Migliorini_PW} conjectured that under non-Abelian Hodge correspondence, the weight filtration coincides with a perverse filtration induced by Hitchin fibration. This is known as the $P=W$ conjecture, they proved it in rank $n=2$. Recently, de Cataldo-Maulik-Shen \cite{decataldo_maulik_shen} proved the conjecture for genus $g=2$ and any rank.

Another interesting aspect of those moduli spaces is the mirror symmetry. Hausel-Thaddeus \cite{Hausel_Thaddeus_mirror01,Hausel_Thaddeus_mirror} conjectured that the moduli space of $\PGL_n$-Higgs bundles and the moduli space of $\SL_n$-Higgs bundles are related by mirror symmetry, see also \cite{hausel_mirror}. This conjecture was proved by Groechenig-Wyss-Ziegler \cite{Groechenig_Wyss_Ziegler_mirror} and a motivic version by Loeser-Wyss \cite{Loeser_Wyss}. Mirror symmetry was also studied in the parabolic case by Biswas-Dey \cite{Biswas_Dey}. Gothen-Oliveira \cite{Gothen_Oliveira} proved a parabolic version of the conjecture, for particular ranks.

An efficient approach to compute cohomological invariant is to count points of algebraic varieties over finite fields. On the Betti side, Hausel and Rodriguez-Villegas \cite{hausel_villegas}  gave a conjectural formula for the mixed-Hodge polynomial of character varieties with one puncture and a central generic monodromy. They proved the $E$-polynomial specialization of the conjecture by counting points over finite fields. With a similar approach, Mereb \cite{Mereb} computed the $E$-polynomial of $\SL_n$ character varieties. Hausel \cite{hausel_mirror} also proposed a conjectural formula for the Hodge polynomial of the associated Dolbeault moduli space. Mozgovoy \cite{mozgovoy_12} extended this conjecture to the motives of the Dolbeault moduli space.

Schiffmann \cite{Schiffmann16}  computed the Poincaré polynomial of the Dolbeault moduli space by counting Higgs bundles over finite fields. In following articles \cite{Mozgovoy_Schiffmann_14,mozgovoy_schiffmann_2020} Mozgovoy-Schiffmann extended this counting to twisted Higgs bundles. Chaudouard-Laumon \cite{chaudouard_laumon_2016} counted Higgs bundles using automorphic forms.

Mellit \cite{Mellit2017} proved that the formula obtained by Schiffmann \cite{Schiffmann16} is equivalent to the Poincaré polynomial specialization of the conjecture of Hausel and Rodriguez Villegas \cite{hausel_villegas}.

Fedorov-Soibelman-Soibelman \cite{Fedorov_Soibelman_Soibelman} computed the motivic class of the moduli stack of semistable Higgs bundles.

\subsubsection{Any number of punctures and arbitrary monodromies}
Logares-Muñoz-Newstead \cite{logares_munoz_newstead} computed the $E$-polynomial of character varieties for $\SL_2$ and small genus $g=1,2$. They consider one puncture with any conjugacy class, without the genericity assumption. They also obtained the Hodge numbers in genus $g=1$.
Logares-Muñoz \cite{Logares_Munoz} extended those results to genus $g=1$ and two punctures. They computed the $E$-polynomials and some Hodge numbers.
Mart\'inez-Muñoz \cite{Martinez_Munoz_14,Martinez_Munoz_14_b} computed the $E$-polynomial of $\SL_2$-character varieties for any genus and any conjugacy class at the puncture. Mart\'inez \cite{martinez2017} then treated the case of $\PGL_2$-character varieties. 

Simpson \cite{Simpson} generalized non-Abelian Hodge theory to character varieties with punctures and arbitrary conjugacy classes. The generalization is even larger as it concerns filtered local systems. They correspond to parabolic Higgs bundles on the Dolbeault side. The moduli space of stable parabolic Higgs bundles was constructed algebraically by Yokogawa \cite{Yokogawa_93}. The moduli spaces were constructed analytically by Konno \cite{Konno} for Higgs fields with nilpotent residues and by Nakajima \cite{Nakajima_1996}. Those analytic constructions provide the non-Abelian Hodge theory as a diffeomorphism. Biquard-Boalch \cite{biquard_boalch_2004} proved a more general wild non-Abelian Hodge theory and constructed the associated moduli spaces. Biquard, García-Prada and Mundet i Riera \cite{Biquard_Garcia_Prada_Riera} generalized filtered non-Abelian Hodge theory to a large family of groups.

On the Dolbeault side of this correspondence, Boden-Yokogawa \cite{Boden_Yokogawa_96} computed the Poincaré polynomial of the moduli space of parabolic Higgs bundles, in rank $n=2$, using Morse theory. Garc\'ia-Prada, Gothen, Muñoz \cite{Garcia_Prada_Gothen_Munoz_parabolic_3} computed the Poincaré polynomial in rank $n=3$.

Hausel, Letellier and Rodriguez-Villegas \cite{aha} made a conjecture for the mixed-Hodge polynomial of character varieties with generic semisimple conjugacy classes at punctures. Counting points of the character variety over finite field they proved the $E$-polynomial specialization. Chuang-Diaconescu-Pan \cite{Chuang_Diaconescu_Pan} and Chuang-Diaconescu-Donagi-Pantev \cite{Chuang_Diaconescu_Donagi_Pantev} proposed a string theoretic interpretation of the conjecture. This string theoretic approach was also applied to wild character varieties by Diaconescu \cite{Diaconescu_17} and Diaconescu-Donagi-Pantev \cite{Diaconescu_Donagi_Pantev}. Another approach uses recursive relations for various genus. It is used by Mozgovoy \cite{mozgovoy_12}, Carlsson and Rodriguez-Villegas \cite{Carlsson_Villegas}. Similarly to this recursive approach, Gonz{\'a}lez-Prieto  \cite{Gonzalez-Prieto} developped a topological quantum field theory associated to character varieties. Fedorov-Soibelman-Soibelman \cite{Fedorov_Soibelman_Soibelman_parabolic} computed the motivic class of the moduli stack of semistable parabolic Higgs bundles.

Mellit \cite{mellit_p} proved the Poincaré polynomial specialization of the conjecture from \cite{aha} by counting parabolic Higgs bundles over finite fields. This result is of the utmost importance for the present article. This is the starting point of the computation of intersection cohomology of the character variety with the closure of any generic conjugacy classes at punctures.

\subsubsection{No punctures}
In absence of punctures, the character variety is singular and corresponds, via non-Abelian Hodge correspondence, to a moduli space of Higgs bundles of degree zero. Counting points over finite fields, Baraglia-Hekmati \cite{baraglia_hekmati} computed the $E$-polynomial of such character varieties in rank $3$. As they are singular it is also interesting to consider their intersection cohomology. Felisetti \cite{felisetti_2021} computed the intersection cohomology in rank $n=2$ and genus $g=2$. Mauri \cite{mauri} generalized to rank $n=2$ and arbitrary genus. Felisetti-Mauri \cite{felisetti_mauri} proved $P=W$ conjecture for intersection cohomology in genus $g=1$ and arbitrary rank $n$, and in genus $g=2$ rank $n=2$. Mauri \cite{Mauri2021} also studied topological mirror symmetry for those varieties, in rank $n=2$.

\subsection{Intersection cohomology of character varieties for punctured Riemann surfaces}
\subsubsection{Poincaré polynomial}
Letellier \cite{zclosure} gave a conjectural formula for the mixed-Hodge polynomial of the character variety $\mathcal{M}_{\boldsymbol{\overline{\mathcal{C}}}}$, with any type of generic conjugacy classes at punctures. This conjecture generalizes the one for semisimple conjugacy classes \cite{aha}. It also involves Hausel-Letellier-Villegas kernel $\mathbb{H}_n^{H L V}$. This kernel lies in
\[
\SymF\left[X_1\right]\otimes \dots \otimes \SymF\left[X_k\right]
\]
with $\SymF\left[X_j\right]$ the space of symmetric functions in the infinite set of variable $X_j$. The definition of the kernel is recalled in \ref{def_omega_hlv}, it uses modified Macdonald polynomials. The Poincaré polynomial specialization of Letellier's conjecture is the following formula
\begin{equation}\label{eq_poincare_intro}
P_c \left(\mcbar; v\right) = v^{d_{\boldsymbol{\mu}}} \left\langle s_{\boldsymbol{\mu'}}, \mathbb{H}_n^{H L V}\left(-1,v\right) \right\rangle.
\end{equation}
$\boldsymbol{\mu}$ encodes the Jordan type of the conjugacy classes (see \eqref{eq_bold_kuple}). $d_{\boldsymbol{\mu}}$ is the dimension of the variety $\mcbar$, the symmetric function $s_{\boldsymbol{\mu'}}$ is a variant of Schur functions, it is defined in \eqref{eq_s_boldmu}. A very interesting feature of this relation is that no matter the $k$-uple of conjugacy classes, the cohomology is encoded in a single object, the kernel $\mathbb{H}_n^{H L V}$. 

Mellit \cite{mellit_p} computed the Poincar{é} polynomial of character varieties with semisimple conjugacy classes. Let $\boldsymbol{\mathcal{S}}=\left(\mathcal{S}_1,\dots,\mathcal{S}_k\right)$ a generic $k$-uple of conjugacy classes. The Jordan type of this $k$-uple is determined by $k$ partitions $\nu^1,\dots,\nu^k$. The parts of the partition $\nu^j$ are the multiplicities of the distinct eigenvalues of $\mathcal{S}_j$. As explained in \ref{subsect_poincare_ss}, Mellit's result is a particular case of the Poincaré polynomial specialization of the conjecture:
\begin{equation}\label{eq_mellit_intro}
P_c \left(\mathcal{M}_{\boldsymbol{\mathcal{S}}};v\right) = v^{d_{\boldsymbol{\nu}}} \left\langle h_{\boldsymbol{\nu}},\mathbb{H}_n^{H L V}(-1,v)\right\rangle.
\end{equation}
With $h_{\boldsymbol{\nu}}$ the symmetric function 
\[
h_{\boldsymbol{\nu}}:=h_{\nu^1}[X_1]\dots h_{\nu^k}[X_k].
\]
The complete symmetric functions $\left(h_{\lambda}[X]\right)_{\lambda \in \Part_n}$ form a basis of the space of symmetric functions of degree $n$. The set of partitions of an integer $n$ is denoted by $\Part_n$. The transition matrices in the space of symmetric functions are well known, for instance they are in Macdonald book \cite{mcd}. Hence we can express $s_{\boldsymbol{\mu'}}$ in terms of $h_{\boldsymbol{\nu}}$. In this article, in order to compute the Poincaré polynomial of a general character variety $\mcbar$, the combinatoric relation between previous symmetric functions is understood in terms of geometric relation between $\mcbar$ and $\mathcal{M}_{\boldsymbol{\mathcal{S}}}$. Letellier obtained such a relation, but between a character variety $\mcbar$ and its resolution.

\subsubsection{Springer theory and resolutions of character varieties}

Logares-Martens \cite{logares_martens} constructed Grothendieck-Springer resolutions for moduli spaces of parabolic Higgs bundles. Letellier \cite{zclosure} constructed resolution of singularities of character varieties 
\[
\mlps\to\mcbar.
\]
Symplectic resolutions of character varieties were also studied in details by Schedler-Tirelli \cite{schedler_tirelli}. The construction of $\mlps$ is recalled in \ref{def_resol_char}, it relies on Springer theory. Springer \cite{Springer1976} proved a correspondence between unipotent conjugacy classes and representations of Weyl groups. Following work of Lusztig \cite{Lusztig81} for the general linear group, Borho-MacPherson \cite{borho} obtained Springer correspondence in terms of intersection cohomology.

Let us briefly recall their result for the Springer resolution of the unipotent locus in $\GL_n$. Let $B$ the subgroup of upper triangular matrices, $U$ the subgroup of $B$ with $1$ on the diagonal. $T$ is the subgroup of diagonal matrices so that $B=T U$. Let $\mathcal{U}$ the set of unipotent elements in $\GL_n$, \textit{i.e.} the set of matrices with all eigenvalues equal to $1$. Then $\mathcal{U}$ is stratified by conjugacy classes $\left(\mathcal{C}_{\lambda}\right)_{\lambda\in\Part_n}$ with $\lambda$ the partition of $n$ with parts specifying the size of the Jordan blocks. Let
\[
\widetilde{\mathcal{U}}=\left\lbrace\right (X, g B)\in\mathcal{U}\times \GL_n/B \left| g\inv X g \in U  \right.  \rbrace 
\]
the projection to the first factor $\widetilde{\mathcal{U}}\to\mathcal{U}$ is a resolution of singularities. Borho-Macpherson approach to Springer theory provides the following relation between cohomology of the resolution $\widetilde{\mathcal{U}}$ and intersection cohomology of the closure of the strata of $\mathcal{U}$
\[
H_c^{r+\dim \widetilde{\mathcal{U}}}\left(\widetilde{\mathcal{U}}, \qlbar \right) \cong \bigoplus_{\lambda \in \Part_n} V_{\lambda} \otimes I H_c^{r+\dim \mathcal{C}_{\lambda}}\left(\overline{\mathcal{C}}_{\lambda},\qlbar \right).
\]
$V_{\lambda}$ is the irreducible representation of the symmetric group indexed by the partition $\lambda$. The indexing is as in Macdonald's book \cite{mcd}, so that $V_{(n)}$ is the trivial representation and $V_{(1^n)}$ the sign. In terms of Poincaré polynomial previous relation becomes
\[
v^{-\dim \widetilde{\mathcal{U}}} P_c\left(\widetilde{\mathcal{U}},v\right) = \sum_{\lambda\in \Part_n}\left( \dim V_{\lambda} \right) v^{-\dim \mathcal{C}_{\lambda}} P_c\left(\overline{\mathcal{C}}_{\lambda},v\right).
\]
Interestingly, this relation between $v^{-\dim \widetilde{\mathcal{U}}} P_c\left(\widetilde{\mathcal{U}},v\right)$ and $v^{-\dim \mathcal{C}_{\lambda}}P_c\left(\overline{\mathcal{C}}_{\lambda},v\right)$ is exactly the base change relation expressing the symmetric function $h_{1^n}$ in terms of Schur functions $\left(s_{\lambda}\right)_{\lambda\in\Part_n}$
\[
h_{1^n} = \sum_{\lambda\in\Part_n} \left( \dim V_{\lambda} \right) s_{\lambda}.
\]
In this simple example, a base change relation between complete symmetric functions and Schur functions has a geometrical interpretation in terms of Springer resolutions. 

For character varieties the idea is similar but a more general theory is necessary. It is provided by Lusztig parabolic induction \cite{Lusztig1984,lusztig_character_sheaves,Lusztig1986}. Letellier applied this theory to obtain relations between cohomology of the resolution $\mlps$ and intersection cohomology of character varieties $\mathcal{M}_{\overline{\mathcal{C}}_{\boldsymbol{\rho},\boldsymbol{\sigma}}}$ (see  \ref{eq_bold_kuple}  and \ref{notations_jmus} for the definition of the $k$-uple of conjugacy classes $\mathcal{C}_{\boldsymbol{\rho},\boldsymbol{\sigma}}$). This was used to prove that various formulations of the conjecture are equivalent \cite[Proposition 5.7]{letellier_gl}.
In terms of Poincaré polynomial the relation reads
\begin{equation}\label{eq_intro_pc_resol}
v^{-d_{\boldsymbol{\mu}}} P_c\left(\mlps,t \right) = \sum_{\boldsymbol{\rho\preceq\mu}} \left(\dim A_{\boldsymbol{\mu',\rho}}\right) v^{-d_{\boldsymbol{\rho}} } P_c\left( \mathcal{M}_{\overline{\mathcal{C}}_{\boldsymbol{\rho},\boldsymbol{\sigma}}} , v \right).
\end{equation}
This geometric relation is discussed in details in \ref{sect_poincare_comput}, it corresponds to a  combinatoric relation between various basis of symmetric functions:
\begin{equation}\label{eq_intro_change_basis}
h_{\boldsymbol{\mu'}} = \sum_{\boldsymbol{\rho\preceq\mu}} \left(\dim A_{\boldsymbol{\mu',\rho}}\right) s_{\boldsymbol{\rho}}.
\end{equation}
In order to generalize Mellit result from semisimple conjugacy classes to any Jordan type, the geometric interpretation of \eqref{eq_intro_change_basis} should involve a character variety with semisimple monodromies $\mathcal{M}_{\boldsymbol{S}}$ instead of a resolution $\mlps$. It will appear that the Poincaré polynomial of the resolution $\mlps$ is equal to the Poincaré polynomial of a character variety with semisimple monodromie $\mathcal{M}_{\boldsymbol{\mathcal{S}}}$. Together with Mellit's result \eqref{eq_mellit_intro}, this implies
\[
v^{-d_{\boldsymbol{\mu}}} P_c\left(\mlps,v \right) = v^{-d_{\boldsymbol{\mu}}} P_c \left(\mathcal{M}_{\boldsymbol{\mathcal{S}}},v\right) = \left\langle h_{\boldsymbol{\mu'}}, \mathbb{H}_n^{H L V}(-1,v) \right\rangle 
\]
Relations \eqref{eq_intro_pc_resol} \eqref{eq_intro_change_basis} can be inverted so that the Poincaré polynomial of a character variety with any type of monodromies can be expressed as Poincaré polynomial of character varieties with semisimple monodromies. This is exactly what is necessary to obtain the general formula \eqref{eq_poincare_intro} from Mellit's result for semisimple conjugacy classes \eqref{eq_mellit_intro}.

To summarize, computing the Poincaré polynomial for intersection cohomology of character varieties requires three elements:
\begin{itemize}
    \item Mellit's result for character varieties with semisimple monodromies \eqref{eq_mellit_intro}.
    \item Letellier's relation \eqref{eq_intro_pc_resol} between cohomology of the resolution $\mlps$ and intersection cohomology of character varieties $\mcbar$.
    \item Relation between cohomology of the resolution $\mlps$ and cohomology of a character variety with semisimple monodromies $\mathcal{M}_{\boldsymbol{\mathcal{S}}}$.
\end{itemize}
The last point is studied in Section \ref{sect_diffeo_res} where a diffeomorphism between the resolution $\mlps$ and a character variety with semisimple monodromies $\mathcal{M}_{\boldsymbol{\mathcal{S}}}$ is detailed so that the Poincaré polynomial coincide. Constructing the diffeomorphism requires analytical technics. They are detailed in \ref{sect_global_moduli}, they rely on the  filtered version of non-Abelian Hodge theory and Riemann-Hilbert correspondence. 
Those correspondences are due to Simpson \cite{Simpson}. The moduli spaces providing non-Abelian Hodge theory as a diffeomorphism were constructed by Konno \cite{Konno}, Nakajima \cite{Nakajima_1996} and Biquard-Boalch \cite{biquard_boalch_2004} in the more general setting of wild non-Abelian Hodge theory. Filtered version of Riemann-Hilbert correspondence is described as a diffeomorphism by Yamakawa \cite{Yamakawa}. A filtered version of non-Abelian Hodge theory was developped for a large family of groups by Biquard, García-Prada and Mundet i Riera \cite{Biquard_Garcia_Prada_Riera}. In Section \ref{sect_diffeo_res} this is used to construct a diffeomorphism between $\mlps$ and $\mathcal{M}_{\boldsymbol{\mathcal{S}}}$, see Theorem \ref{th_diffeo_resol_ss}. The proof of the Poincaré polynomial specialization of Letellier's conjecture is achieved in \ref{sect_poincare_comput}:
\begin{theorem}\label{th_poincaré_intro}
Consider a generic $k$-uple of conjugacy classes ${\mathcal{C}}_{\boldsymbol{\mu},\boldsymbol{\sigma}}$ (notations are introduced in \eqref{eq_bold_kuple}). the Poincaré polynomial for compactly supported intersection cohomology of the character variety $\mathcal{M}_{\overline{{\mathcal{C}}}_{\boldsymbol{\mu},\boldsymbol{\sigma}}}$ is
\[
P_c \left(\mathcal{M}_{\overline{{\mathcal{C}}}_{\boldsymbol{\mu},\boldsymbol{\sigma}}},v\right) = v^{d_{\boldsymbol{\mu}}}\left\langle s_{\boldsymbol{\mu'}}, \mathbb{H}_n^{H L V}(-1,v) \right\rangle .
\]
\end{theorem}

In addition to provide a combinatorial relation between Poincaré polynomials, a fundamental aspect of Springer theory and Lusztig parabolic induction is the action of Weyl group on cohomology spaces.

\subsubsection{Weyl group action on the cohomology of character varieties}
The construction of resolutions of character varieties relies on Springer resolutions and Lusztig parabolic induction. Therefore there is a Weyl group action on the cohomology of resolutions of character varieties (see Letellier \cite{zclosure}). It is interesting to notice that the Weyl group only acts on the cohomology and not on the variety itsel. Another Weyl group action on the cohomology of character varieties and their resolutions is constructed by Mellit \cite{mellit_cell} when $k-1$ among $k$ conjugacy classes are semisimple. It is called \textit{monodromic} Weyl group action.

As explained in previous section, in order to compute the intersection cohomology of character varieties for any conjugacy classes, we construct a diffeomorphism between a resolution $\mlps$ and a character variety with semisimple monodromies $\mathcal{M}_{\boldsymbol{\mathcal{S}}}$. This diffeomorphism allows to move the Springer-like Weyl group action on the cohomology of the resolution, to a Weyl group action on the cohomology of the character varieties with semisimple monodromies $\mathcal{M}_{\boldsymbol{\mathcal{S}}}$. This action is enough for our purpose of computation of the Poincaré polynomial. Moreover, it also provides the $\boldsymbol{\eta}$-twisted Poincaré polynomials, \textit{i.e.} the trace of any elements of the Weyl group on the cohomology spaces, see Definition \ref{def_twisted_Poincare}. Considering a $k$-uple of generic semisimple conjugacy classes $\boldsymbol{\mathcal{S}}=\left(\mathcal{S}_1,\dots,\mathcal{S}_k\right)$, the relative Weyl group is the group permuting eigenvalues with the same multiplicity in a given class $\mathcal{S}_j$. Next theorem is proved in \ref{subsect_twisted_poincare}.
\begin{theorem}\label{th_twisted_intro}
For any $\boldsymbol{\eta}$ conjugacy class in the relative Weyl group, the $\boldsymbol{\eta}$-twisted Poincaré polynomial of the character variety $\mathcal{M}_{\boldsymbol{\mathcal{S}}}$ is
\[
 P_c^{\boldsymbol{\eta}}\left( \mathcal{M}_{\boldsymbol{\mathcal{S}}} ,v\right) :=\sum_{r} \tr\left(\boldsymbol{\eta},H_c^r(\mathcal{M}_{\boldsymbol{\mathcal{S}}}, \qlbar)\right)v^r = (-1)^{r(\boldsymbol{\eta})}v^{d_{\boldsymbol{\mu}}}\left\langle \widetilde{h}_{\boldsymbol{\eta}},\mathbb{H}_n^{H L V}\left(-1,v \right)  \right\rangle.
\]
$r(\boldsymbol{\eta})$ and the symmetric functions $\widetilde{h}_{\boldsymbol{\eta}}$ are defined in \ref{notations_k_h_tilde}.
\end{theorem}

However a more satisfying approach would be to directly construct a monodromic Weyl group action on the cohomology of character varieties with semisimple monodromies. Like the one constructed by Mellit for the $k$-th monodromy \cite{mellit_cell}.

\subsection{Plan}
Section \ref{sect_background} contains reminder and notations about intersection cohomology, symmetric functions and Springer theory. 

The construction of character varieties and their resolutions is recalled in  \ref{sect_background_char}. This section also includes discussion about previous results and conjectures for cohomology of character varieties.

The diffeomorphism between a resolution $\mlps$ and a character variety with semisimple monodromies $\mathcal{M}_{\boldsymbol{\mathcal{S}}}$ is constructed in \ref{sect_diffeo_res}. First for a particularly interesting example (sphere with four punctures and rank $2$) using only algebraic tools. The general case then relies on analytic technics such as non-Abelian Hodge theory and Riemann-Hilbert correspondence. 

The  Poincaré polynomial for intersection cohomology is computed in \ref{sect_computation_pc}. Its twisted version, recording traces of Weyl group action is also given.

\section*{Aknowledgement}

This work is part of my PhD thesis. First of all I want to thank my advisors Emmanuel Letellier and Fernando Rodriguez Villegas for introducing me to the study of character varieties.

Much of this work rely on previous results from Anton Mellit. I am very grateful to him for interesting discussions.

I want to thank Philip Boalch, a fundamental idea in this work was hinted by him. He explained me why non-Abelian Hodge theory could be useful and pointed out the relevant references.

All my gratitude goes to Duiliu-Emanuel Diaconescu and Olivier Schiffmann for their detailed reading of the thesis.

Many thanks to Olivier Biquard, Oscar García-Prada, Penka Georgieva, Tamás Hausel and Florent Schaffhauser for accepting to be part of the jury.

I want to thank Carlos Simpson for pointing out relevant references.

I am grateful to Mirko Mauri for interesting suggestions and references.

\section{Geometric and combinatoric background}\label{sect_background}
\subsection{Perverse sheaves and intersection cohomology}\label{sect_perv_ic}
In this section classical results about perverse sheaves and intersection cohomology are stated. The constructions come from Beilinson, Bernstein, Deligne and Gabber  \cite{bbd}.

$\KK$ is either $\C$ or an algebraic closure $\Fqbar$ of a finite field $\Fq$ with $q$ elements. $X$ is an algebraic variety over $\KK$. Let $l$ be a prime different from the characteristic of $\KK$, the constant $l$-adic sheaf on $X$ with coefficients in $\qlbar$ is denoted by $\kappa_X$ or just $\kappa$ when the context is clear.

\begin{notations}
The category of $\kappa$-constructible sheaves on $X$ is denoted by $\dcb{X}$. Its objects are represented by complexes of sheaves $K$ such that the cohomology sheaves $\mathcal H ^i K$ are  $\kappa$-constructible sheaves on $X$ and finitely many of them are non-zero. For $Y$ a variety over $\KK$ and $f:X\to Y$ a morphism one has the usual functors 
\begin{eqnarray*}
f^* , f^! : \dcb Y  \to \dcb X \\
f_*, f_! : \dcb X \to \dcb Y
\end{eqnarray*}
For $m$ an integer $K[m]$ is the shifted complex such that $\hi K[m] = \mathcal{H}^{i+m}K$. For $x$ a point in $X$, the stalk at $x$ of the $i$-th cohomology sheave of the complex $K$ is denoted by $\hi_x K $ . The structural morphism of $X$ is $p:X\to\Spec\KK$. The $k$-th cohomology space of $X$ with coefficients in $\kappa$ is
\[
H^k(X,\kappa) := \mathcal{H}^k p_* \kappa_X
\]
and the $k$-th compactly supported intersection cohomology space of $X$ is
\[
H^k_c (X,\kappa) := \mathcal{H}^k p_! \kappa_X.
\]
The Verdier dual operator is denoted by $D_X:\dcb X \to \dcb X$.
\end{notations}

\begin{definition}[Perverse sheave]
A perverse sheave is an object $K$ in $\dcb X$ such that for all $i\in\NN$
\begin{eqnarray*}
\dim \left( \supp \hi K\right) &\le& -i \\
\dim \left( \supp \hi D_X K\right) &\le & -i .
\end{eqnarray*}
The category of perverse sheaves on $X$ is denoted by $\mathcal{M}(X)$, it is an abelian category.
\end{definition}

\begin{definition}[Intersection complex]\label{def_intersection_complex}
Let $Y\hookrightarrow X$ a closed embedding and $j:U\hookrightarrow Y$ an open embedding. Assume $U$ is smooth, irreducible and $\overline U = Y$. Let $\xi$ be a local system on $U$. $\underline{\mathcal{I C}}_{Y,\xi}^\bullet$ is the unique perverse sheave $K$ on $Y$ characterized by
\begin{eqnarray*}
\hi K &=& 0 \text{ \ if \ } i < -\dim Y \\
\mathcal{H}^{-\dim Y} K_{\left| U\right.}&=&\xi \\
\dim \left(\supp \hi K \right) &<& -i \text{ \ if \ } i > -\dim Y \\
\dim \left(\supp \hi D_Y K \right) &<& -i \text{ \ if \ } i > -\dim Y.
\end{eqnarray*}
We also denote $\underline{\mathcal{I C}}_{Y,\xi}^\bullet$ its extension $j_*\underline{\mathcal{I C}}_{Y,\xi}^\bullet$.  The intersection complex defined by Goresky-MacPherson \cite{Goresky_MacPherson} and Deligne is obtained by shifting this perverse sheaf
\[
\ic{Y,\xi}:=\uic{Y,\xi}\left[-\dim Y\right].
\]
\end{definition}

\begin{remark}
The intersection complex does not depend on the choice of smooth open subset in $Y$. When the local system $\xi$ is not specified, it is chosen to be the constant sheaf $\kappa_U$ and $\uic X := \uic {X,\kappa_U}$.
\end{remark}

\begin{definition}[Intersection cohomology]
Let $p:X\to \Spec \KK$ the structural morphism and $k$ an integer. The $k$-th intersection cohomology space of $X$ is
\[
IH^k(X,\kappa) := \mathcal{H}^k p_* \ic X
\]
and the $k$-th compactly supported intersection cohomology space of $X$ is
\[
\IH^k_c (X,\kappa) := \mathcal{H}^k p_! \ic X.
\]

\end{definition}

For $\KK=\C$, Saito \cite{saito1986} proved that the intersection cohomology spaces carry a mixed-Hodge structure. Thus there exists on $I H_c^k(X,\mathbb{Q})$ an increasing finite filtration called the weight filtration and denoted by $W^k_{\bigcdot}$  such that the complexified quotient $\mathbb{C}\otimes_{\mathbb{Q}}W^k_{r}/ W^k_{r-1}$ carries a pure Hodge structure of weight $r$. The Hodge numbers of this structure are denoted $h_c ^{i,j,k}(X)$ and satisfy $i+j=r$.
\begin{definition}\label{def_polynomial_coho}
The mixed-Hodge structure is encoded in the mixed-Hodge polynomial:
\begin{equation}
    I H_c\left(X;x,y,v\right):=\sum_{i,j,k}  h_c^{i,j,k}(X) x^i y^j v^k.
\end{equation}
This polynomial has two important specialisations, the Poincaré polynomial
\begin{equation}
    P_c(X;t):=I H_c\left(X;1,1,v\right)=\sum_k \dim I H_c^k(X,\kappa) v^k
\end{equation}
and the $E$-polynomial
\begin{equation}
    E_c (X;x,y):=I H_c \left(X;x,y,-1\right).
\end{equation}
\end{definition}

\subsection{Symmetric functions}
In this section facts about symmetric functions are recalled. Symmetric functions is a convenient language to perform computations of cohomology of character varieties.
\subsubsection{General notations}
\begin{notations}[Partitions]\label{notations_partitions}
A partition of an integer $n\in\NN$ is a decreasing sequence of non-negative integers 
\[
\lambda=(\lambda_1,\ \lambda_2,\dots, \lambda_{l(\lambda)}) \text{ with } |\lambda|:=\lambda_1+\lambda_2+\dots+\lambda_{l(\lambda)}=n.
\]
The length of $\lambda$ is the number $l(\lambda)$ of non-zero terms. The set of partitions of $n$ is denoted by $\Part_n$ and 
\[
\Part:=\bigcup_{n\in\NN } \Part_n.
\]
The dominance ordering on $\Part$ is defined by $\lambda \preceq \mu$ if and only if
$|\lambda|=|\mu|$ and 
\[
\sum_{i=1}^{k}\lambda_i \le \sum_{i=1}^{k}\mu_i \text{ for all } k\in \NN
\]
For $\lambda=(\lambda_1,\dots,\lambda_l)$ a partition then
\[
\Part_{\lambda} := \Part_{\lambda_1}\times\dots\times \Part_{\lambda_l}.
\]
\end{notations}
\begin{notations}[Young diagrams]
The Young diagram of a partition $\lambda$ is the set
\[
\left\lbrace (i,j) \left| 1\le i \le l(\lambda) \text{ and } 1 \le j \le \lambda_i \right. \right\rbrace.
\]
A partition is often identified with its Young diagram so that $(i,j)\in \lambda$ means that $(i,j)$ belongs to the Young diagram of $\lambda$. The transpose of a Young diagram is obtained by permuting $i$ and $j$. The transpose $\lambda'$ of a partition $\lambda$ is the partition with Young diagram the transpose of the Young diagram of $\lambda$. The Young diagram of the partition $\lambda=(5,4,2)$ has the following form
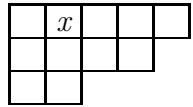
\[
\begin{Young}
      & x& &  &  \cr
      &  &  & \cr
      & \cr
\end{Young}
\]
with $x$ the box $(i,j)=(1,2)$. The arm length of $x$ is number of box right of $x$, here $a(x)=3$. The leg length is the number of box under $x$, here $l(x)=2$.
\end{notations}
\begin{notations}[Symmetric functions]
Let $X=(x_1,x_2,\dots)$ be an infinite set of variable and $\SymF[X]$ be the ring of symmetric functions in $(x_1,x_2,\dots)$. We use the usual notations from Macdonald's book \cite{mcd}. In particular for basis of symmetric functions indexed by partitions: $m_{\lambda}$, $e_{\lambda}$, $h_{\lambda}$, $p_{\lambda}$ and $s_{\lambda}$.

The Hall pairing is denoted by $\left\langle\dots,\dots\right\rangle$ and satisfy
\begin{equation}\label{orthp}
    \left\langle p_\lambda,p_\mu \right\rangle =\delta_{\lambda,\mu}z_\lambda
\end{equation}
$\delta_{\lambda,\mu}$ is $1$ if $\lambda=\mu$ and $0$ otherwise.  $z_\lambda$ is the order of a the stabilizer of a partition of cycle type $\lambda$. Namely
\begin{equation*}
    z_\lambda=\prod_{l=1}^k i_l^{m_{l}}m_{l}! 
\end{equation*}for a partition $\lambda=(\underbrace{i_1,\dots\ i_1}_{m_{1}},\dots\ ,\underbrace{i_k,\dots\ i_k}_{m_{k}})$ .
\end{notations}

\begin{definition}[Adams operator]
Adams operators $p_n$ for $n\in\mathbb{N}_{>0}$ are ring morphisms on $\SymF[X]$ defined by their values on the power sums
\[
p_m\left[p_n[X]\right]:=p_{mn}[X]\; \text{ for } m \in \NN_{>0} \text{ and } n \in \NN.
\]
The following notation is frequently used for Adams operators
\[
F\left[X^n\right]:=p_n\left[F[X]\right]
\]
\end{definition}

\subsubsection{Generating series and plethystics operations}
Fix a positive integer $k$ and consider the space of multivariate symmetric functions over $\mathbb{Q}(q,t)$:
\[
\SymF\left[X_1,\dots,X_k\right]:=\mathbb{Q}(q,t)\otimes \SymF[X_1]\otimes\dots\otimes \SymF[X_k].
\]
Cohomological information about character varieties is naturally encoded by an element of $\SymF\left[X_1,\dots,X_k\right][[s]]$, the ring of series with coefficients in $\SymF\left[X_1,\dots,X_k\right]$. Adams operators extend to ring morphisms of $\SymF\left[X_1,\dots,X_k\right][[s]]$ defined by
\[
p_n\left[f(q,t)F_1\left[X_1\right]\otimes\dots\otimes F_k\left[X_k\right]s^l\right]=f(q^n,t^n)F_1\left[X_1^n\right]\otimes\dots\otimes F_k\left[X_k^n\right]s^{n l}
\]
\begin{definition}[Plethystic exponential and logarithm]
The plethystic exponential $\Exp : s \SymF [X_1,\dots,X_k][[s]] \to \SymF [X_1,\dots,X_k][[s]]$ is defined by
\[
\Exp [G] := \exp\left( \sum_{n\ge 1} \frac{p_n[G]}{n} \right)
\]
The plethystic logarithm $\Log : 1+ s \SymF [X_1,\dots,X_k][[s]] \to \SymF [X_1,\dots,X_k][[s]]$ is defined by
\[
    \Log[1+G]:=\sum_{n\ge 1}\frac{\mu(n)}{n}p_n \left[\log(1+G)\right].
\]
with $\mu$ the usual Mobius function.
Contrarily to the ordinary ones, the plethystic exponential and logarithm are written with an uppercase character.
\end{definition}

\begin{remark}
Note the following relations for plethystics operations
\begin{eqnarray*}
    \Exp[F+G]&=&\Exp[F]\Exp[G]\\
    \Log[(1+F)(1+G)]&=&\Log[1+F]+\Log[1+G] \\
    \Log[\Exp[G]] &=&G
\end{eqnarray*}
The first relations come from Adams operators being ring morphisms, the last from the characterisation of Mobius function.
\end{remark}

\subsubsection{Symmetric functions and characters of the symmetric group}
There is a correspondence between symmetric functions and representations of symmetric groups (see \cite{mcd}). Let $R_n$ be the space of characters of of $\Sym_n$. Consider $R=\bigoplus_{n\in\mathbb{N}} R_n$, it is endowed with a pairing $\left\langle\dots,\dots\right\rangle$ such that $R_m$ is orthogonal to $R_n$ for $m\ne n$. For $V_\chi$, respectively $V_\eta$ representations of $\Sym_n$ with characters $\chi$, respectively $\eta$.
\[
\left\langle \chi,\eta \right\rangle = \dim \Hom_{\Sym_n} \left(V_{\chi}, V_{\eta}\right)
\]
$R$ is endowed with a product called external tensor product. For $V_\chi$, respectively $V_{\eta}$ representations of $\Sym_m$, respectively $\Sym_n$, the space $V_\chi \otimes V_\eta$ is a representation of $\Sym_m \times \Sym_n$. The product $\chi.\eta\in R_{m+n}$ is defined to be the character of the representation $\Ind_{\Sym_m\times\Sym_n}^{\Sym_{m+n}}V_\chi \otimes V_\eta$.

Irreducible characters of the symmetric group $\Sym_n$ are indexed by partitions: $\left(\chi_{\lambda}\right)_{\lambda\in\Part_n}$. Define the characteristic map $ch:R\to \SymF[X]$ by $ch\left( \chi_{\lambda}\right)=s_{\lambda}$. 
In particular the characteristic map sends the sign representation of $\Sym_n$ to the elementary symmetric function $e_n$.

\begin{proposition}\label{prop_characteristic_map} 
The characteristic map $\ch$ is an isomorphism between $R$ and $\SymX$ compatible with the products and the pairings ($\SymX$ being endowed the Hall pairing).
\end{proposition}
\begin{proof}
See Macdonald \cite[I-7]{mcd}.
\end{proof}

\begin{remark}\label{remark_trace_mult}
Let $\chi_V\in R_n$ the character of a representation $V$ of $\Sym_n$. The Schur functions and the power sums have the following representation theoretic interpretation:
\begin{itemize}
    \item $\left\langle s_{\lambda},\ch(\chi_{V})\right\rangle$ is the multiplicity of the irreducible representation $V_{\lambda}$ in the representation $V$.
    \item $\left\langle p_{\mu},\ch(\chi_{V})\right\rangle$ is the trace of an element in $\Sym_n$ with cycle type $\mu$ on the representation $V$.
\end{itemize}
\end{remark}

\begin{lemma}\label{lemma_ind_sign}
For $\nu$ a partition of $n$ let $\epsilon_{\nu}$ the sign representation of 
\[
\Sym_{\nu}=\Sym_{\nu_1}\times\dots\times \Sym_{\nu_l}.
\]
A choice of inclusion $\Sym_{\nu}\subset \Sym_n$ allows to induce $\epsilon_{\nu}$. Then for $\lambda\in\Part_n$
\[
\dim \Hom_{\Sym_n} \left(\Ind_{\Sym_{\nu}}^{\Sym_n} \epsilon_{\nu}, V_{\lambda}\right)=\left\langle e_{\nu},s_{\lambda}\right\rangle = \left\langle h_{\nu},s_{\lambda'}\right\rangle.
\]
\end{lemma}
\begin{proof}
$\dim \Hom_{\Sym_n} \left(\Ind_{\Sym_{\nu}}^{\Sym_n} \epsilon_{\nu}, V_{\lambda}\right)$ is the multiplicity of the irreducible representation $V_{\lambda}$ in $\Ind_{\Sym_{\nu}}^{\Sym_n} \epsilon_{\nu}$.  For $m\in \mathbb{N}_{>0}$ the symmetric function $e_{m}$ is the characteristic of the sign representation of $\Sym_m$. Therefore $e_{\nu}$ is the characteristic of $\Ind_{\Sym_{\nu}}^{\Sym_n} \epsilon_{\nu}$. First equality now follows from Remark \ref{remark_trace_mult}.
To obtain the second equality, notice that $V_{\lambda'}$ is the representation $V_{\lambda}$ twisted by the sign.
\end{proof}

\subsection{Conjugacy classes in the general linear group}\label{sect_conjugacy_classes}
\subsubsection{Notations for conjugacy classes}\label{subsect_notations_adjoint_conj}

For $r$ an integer and $ z \in \C^*$, denote by $J_r (z)$ the Jordan block of size $r$ with eigenvalue~$z$

\[
J_{r}(z):=  \begin{pmatrix}
    z \ & 1 & & & \\
    & z & \ddots & &\\
    &  & \ddots & 1 & \\
    & & & z & 1 \\
    & & & & z
  \end{pmatrix} \in \GL_{r}.
\]
Let $\mu=(\mu_1,\mu_2,\dots,\mu_s)$ a partition of an integer $m$ and let $z \in \C^*$. Denote by $J_{\mu} (z)$ the matrix with eigenvalue $z$ and Jordan blocks of size $\mu_j$.
\[
J_{\mu} (z) :=      
\begin{pmatrix}
    J_{\mu_1}(z) & &  &\\
    & J_{\mu_2}(z) & & \\
    & & \ddots & \\
     & & & J_{\mu_s}(z)
  \end{pmatrix} \in \GL_{m}.
\]
Let $\nu=(\nu_1,\dots,\nu_l) \in \Part_n$ a partition of $n$, introduce the following notation
\[
\Part_{\nu} := \Part_{\nu_1}\times\Part_{\nu_2} \times\dots\times \Part_{\nu_l}.
\]
Consider a diagonal matrix $\sigma$
\begin{equation}\label{eq_def_sigma}
\sigma = \begin{pmatrix}
    \sigma_1 \Id_{\nu_1} & & & \\
    & \sigma_2 \Id_{\nu_2} & & \\
    &  & \ddots & \\
    &  &   &  \sigma_l \Id_{\nu_l}
\end{pmatrix}
\end{equation}
with $\sigma_i \ne \sigma_j$ for $i\ne j$, so that $\nu_i$ is the multiplicity of the eigenvalue $\sigma_i\in \CC^*$. Let $\underline{\mu}= \left(\mu^{1},\dots,\mu^{l}\right)\in \Part_{\nu}$. 
\begin{notations}\label{notations_jmus}
Denote by $\mathcal{C}_{\underline{\mu},\sigma}$ the conjugacy class of the matrix:
\[
J_{\underline{\mu},\sigma}:=\begin{pmatrix}
   J_{\mu^{1}}(\sigma_1 ) & & & \\
    & J_{\mu^{2}}(\sigma_2)  & & \\
    &  & \ddots & \\
    &  &   &  J_{\mu^{l}}(\sigma_l).
\end{pmatrix}
\]

\end{notations}
We recall a well-known proposition.
\begin{proposition}
The Zariski closure of the conjugacy class $\mathcal{C}_{\underline{\mu},\sigma}$ is
\[
\overline{\mathcal{C}}_{\underline{\mu},\sigma} = \bigcup_{\underline{\rho}\preceq\underline{\mu} } \mathcal{C}_{\underline{\rho},\sigma}
\]
the union is over $l$-uple $\underline{\rho}=\left(\rho^{1},\dots,\rho^{l} \right)$ with
$
\rho^{j}\preceq \mu^{j} \; \text{ for all } 1\le j \le l
$. The dominance order on partition was recalled in \ref{notations_partitions}. \end{proposition}

\subsubsection{Resolutions of Zariski closure of conjugacy classes}\label{subsect_resolution_conjugacy}
In this section the construction of resolutions of closures of conjugacy classes is recalled. 
This construction comes from Kraft-Procesi \cite{Kraft_Procesi_80}, Nakajima \cite{Nakajima_1998,nakajima_01}, Crawley-Boevey \cite{Crawley-Boevey_parabolic, crawley-boevey2003} and Shmelkin \cite{shmelkin2009} (see also Letellier \cite{letellier_gl}).

Consider a conjugacy class $\mathcal{C}_{\underline{\mu},\sigma}$. Notations are introduced in previous section, $\sigma$ in $\GL_n$ is a diagonal matrix like in \eqref{eq_def_sigma}, denote by $M$ its centralizer in $\GL_n$.
\[
M=\begin{pmatrix}
   \GL_{\nu_1} & 0 & \\
   0 & \GL_{\nu_2} &  \\
   \vdots & 0 & \ddots 
\end{pmatrix}
\]
$\underline{\mu}=(\mu^{1},\dots,\mu^{l})$ with $\mu^{i}$ a partition of $\nu_i$. The transposed partition is denoted by ${\mu^{i}}'=\left({\mu^{i}_1}',{\mu^{i}_2}',\dots \right) $. Let $L$ the subgroup of $\GL_n$ formed by block diagonal matrices with blocks of size ${\mu^{i}_r}'$, it is a subgroup of $M$ with the following form
\[
L=\left(
\begin{array}{ccccc}
\overbrace{
\begin{array}{ccccc}
  \GL_{{\mu^{1}_1}' }  & 0 &   \\
   0  &  \GL_{{\mu^{1}_2}' } &  \\
 \vdots & 0 & \ddots 
\end{array}
}^{\nu_1}  &     &   \\
  &  \overbrace{
\begin{array}{ccccc}
  \GL_{{\mu^{2}_1}' }  & 0 &   \\
   0  &  \GL_{{\mu^{2}_2}' } &  \\
 \vdots & 0 & \ddots 
\end{array}
}^{\nu_2}  &   \\
  &    &  \ddots 
\end{array} 
\right).
\]
\begin{notations}\label{notations_group_part}
For $\nu=(\nu_1,\dots,\nu_l)$ a partition let  \[
\Sym_{\nu}= \Sym_{\nu_1}\times\dots\times \Sym_{\nu_l} \; \text{ and } \; \GL_{\nu}:=\GL_{\nu_1}\times\dots\times\GL_{\nu_l}.
\]
For $\underline{\rho}=(\rho^1,\dots,\rho^l)\in \Part_{\nu}$
\[
\GL_{\underline{\rho}}:=\GL_{\rho^1}\times\dots\GL_{\rho^l}=\prod_{r,s}\GL_{\rho^r_s}
\]
and
\[
\Sym_{\underline{\rho}} :=\Sym_{\rho^1}\times\dots\Sym_{\rho^l}=\prod_{r,s}\Sym_{\rho^r_s}
\]
Then the previously introduced Levi subgroups satisfy $M\cong \GL_{\nu}$ and $L\cong \GL_{\underline{\mu}'}$.
\end{notations}
Denote by $P$ the parabolic subgroup of blocks upper triangular matrices having $L$ as a Levi factor, $P=L U_P$ with
\[
U_P =\left(
\begin{array}{ccccc}
\overbrace{
\begin{array}{ccccc}
  \Id_{{\mu^{1}_1}' }  & * &   \\
   0  &  \Id_{{\mu^{1}_2}' } &  \\
 \vdots & 0 & \ddots 
\end{array}
}^{\nu_1}  &  \text{\Large{*}} &   \\
  &  \overbrace{
\begin{array}{ccccc}
  \Id_{{\mu^{2}_1}' }  & * &   \\
   0  &  \Id_{{\mu^{2}_2}' } &  \\
 \vdots & 0 & \ddots 
\end{array}
}^{\nu_2}  &   \\
  &    &  \ddots 
\end{array} 
\right).
\]
Now we can construct a resolution of singularities of $\overline{\mathcal{C}}_{\underline{\mu},\sigma}$
\[
\mathbb{\widetilde{X}}_{L,P,\sigma} :=\left\lbrace (X, g P) \in \GL_n \times  \GL_n / P \left| g\inv X g \in \sigma U_P \right.\right\rbrace 
\]
\begin{proposition}[Resolutions of Zariski closures of conjugacy classes]\label{prop_resol_closure}
The image of the projection to the first factor $\mathbb{\widetilde{X}}_{L,P,\sigma}\to \GL_n   $ is the Zariski closure of the conjugacy class $\mathcal{C}_{\underline{\mu},\sigma}$. Moreover the following map is a resolution of singularities
\[
\map{p^{\sigma}}{\mathbb{\widetilde{X}}_{L,P,\sigma}}{\overline{\mathcal{C}}_{\underline{\mu},\sigma}}{(X,g P)}{X}.
\]
\end{proposition}

\subsection{Resolutions of closures conjugacy classes and Weyl group actions}\label{sect_resol_cc}
\subsubsection{Borho-MacPherson approach to Springer theory}\label{subsect_springer_BM}
The approach of Borho-MacPherson \cite{borho} to Springer theory is recalled, it relies on perverse sheaves. It follows work of Lusztig \cite{Lusztig81} for the general linear group. $G$ is a reductive group over $\KK$ and $B$ a Borel subgroup of $G$. There is a decomposition $B=T U$ with $T$ a maximal torus and $U$ the unipotent radical of $B$. Consider the Grothendieck-Springer resolution
\[
\widetilde{G}=\left\lbrace (X,g B)\in G \times G/B \left| g\inv X g\in B  \right.\right\rbrace.
\]
Denote by $G^{\reg}\subset G$ the subset of regular semi-simple elements and
\[
\widetilde{G}^{\reg}=\left\lbrace (X,g B)\in G^{\reg} \times G/B \left| g\inv X g\in B  \right.\right\rbrace.
\]
Let $T^{\reg}:=G^{\reg}\cap T$, one has the following isomorphism
\[
\begin{array}{ccc}
      \left\lbrace (X,g T)\in G\times G/T \left| g\inv X g \in T^{\reg}\right.\right\rbrace& \to &\widetilde{G}^{\reg} \\
    (X,g T) & \mapsto & (X, g B) 
\end{array}
\]
The Weyl group $W=N_G(T)/T$ acts on the left hand side by 
\[
w.(X,g T) := (X, g \Dot{w}\inv T)
\]
so that it acts on $\widetilde{G}^{\reg}$.
Consider the following map
\[
\map{p^G}{\widetilde{G}}{G}{(X,g B)}{X}.
\]
Denote by $p^{\reg}$ its restriction to $\widetilde{G}^{\reg}$. Then $p^{\reg}$ is a Galois cover with group $W$ so that $W$ acts on $p^{\reg}_!\kappa_{\widetilde{G}^{\reg}}$.
Denote by $\mathcal{U}\subset G$ the subset of unipotent elements and
\[
\widetilde{\mathcal{U}}=\left\lbrace (X,g B)\in \mathcal{U} \times G/B \left| g\inv X g\in U  \right.\right\rbrace.
\]
Consider the following diagram, both squares are cartesian
\[
\begin{tikzcd}
 \widetilde{\mathcal{U}} \arrow[r,hookrightarrow] \arrow[d,"p^{\mathcal{U}}",swap] & \widetilde{G} \arrow[d,"p^G",swap] \arrow[r,hookleftarrow,"\widetilde{i}"] &  \widetilde{G}^{\reg} \arrow[d,"p^{\reg}"] \\
 \mathcal{U} \arrow[r,hookrightarrow] & G \arrow[r,hookleftarrow,"i"'] & G^{\reg}.
\end{tikzcd}
\]
\begin{proposition}[Borho-MacPherson \cite{borho}, 2.6]
There is a natural action of the Weyl group $W$ on $p^G_! \kappa_{\widetilde{G}}$ and on $p^{\mathcal{U}}_! \kappa_{\widetilde{\mathcal{U}}}$. Moreover
\[
i^* p^G_! \kappa_{\widetilde{G}} \cong p^{\reg}_!\kappa_{\widetilde{G}^{\reg}}
\]
and this isomorphism is compatible with the $W$-action.
\end{proposition}

To study characters varieties, this construction is used when $G$ is either $\GL_n$ or a Levi subgroup of a parabolic subgroup of $\GL_n$.
\begin{example}
When $G=\GL_n$, the Weyl group is isomorphic to a symmetric group $\sn$. The irreducible representations of the symmetric group $\sn$ are indexed by partitions of $n$. For $\lambda\in\Part_n$ the associated irreducible representation is $V_\lambda$. The trivial representation is $V_{(n)}$ and $V_{(1^n)}$ is the sign representation. Then there is a nice description of the left $W$-action on $p^{\mathcal{U}}_!\kappa_{\widetilde{\mathcal{U}}}$
\[
p^{\mathcal{U}}_! \kappa_{\widetilde{\mathcal{U}}}\left[\dim{ \widetilde{\mathcal{U}}}\right] =\bigoplus_{\lambda\in \Part_n}  V_{\lambda}\otimes\underline{\mathcal{I C}}^{\bigcdot}_{\overline{\mathcal{C}}_{\lambda}}.
\]
With $\mathcal{C}_{\lambda}$ the unipotent class with Jordan type $\lambda$. With notations from previous section $\mathcal{C}_{\lambda}=\mathcal{C}_{\lambda,1}$.
\end{example}

\begin{example}
Using Notations \ref{notations_group_part}, for a Levi subgroup $M$ of a parabolic subgroup of $\GL_n$ with 
\[
M\cong \GL_{\nu}
\]

the Weyl group $W_M = N_M ( T )/T$ is isomorphic to $\Sym_{\nu}$. Let $\mathcal{U}_M\subset M$ the subset of unipotent element in $M$ and $\widetilde{\mathcal{U}}_M$ its Springer resolution. The result for $\GL_n$  generalizes to
\begin{equation}\label{eq_Springer_NL}
p_!^{\mathcal{U}_M} \kappa_{\widetilde{\mathcal{U}}_M} \left[\dim \widetilde{\mathcal{U}}_M\right]=\bigoplus_{\underline{\rho}\in \Part_{\nu}} V_{\underline{\rho}}\otimes \underline{\mathcal{I C}}^{\bigcdot}_{\overline{\mathcal{C}^M_{\underline{\rho}}}}
\end{equation}
with $\mathcal{C}^M_{\underline{\rho}}$ the unipotent conjugacy class in $M$ defined for $\underline{\rho}=\left(\rho^{1},\dots,\rho^{l}\right)$ by
\[
\mathcal{C}^M_{\underline{\rho}}:=\mathcal{C}_{\rho^{1}}\times\dots\times\mathcal{C}_{\rho^{l}}\subset \GL_{\nu_1}\times\dots\times\GL_{\nu_l}.
\]
and $V_{\underline{\rho}}$ the following irreducible representation of $W_M$
\[
V_{\underline{\rho}}:=V_{\rho^1}\otimes\dots\otimes V_{\rho^l}.
\]
\end{example}

\subsubsection{Parabolic induction}
In this section Lusztig parabolic induction is recalled \cite{Lusztig1984,lusztig_character_sheaves,Lusztig1986}. Most results hold for any reductive algebraic group $G$, for our purpose we assume $G$ is either $\GL_n$ or a Levi factor of a parabolic subgroup of $\GL_n$. Let $P$ be a parabolic subgroup of $G$ with Levi decomposition $P=L U_P$.  The projection to $L$ with respect to this decomposition is $\pi_{P}: L U_P\to L$. Consider the diagram
\begin{equation}\label{parabolic_induction_diagramm}
 \begin{tikzcd}
L & \arrow[l,swap,"\rho"] V_1 \arrow[r,"\rho'"] & V_2 \arrow[r,"\rho''"] & G
 \end{tikzcd}
\end{equation}
with
\begin{eqnarray*}
V_1&=&\left\lbrace (x,g)\in G\times G \left| g\inv x g \in L U_P\right.\right\rbrace \\
 V_2 &=& \left\lbrace (x, g P)\in G\times G/P \left|  g\inv x g \in L U_P\right.\right\rbrace 
\end{eqnarray*}
\begin{eqnarray*}
\rho(x,g)&=&\pi_{P}( g\inv x g ) \\
\rho'(x,g)&=& (x,g P) \\
 \rho''(x,g P)&=&x 
\end{eqnarray*}
Parabolic induction is a functor $\Ind_{L\subset P}^{G}$ from the category of $L$-equivariant perverse sheaves on $L$ to the derived category of $G$-equivariant $\kappa$-constructible sheaves on $G$.
Take $K$ an $L$-equivariant perverse sheaf on $L$. The morphism $\rho$ is smooth with connected fibers of dimension $m=\dim G +\dim U_P$. Therefore the shifted pull-back $\rho^* K[m]$ is an $L$-equivariant perverse sheaf on $V_1$. Hence there exists a perverse sheaf $\widetilde{K}$ on $V_2$, unique up to isomorphism, such that $\rho'^* \widetilde{K}[\dim P] \cong \rho^* K[m]$. Then the parabolic induction of $K$ is defined by $\Ind_{L \subset P}^{G} K := \rho''_{!}\widetilde{K}$.

\begin{example}
The Springer complex $p^G_! \kappa_{\widetilde{G}}$ is nothing but $\Ind_{T\subset B}^G \kappa_T$ and the $W$-action on this complex is a particular case of a more general situation studied by Lusztig \cite{Lusztig1986}.
\end{example}
\begin{example}\label{example_induction_resolution}
Parabolic induction also relates to the resolutions of closures of conjugacy classes from \ref{subsect_resolution_conjugacy}.
Consider the following diagram with the first line being the diagram of parabolic induction
\begin{equation*}
     \begin{tikzcd}
L & \arrow[l,swap] V_1 \arrow[r] & V_2 \arrow[r] & \GL_n \\
\left\lbrace \sigma \right\rbrace \arrow[u,hookrightarrow] &\arrow[l]\arrow[r] \arrow[u,hookrightarrow]  \widehat{\mathbb{X}}_{L,P, \sigma} & \arrow[r,swap,"p^{\sigma}"] \arrow[u,hookrightarrow]  \mathbb{\widetilde{X}}_{L,P, \sigma} &\arrow[u,hookrightarrow] \overline{\mathcal{C}}_{\underline{\mu},\sigma}
 \end{tikzcd}
\end{equation*}
with
\[
\mathbb{\widehat{X}}_{L,P,\sigma} :=\left\lbrace (X, g ) \in \GL_n \times  \GL_n  \left| g\inv X g \in \sigma U_P \right.\right\rbrace.
\]
Then
\[
p_!^{\sigma} \kappa_{\mathbb{\widetilde{X}}_{L,P,\sigma}}\left[\dim \mathbb{\widetilde{X}}_{L,P,\sigma} \right] \cong \Ind_{L\subset P}^{\GL_n}\kappa_{\left\lbrace\sigma\right\rbrace}.
\]
with $\kappa_{\left\lbrace\sigma\right)}$ the constant sheaf with support $\left\lbrace \sigma \right\rbrace$.
\end{example}

\begin{proposition}[Lusztig \cite{lusztig_character_sheaves} I-4.2]\label{prop_trans_ind}
Let $P,Q$ be parabolic subgroups of $G$ with Levi decompositions $P=L U_P$, $Q=M U_Q$ such that $P\subset Q$ and $L\subset M$, then $P\cap M$ is a parabolic sugroup of $M$ with $L$ as a Levi subgroup. Let $K$ a $L$-equivariant perverse sheaf on $L$ such that $\Ind_{L\subset P\cap M}^{M} K$ is a perverse sheaf on $M$. Then
\begin{equation*}
    \Ind_{L\subset P}^{G} K\cong \Ind_{M\subset G}^{G}\left( \Ind_{L\subset P\cap M}^{M} K\right).
\end{equation*}
\end{proposition}

Let us detail the implication of this proposition for Springer complexes.  As in previous section, $G=\GL_n$, $B$ is a Borel subgroup of $G$ and $T$ a maximal torus in $B$.  $M$ is a Levi factor of $P$ a parabolic subgroup of $G$ containing $B$, it has the following form for some $\nu\in\Part_n$
\[
M\cong \GL_{\nu}.
\]
By transitivity of the parabolic induction from previous proposition
\begin{equation}\label{eq_trans_ind_TL}
  \Ind_{T\subset B}^G\kappa_{T}\cong  \Ind_{M\subset P}^G\Ind_{T\subset B\cap M}^M \kappa_{T}.  
\end{equation}
The left hand side is the Springer complex for $G$ so that it carries a $W$-action, this action restricts to a $W_M$-action as $W_M\subset W$. Similarly $\Ind_{T\subset B\cap M}^M \kappa_{T}$ carries a $W_M$-action as it is isomorphic to 
the Springer complex for $M$. Under the parabolic induction functor $\Ind_{M\subset P}^{G}$, this $W_M$-action on $\Ind_{T\subset B\cap M}^M \kappa_{T}$ induces a $W_M$-action on $\Ind_{M\subset P}^{G}\Ind_{T\subset B\cap M}^M \kappa_{T}$. Lusztig \cite[2.5]{Lusztig1986} proved that both $W_M$-action coincide under the isomorphism \eqref{eq_trans_ind_TL}. With Example \ref{example_induction_resolution}, this implies next theorem
\begin{theorem}\label{th_parabolic_induction_constant}
Consider the resolution of the closure of an adjoint orbit ${p^{\sigma}}:{\mathbb{\widetilde{X}}_{L,P,\sigma}}\to{\overline{\mathcal{C}}_{\underline{\mu},\sigma}}$ as in Proposition \ref{prop_resol_closure}, then
\[
p_!^{\sigma} \kappa_{\mathbb{\widetilde{X}}_{L,P,\sigma}}\left[\dim \mathbb{\widetilde{X}}_{L,P,\sigma} \right] \cong \bigoplus_{\underline{\rho}\in\Part_{\nu} } \Hom_{W_M} \left( \Ind_{W_L}^{W_M} \epsilon, V_{\underline{\rho}} \right)\otimes \uic{\overline{\mathcal{C}}_{\underline{\rho},\sigma} } 
\]
with $\epsilon$ the sign representation of $W_L$ and $V_{\underline{\rho}}:=V_{\rho^1}\otimes\dots\otimes V_{\rho^l}$ an irreducible representation of $W_M$.
\end{theorem}

\subsubsection{Relative Weyl group actions on multiplicity spaces}\label{subsect_relative_weyl_class}
An interesting feature of the multiplicity spaces $\Hom_{W_M} \left( \Ind_{W_L}^{W_M} \epsilon, V_{\underline{\rho}} \right)$ is that they carry a relative Weyl group action. Before describing this action, we recall a general result about symmetric group, see Letellier \cite[6.1, 6.2]{letellier_gl}. 
\begin{notations}
A type is a sequence $\omega=(d_1,\omega^1)\dots(d_l,\omega^l)$ with $d_j$ positive integers and $\omega^j$ partitions of any integer. The degree of $\omega$ is:
\[
|\omega|:=\sum_{i=1}^l d_i \left|\omega^i\right|
\]
The Schur function associated to a type $\omega$ is
\[
s_{\omega}:= s_{\omega^1}\left[X^{d_1}\right]\dots s_{\omega^l}\left[X^{d_l}\right].
\]
and
\begin{equation}\label{eq_r_omega}
r(\omega):= \sum_{i=1}^l (d_i-1)|\omega^i|.    
\end{equation}
\end{notations}

\begin{definition}[Twisted Littlewood-Richardson coefficients]\label{def_twisted_c}
As the usual Schur funtions $(s_{\rho})_{\rho\in \Part_n}$ form  a basis of $\SymF_n [X]$, for $\omega$ a type of degree $n$, there exist coefficients $c_{\omega}^{\rho}$ such that
\[
s_{\omega} = \sum_{\rho\in\Part_n} c^{\rho}_{\omega} s_{\rho}. 
\]
Coefficients $c_{\omega}^{\rho}$ are called the twisted Littlewood-Richardson coefficients.
\end{definition}
\begin{lemma}\label{lemma_transpose_r_omega}
Let $\omega'$ the transpose of $\omega$, i.e $\omega'=\left(d_1,{\omega^1}'\right)\dots\left(d_l,{\omega^l}'\right)$. Then
\[
c^{\rho'}_{\omega'}=(-1)^{r(\omega)} c^{\rho}_{\omega}.
\]
with $r(\omega)$ defined in \eqref{eq_r_omega}.
\end{lemma}
\begin{proof}
This follows from a computation in the ring of symmetric functions using the basis of power sums, see Letellier \cite[6.2.4]{letellier_gl}.
\end{proof}

Let us recall the interpretation of Littlewood-Richardson coefficients $c_{\omega}^{\rho}$ in terms of representations of symmetric group.  The type $\omega$ defines an irreducible representation $V_{\omega}$ of the group  $\Sym_{\omega}:=\prod_{i=1}^l \Sym_{|\omega^i|}^{d_i}$.
\[
V_{\omega}:=\otimes_{i=1}^l V_{\omega^i}^{\otimes d_i}
\]
with $V_{\omega^i}$ the representation of $\Sym_{|\omega^i|}$ indexed by the partition $\omega^i$. Denote by $f_{\omega}$ the morphism $\Sym_{\omega}\to \GL(V_{\omega})$ induced by the representation $V_{\omega}$. Introduce the relative Weyl group
\[
W_{\Sym_n} \left(\Sym_{\omega},V_{\omega}\right) = \left\lbrace \nu \in N_{\Sym_n}(\Sym_{\omega}) \left| f_{\omega}(\nu\inv \dots \nu) = f_{\omega}(\dots)\right. \right\rbrace/ \Sym_{\omega}
\]
This is the group of permutations of the blocks of $\Sym_{\omega}$ corresponding to the same representation $V_{\omega^i}$.
\begin{proposition}[Letellier \cite{letellier_gl} Proposition 6.2.5]\label{proposition_tr_relative_sym}
For $\rho\in \Part_n$ and $V_{\rho}$ the associated representation of $\Sym_n$. For $\omega$ a type, the relative Weyl group $W_{\Sym_n} \left(\Sym_{\omega},V_{\omega}\right)$ acts on
\[
\Hom_{\Sym_n} \left( \Ind_{\Sym_{\omega}}^{\Sym_n} V_{\omega}, V_{\rho} \right).
\]
Let $w\in W_{\Sym_n} \left(\Sym_{\omega},V_{\omega}\right) $ acting by cyclic permutation of the $d_i$ blocks with representation $V_{\omega^i}$ for $1\le i \le l$. Then
\[
\tr \left( w, \Hom_{\Sym_n} \left( \Ind_{\Sym_{\omega}}^{\Sym_n} V_{\omega}, V_{\rho} \right) \right) = c_{\omega}^{\rho}.
\]
\end{proposition}
\begin{remark}
Assume the type $\omega$ has the following form
\[
\omega=(\lambda_1, (1))\dots (\lambda_l, (1)) \text{ with } \lambda=(\lambda_1,\dots,\lambda_l)\in\Part_n.
\]
Then $s_{\omega}=p_{\lambda}$ and  for $\rho\in\Part_n$
\[
c^{\rho}_{\omega}= \chi^{\rho}_{\lambda}.
\]
Notice that $W_{\Sym_n}\left(\Sym_n, V_{\omega}\right)\cong \Sym_n$ and the element $w$ associated to $\omega$ has cycle type $\lambda$. Therefore the proposition implies that as a $W_{\Sym_n}\left(\Sym_n, V_{\omega}\right)$ representation 
\[
\Hom_{\Sym_n} \left( \Ind_{\Sym_{\omega}}^{\Sym_n} V_{\omega}, V_{\rho} \right)\cong V_{\rho}.
\]
\end{remark}

With this general result about symmetric group, we go back to the Weyl groups relative to resolution of conjugacy classes.

\begin{definition}[Relative Weyl group]
For $L$ a Levi subgroup of $M$, The relative Weyl group is
\[
W_{M}(L) := N_{M} (L)/L.
\]
\end{definition}

Take $L$ and $M$ similarly to Section \ref{subsect_resolution_conjugacy}. Denote by $(m^{i}_1,\dots,m^{i}_{k_{i}})$ the multiplicity of the parts of ${\mu^{i}}'$ so that it has the following form
\[
{\mu^{i}}'=\left( \underbrace{a^{i}_1,\dots,a^{i}_1}_{m^{i}_1},\underbrace{a^{i}_2,\dots,a^{i}_2}_{m^{i}_2} ,\dots, \underbrace{a^{i}_{k_i},\dots,a^{i}_{k_i}}_{m^{i}_{k_i}} \right).
\]
Then with notations \ref{notations_group_part}  $L\cong \GL_{\underline{\mu}'}$ and the relative Weyl group is
\[
W_{M}(L) \cong \prod_{ \substack{  1 \le i \le l \\ 1 \le r \le k_{i}  } } \Sym_{m^{i}_{r}}.
\]
When $M=\GL_n$ then the relative Weyl group is the group of permutations of same-sized blocks of $L$.

\begin{notations}\label{notations_h_tilde}
Conjugacy classes in $W_M(L)$ are indexed by elements 
\[
\eta=(\eta^{i,r})_{  \substack{  1 \le i \le l\\ 1 \le r \le k_{i}  } }  \in \prod_{ \substack{  1 \le i \le l \\ 1 \le r \le k_{i}  } } \Part_{m^{i}_{r}}.
\]
A conjugacy class then determined $l$ types
$\omega_{\eta^i}$ with parts $\left(\eta^{i,r}_s,(1^{a^i_r})\right)_{  \substack{  1\le r \le k_i \\ 1\le s \le l(\eta^{i,r})} } $.

Note that
\[
 s_{\omega_{\eta^i}'}= \prod_{r=1}^{k_i} \prod_{s=1}^{l(\eta^{i,r})} h_{a^i_r}\left[X^{\eta^{i,r}_s}\right]
\]
Following notations will be convenient to compute Weyl group actions on the cohomology of character varieties.
\[
\widetilde{h}_{\eta}:=\prod_{i=1}^l s_{\omega_{\eta^i}'}
\]
and
\[
r(\eta):=\sum_{i=1}^l r(\omega_{\eta^i}).
\]
with $r(\omega_{\eta^i})$ defined by \eqref{eq_r_omega}.
\end{notations}

Those data describe the $W_M(L)$ action on the multiplicity spaces, Proposition \ref{proposition_tr_relative_sym} implies:
\begin{theorem}
Let $\epsilon_{\underline{\mu}'}$ the sign representation of $W_L$ and $\underline{\rho}\in\Part_{\nu}$. The relative Weyl group $W_M (L)$ acts on $\Hom_{W_M} \left( \Ind_{W_L}^{W_M} \epsilon_{\underline{\mu}'}, V_{\underline{\rho}} \right)$. The trace of the action of an element with conjugacy class indexed by $\eta \in \prod_{ \substack{  1 \le i \le l_j \\ 1 \le r \le k_{i}  } } \Part_{m^{i}_{r}}$ is
\[
\tr\left(\eta , \Hom_{W_M} \left( \Ind_{W_L}^{W_M} \epsilon_{\underline{\mu}'}, V_{\underline{\rho}} \right) \right) = \prod_{i=1}^{l} c^{\rho^i}_{w_{\eta^i}}.
\]
\end{theorem}

\section{Background on character varieties for punctured Riemann surfaces}\label{sect_background_char}
\subsection{Character varieties and their resolutions}

\subsubsection{Construction of character varieties}
Let $\Sigma$ be a compact Riemann surface of genus $g$. Consider the punctured Riemann surface $\Sigma^0=\Sigma\setminus\left\lbrace p_1,\dots,p_k\right\rbrace$ where $p_j$ are distinct points on $\Sigma$ called punctures. The field $\mathbb{K}$ is either $\C$ or an algebraic closure $\Fqbar$ of a finite field $\Fq$ with $q$ elements. Fix a non negative integer $n$. We are concerned by $n$-dimensional $\mathbb{K}$-representations of the fundamental group of $\Sigma^0$ with prescribed monodromy around the punctures. 

For each puncture, specify a conjugacy class $\mathcal{C}_{\underline{\mu}^j,\sigma^j}$. The notations are the same as in previous section, with the addition of an upper index $1\le j\le k$ labelling the punctures. $\sigma^j$ is a diagonal matrix with diagonal coefficients
\[
(\underbrace{\sigma^j_1,\dots,\sigma^j_{1}}_{\nu^j_1},\dots, \underbrace{\sigma^j_{l_j},\dots,\sigma^j_{l_j}}_{\nu^j_{l_j}})
\]
and $\sigma^j_r \ne \sigma^j_s$ for $r\ne s$.
Moreover, $\underline{\mu}^j=\left(\mu^{j,1},\dots,\mu^{j,l_j}\right)$ with $\mu^{j,r}\in \Part_{\nu^j_r}$ the partition giving the size of the Jordan blocks of the eigenvalue $\sigma^j_r$.

A bold symbol is used to represent $k$-uples:
\begin{equation}\label{eq_bold_kuple}
    \begin{array}{ccc}
\boldsymbol{\mu}& := & \left( \underline{\mu}^1,\dots, \underline{\mu}^k\right)  \\
\boldsymbol{\sigma} & := & \left(\sigma^1,\dots,\sigma^k\right) \\
\mathcal{C}_{\boldsymbol{\mu},\boldsymbol{\sigma}} & := & \left(\mathcal{C}_{\underline{\mu}^1,\sigma^1}, \dots, \mathcal{C}_{\underline{\mu}^k,\sigma^k} \right)
    \end{array}
\end{equation}

The representations of the fundamental group of $\Sigma^0$ with monodromy around $p_j$ in the closure  $\overline{\mathcal{C}}_{\underline{\mu}^j,\sigma^j}$ form the following affine variety
\begin{multline*}
\mathcal{R}_{\overline{ \mathcal{C}  }_{ \boldsymbol{\mu} ,\boldsymbol{\sigma}  } }:=\left\lbrace\left(A_1,B_1,\dots,A_g,B_g,X_1,\dots,X_k\right)\in \GL_n^{2g}\times \overline{\mathcal{C}}_{\underline{\mu}^1,\sigma^1}\times\dots\times \overline{\mathcal{C}}_{\underline{\mu}^k,\sigma^k}\right| \\ \left. A_1 B_1 A_1\inv B_1\inv \dots A_g B_g A_g\inv B_g\inv X_1\dots X_k =\Id \right\rbrace.
\end{multline*}
The group $\GL_n$ acts by simultaneous conjugation on $\mathcal{R}_{\overline{ \mathcal{C}  }_{ \boldsymbol{\mu} ,\boldsymbol{\sigma}  } }$
\begin{equation*}
    g.\left(A_1,\dots,B_g,X_1,\dots,X_k\right)=\left(g A_1 g\inv,\dots ,g B_g g\inv,g X_1 g\inv ,\dots,g X_k g\inv\right).
\end{equation*}
The center of $\GL_n$ acts trivially so this action factors through an action of $\PGL_n$.
\begin{definition}[Character variety]
The character variety we are interested in is the following GIT quotient
\begin{equation*}
    \mathcal{M}_{\overline{ \mathcal{C}  }_{ \boldsymbol{\mu} ,\boldsymbol{\sigma}  } } :=\mathcal{R}_{\overline{ \mathcal{C}  }_{ \boldsymbol{\mu} ,\boldsymbol{\sigma}  } } //\PGL_n :=\Spec \KK\left[  \mathcal{R}_{\overline{ \mathcal{C}  }_{ \boldsymbol{\mu} ,\boldsymbol{\sigma}  } }  \right]^{\PGL_n}.
\end{equation*}
It is an affine variety with regular functions the $\PGL_n$-invariants functions on $\mathcal{R}_{\overline{ \mathcal{C}  }_{ \boldsymbol{\mu} ,\boldsymbol{\sigma}  } }$.
\end{definition}
Under some generecity assumptions, the $\PGL_n$ action is free.
\begin{definition}[Generic conjugacy classes]\label{def_generic_class}
Denote $\Delta(\sigma^j)$ the multiset of eigenvalues of $\sigma^j$ repeated according to multiplicities. $\sigma^j_r$ appears exactly $\nu^j_r$ times in the multiset $\Delta(\sigma^j)$. The $k$-uple of conjugacy classes $\mathcal{C}_{\boldsymbol{\mu},\boldsymbol{\sigma}}$ is generic if and only if it satisfy the two following conditions
\begin{enumerate}
    \item \[ \prod_{j=1}^k\prod_{\alpha \in \Delta (\sigma^j )}\alpha = 1 \]
    \item For any $r \le n-1$, for all $(R_1,\dots, R_k)$ with $R_j\subset \Delta(\sigma^j)$ of size $r$
    \[
    \prod_{j=1}^k \prod_{\alpha\in R_j} \alpha \ne 1
    \]
\end{enumerate}
Throughout the thesis, every character varieties considered are assumed to have generic conjugacy classes at the punctures.
\end{definition}

\begin{definition}
Let $\mathcal{R}_{ \mathcal{C}_{\boldsymbol{\mu},\boldsymbol{\sigma}} }:=\mathcal{R}_{ \overline{\mathcal{C}}_{\boldsymbol{\mu},\boldsymbol{\sigma}} } \cap \left( \GL_n(\KK)^{2 g}\times \prod_{j=1}^k \mathcal{C}_{\underline{\mu}^j,\sigma^j} \right)$ and $\mathcal{M}_{ \mathcal{C}_{\boldsymbol{\mu},\boldsymbol{\sigma}} } $ the image of $\mathcal{R}_{ \mathcal{C}_{\boldsymbol{\mu},\boldsymbol{\sigma}} }$ in $\mathcal{R}_{ \overline{\mathcal{C}}_{\boldsymbol{\mu},\boldsymbol{\sigma}} }$.
\end{definition}

We recall a proposition from \cite{zclosure}, and \cite{aha} for the semisimple case.

\begin{proposition}\label{prop_smooth}
If $\mathcal{C}_{\boldsymbol{\mu},\boldsymbol{\sigma}}$ is generic then $\mathcal{R}_{\mathcal{C}_{\boldsymbol{\mu},\boldsymbol{\sigma}}}$ is non-singular, when non-empty its dimension is
\[
\dim \mathcal{R}_{\mathcal{C}_{\boldsymbol{\mu},\boldsymbol{\sigma}}} = 2 g n^2 -n^2 +1 +\sum_{j=1}^k \dim \mathcal{C}_{\underline{\mu}^j,\sigma^j}.
\]
\end{proposition}
\begin{proof}
The proof combines the one of theorem 2.2.5 in \cite{hausel_villegas} and proposition 5.2.8 in \cite{etingof}.
\end{proof}

\begin{proposition}[Stratification of $\mathcal{M}_{\overline{\mathcal{C}}_{\boldsymbol{\mu},\boldsymbol{\sigma}}}$, \cite{zclosure} Corollary 3.6 ]
We assume $\mathcal{C}_{\boldsymbol{\mu},\boldsymbol{\sigma}}$ is generic. The stratification of Zariski closure of conjugacy classes induces a stratification of the character variety:
\[
\mathcal{M}_{\overline{\mathcal{C}}_{\boldsymbol{\mu},\boldsymbol{\sigma}}}=\bigsqcup_{\boldsymbol{\rho}\preceq\boldsymbol{\mu} } \mathcal{M}_{\overline{\mathcal{C}}_{\boldsymbol{\rho},\boldsymbol{\sigma}}}.
\]
The union is over $\boldsymbol{\rho} = \left(\underline{\rho}^{1},\dots,\underline{\rho}^{k}  \right)$ with $\underline{\rho}^j=\left(\rho^{j,1},\dots,\rho^{j,l_j}\right)$ such that
\[
\rho^{j,i}\preceq \mu^{j,i},\; \text{ for all }1\le j \le k,\; 1\le i \le l_j
\]
with $\preceq$ the dominance order on $\Part_{\nu^j_i}$.

Moreover if $\mathcal{M}_{\overline{\mathcal{C}}_{\boldsymbol{\mu},\boldsymbol{\sigma}}}$ is non empty, then $\mathcal{M}_{\mathcal{C}_{\boldsymbol{\mu},\boldsymbol{\sigma}}}$ is also non empty. Therefore when $\mathcal{M}_{\overline{\mathcal{C}}_{\boldsymbol{\mu},\boldsymbol{\sigma}}}$ is non empty, its dimension is 
\begin{equation}\label{eq_dim_MC}
\dim \mathcal{M}_{\overline{\mathcal{C}}_{\boldsymbol{\mu},\boldsymbol{\sigma}}} = d_{\boldsymbol{\mu}}:=n^2 (2 g-2) +2 +\sum_{j=1}^k \dim \mathcal{C}_{\underline{\mu}^j,\sigma^j }.
\end{equation}
\end{proposition}

\subsubsection{Resolutions of character varieties}\label{subsect_resolution_character}

The resolutions of conjugacy classes introduced in \ref{subsect_resolution_conjugacy} induce resolutions of character varieties. As before we consider a generic $k$-uple of conjugacy classes
\[
\mathcal{C}_{\boldsymbol{\mu},\boldsymbol{\sigma}} = \left(\mathcal{C}_{\underline{\mu}^1,\sigma^1}, \dots, \mathcal{C}_{\underline{\mu}^k,\sigma^k} \right)
\]
and upper indices $1\le j \le k$ label the puncture. $\sigma^j$ is a diagonal matrix with diagonal coefficients
\[
(\underbrace{\sigma^j_1,\dots,\sigma^j_{1}}_{\nu^j_1},\dots, \underbrace{\sigma^j_{l_j},\dots,\sigma^j_{l_j}}_{\nu^j_{l_j}}).
\]
Let $M^j:= Z_{\GL_n} (\sigma^j)$ then with Notation \ref{notations_group_part}
\[
M^j\cong \GL_{\nu^j}
\]
$\mu^{j,i}\in \Part_{\nu^{j}_i}$ gives the size of the Jordan blocks of $\mathcal{C}_{\underline{\mu}^j,\sigma^j}$ relative to the eigenvalue $\sigma^j_i$. Denote by ${\mu^{j,i}}'=\left({\mu^{j,i}_1}',{\mu^{j,i}_2}',\dots \right)$ the transposed partition. Let $L^j\subset M^j$ the subgroup of block-diagonal matrices as in \ref{subsect_resolution_conjugacy}
\[
L^j \cong \underbrace{\GL_{{\mu^{j,1}_1}'}\times \GL_{{\mu^{j,1}_2}'} \times \dots }_{\subset \GL_{\nu^j_1}}\times \dots \times \underbrace{\GL_{{\mu^{j,l_j}_1}'}\times \GL_{{\mu^{j,l_j}_2}'} \times \dots}_{\subset \GL_{\nu^j_{l_j}}}.
\]
Let $\mathbb{\widetilde{X}}_{L^j,P^j,\sigma^j}$ a resolution of $\overline{\mathcal{C}}_{\underline{\mu}^j,\sigma^j}$ as constructed in \ref{subsect_resolution_conjugacy} 
\[
\mathbb{\widetilde{X}}_{\boldsymbol{L},\boldsymbol{P},\boldsymbol{\sigma}} := \prod_{1\le j\le k} \mathbb{\widetilde{X}}_{L^j,P^j,\sigma^j}.
\]
Letellier \cite{zclosure} constructed resolutions of singularities for character varieties.
\begin{definition}[Resolutions of character varieties]\label{def_resol_char}
Define
\begin{multline}
   \mlps := \left\lbrace(A_i, B_i)_{1\le i \le g}, (X_j, g_j P^j)_{1\le j \le k} \in \GL_n^{2 g}\times \mathbb{\widetilde{X}}_{\boldsymbol{L},\boldsymbol{P},\boldsymbol{\sigma}} \right. \\
   \left| A_1 B_1 A\inv_1 B\inv_1 \dots B\inv_g X_1\dots X_k = \Id \right\rbrace//\PGL_n.
\end{multline}
The maps $p^{\sigma^j} : \mathbb{\widetilde{X}}_{L^j,P^j,\sigma^j} \to \overline{\mathcal{C}}_{\underline{\mu}^j,\sigma^j}$ induce a map
\[
p^{\boldsymbol{\sigma}} :    \mlps \to \mathcal{M}_{\overline{\mathcal{C}}_{\boldsymbol{\mu},\boldsymbol{\sigma} }  }
\]
this map is a resolution of singularity.
\end{definition}

The description of the resolutions of closures of conjugacy classes from Theorem \ref{th_parabolic_induction_constant} extends to the resolution of character varieties.
\begin{notations}
Notations for resolution of closure of conjugacy classes are extended to $k$-uple.
The Weyl group of $M^j$ is $W_{M^j}=N_{M^j}(T)/T$ then
\[
W_{M^j}\cong \Sym_{\nu^j}
\]
its irreducible representations are labelled by $\underline{\rho}^j=\left(\rho^{j,1},\dots,\rho^{j,l_j}\right)\in \Part_{\nu^j}$.

The Weyl group of $L^j$ is $W_{L^j}=N_{L^j}(T)/T$, it is a subgroup of $W_{M^j}$
\[
W_{L^j} \cong \underbrace{\Sym_{{\mu^{j,1}_1}'}\times \Sym_{{\mu^{j,1}_2}'} \times \dots }_{\subset \Sym_{\nu^j_1}}\times \dots \times \underbrace{\Sym_{{\mu^{j,l_j}_1}'}\times \Sym_{{\mu^{j,l_j}_2}'} \times \dots}_{\subset \Sym_{\nu^j_{l_j}}}.
\]
The sign representation for this Weyl group is denoted by $\epsilon_{{\underline{\mu}^j}'}$ to remind the form of the Weyl group $W_{L^j}\cong \Sym_{{\underline{\mu}^j}'}$.

Define $W_{\boldsymbol{L}} := \prod_{j=1}^k W_{L^j} $ and similarly $W_{\boldsymbol{M}}:= \prod_{j=1}^k W_{M^j}$.
The parameter $\boldsymbol{\rho}=(\underline{\rho}^1,\dots,\underline{\rho}^k)\in \Part_{\nu^1}\times \dots \times \Part_{\nu^k}$ indexes irreducible representations of $W_{\boldsymbol{M}}$
\[
V_{\boldsymbol{\rho}}=\bigotimes_{j=1}^k V_{\underline{\rho}^j}.
\]
Let $\epsilon_{\boldsymbol{\mu'}}$ be the sign representation of $W_{\boldsymbol{L}}$.
\end{notations}

Next theorem is a particular case of a result of Letellier \cite[Theorem 5.4]{zclosure}
\begin{theorem}\label{th_coho_resol_char}
\[
p^{\boldsymbol{\sigma}}_! \kappa [d_{\boldsymbol{\mu}}] \cong \bigoplus_{\boldsymbol{\rho\preceq \mu} }  A_{\boldsymbol{\mu'},\boldsymbol{\rho} }\otimes \uic{\mathcal{M}_{\overline{\mathcal{C}}_{\boldsymbol{\rho},\boldsymbol{\sigma}}  }  }
\]
and in terms of cohomology:
\begin{equation}\label{eq_cohomology_resol}
H^{i+d_{\boldsymbol{\mu}}}_c
\left(\mlps,\qlbar \right) 
\cong 
\bigoplus_{\boldsymbol{\rho\preceq \mu}}
A_{\boldsymbol{\mu'},\boldsymbol{\rho} }\otimes
I H_c^{i+d_{\boldsymbol{\rho}} } \left(\mathcal{M}_{\overline{\mathcal{C}}_{\boldsymbol{\rho},\boldsymbol{\sigma}}  } ,\qlbar \right).
\end{equation}
The multiplicity space is given by
\[
A_{\boldsymbol{\mu'},\boldsymbol{\rho} }:= \Hom_{W_{\boldsymbol{M}}}\left(\Ind_{W_{\boldsymbol{L}}}^{W_{\boldsymbol{M}}}\epsilon_{\boldsymbol{\mu'}}, V_{\boldsymbol{\rho}} \right) \cong \bigotimes_{j=1}^k \Hom_{W_{M^j}}\left(\Ind_{W_{L^j}}^{W_{M^j}} \epsilon_{{\underline{\mu}^j}'} , V_{\underline{\rho}^j}\right).
\]
\end{theorem}

\subsubsection{Relative Weyl group actions}\label{subsect_relative_weyl_char}

The relative Weyl group actions on the cohomology of resolutions of closure of conjugacy classes give rise to relative Weyl group actions on the cohomology of resolutions of character varieties.

\begin{notations}\label{notations_k_h_tilde}
The relative Weyl group is
\[
W_{\boldsymbol{M}}\left(\boldsymbol{L}\right):=\prod_{j=1}^k W_{M^j}(L^j)
\]
with $W_{M^j}(L^j)$ the relative Weyl groups described in \ref{subsect_relative_weyl_class}.
Conjugacy classes in $W_{\boldsymbol{M}}(\boldsymbol{L})$ are labelled by elements
\[
\boldsymbol{\eta} = (\eta^j)_{1\le j\le k}
\]
with $\eta^j $ 
as in \ref{subsect_relative_weyl_class} with an additional index $j$ for the puncture:
\[
\eta^j =(\eta^{j,i,r})_{  \substack{  1 \le i \le l_j\\ 1 \le r \le k_{j,i}  } }  \in \prod_{ \substack{  1 \le i \le l_j \\ 1 \le r \le k_{j,i}  } } \Part_{m^{j,i}_{r}}.
\]
Notations \ref{notations_h_tilde} extend to $k$-uple:
\[
\widetilde{h}_{\boldsymbol{\eta}}:=\prod_{j=1}^k\prod_{i=1}^{l_j} s_{ {\omega_{\eta^{j,i}}}' }\left[ X_j\right]
\]
and
\[
r(\boldsymbol{\eta}):=\sum_{j=1}^k\sum_{i=1}^{l_j} r(\omega_{\eta^{j,i}} ).
\]
\end{notations}

\begin{theorem}\label{th_tr_resol}
Let $\mathcal{C}_{\boldsymbol{\mu},\boldsymbol{\sigma}}$ a generic $k$-uple of conjugacy classes and $\mlps$ the resolution of $\mathcal{M}_{\overline{\mathcal{C}}_{\boldsymbol{\mu},\boldsymbol{\sigma}}}$. The relative Weyl group $W_{\boldsymbol{M}}(\boldsymbol{L})$ acts on the cohomology of $\mlps$. The trace of an element in the conjugacy class indexed by $\boldsymbol{\eta}$ is 
\[
\tr\left(\boldsymbol{\eta}, H_c^{i+d_{\boldsymbol{\mu}}}\left(\mlps ,\kappa \right)\right) = \sum_{\boldsymbol{\rho\preceq \mu}} \tr\left(\boldsymbol{\eta}, A_{\boldsymbol{\mu'},\boldsymbol{\rho} }\right) H_c^{i+d_{\boldsymbol{\rho}}}\left( \mathcal{M}_{\overline{\boldsymbol{\mathcal{C}}}_{\boldsymbol{\sigma},\boldsymbol{\rho}}  } ,\kappa \right).
\]
with
\[
\tr\left(\boldsymbol{\eta}, A_{\boldsymbol{\mu'},\boldsymbol{\rho} }\right) = \prod_{j=1}^k\prod_{i=1}^{l_j} c^{\rho^{j,i}}_{\eta^{j,i}}.
\]
\end{theorem}

\subsection{Cohomology of character varieties: some results and conjectures}\label{sect_cohomology_char}
\subsubsection{Conjectural formula for the mixed-Hodge polynomial}
Hausel, Letellier and Rodriguez-Villegas \cite{aha} introduced a generating function conjecturally encoding mixed-Hodge structure on the cohomology of character varieties. As before $g$ is  a non-negative integer, the genus, and $k$ a positive integer, the number of punctures.
\begin{definition}[Generating function $\Omega$ and Hausel-Letellier-Villegas kernel]\label{def_omega_hlv}
The $k$-points, genus $g$ Cauchy function is defined by
\begin{equation}
    \Omega^g_k(z,w):=\sum_{\lambda\in\Part}\mathcal{H}_\lambda(z,w)\prod_{i=1}^k\Ht_\lambda\left[X_i,z^2,w^2\right]s^{|\lambda|}
\end{equation}
with
\begin{equation}
\mathcal{H}_\lambda(z,w):=\prod\frac{\left(z^{2a+1}-w^{2l+1}\right)^{2g}}{\left(z^{2a+2}-w^{2l}\right)\left(z^{2a}-w^{2l+2}\right)}.
\end{equation}
The degree $n$ Hausel-Letellier-Villegas kernel is defined by
\[
\mathbb{H}^{H L V}_n(z,w):=\left.(z^2-1)(1-w^2) \Log \Omega_k^g(z,w)\right|_{s^n}.
\]
\end{definition}
The generating function $\Omega_k^g(z,w)$ belongs to the lambda ring $\SymF\left[X_1,\dots,X_k\right]\left[\left[s\right]\right]$. This Cauchy function is known to encode cohomological information about character varieties and quiver varieties, let us recall these various conjectures and theorems.

When the conjugacy classes are semisimple Hausel, Letellier, Rodriguez-Villegas stated a conjecture for the mixed-Hodge polynomial of the character variety \cite{aha}. They proved the specialisation corresponding to the $E$-polynomial. Letellier generalized this conjecture to arbitrary types and intersection cohomology.

Let $\mathcal{C}_{\boldsymbol{\mu},\boldsymbol{\sigma}}$ a $k$-uple of generic conjugacy classes. Then $\boldsymbol{\mu}=(\underline{\mu}^1,\dots,\underline{\mu}^k)$ with $\underline{\mu}^j=\left(\mu^{j,1},\dots,\mu^{j,l_j}\right)$. The transposition of the partition $\mu^{j,i}\in \Part_{\nu^j_i}$ is denoted by ${\mu^{j,i}}'$ and
\begin{equation}\label{eq_s_boldmu}
s_{\boldsymbol{\mu'}}:=\prod_{j=1}^k \prod_{i=1}^{l_j} s_{ {\mu^{j,i}}' }[X_j]
\end{equation}
\begin{conjecture}[Letellier \cite{zclosure}, Conjecture 1.5]\label{conjecture_ic_char}
For $\mathcal{C}_{\boldsymbol{\mu},\boldsymbol{\sigma}}$ a generic $k$-uple of conjugacy classes, the mixed-Hodge polynomial of the character variety $\mathcal{M}_{\overline{\mathcal{C}}_{\boldsymbol{\mu},\boldsymbol{\sigma}}}$ is
\[
IH_c (  \mathcal{M}_{\overline{\mathcal{C}}_{\boldsymbol{\mu},\boldsymbol{\sigma}}}, q, v) = \left( v \sqrt{q} \right)^{d_{\boldsymbol{\mu}}} \left\langle s_{\boldsymbol{\mu'}}, \mathbb{H}_n^{H L V}\left(\frac{-1}{\sqrt{q}},v\sqrt{q}\right) \right\rangle
\]
with $q=x y$. In particular after specializing to the Poincaré polynomial
\begin{equation}\label{eq_poincare_poly}
    P_c (  \mathcal{M}_{\overline{\mathcal{C}}_{\boldsymbol{\mu},\boldsymbol{\sigma}}} ,  v) =  v^{d_{\boldsymbol{\mu}}} \left\langle s_{\boldsymbol{\mu'}}, \mathbb{H}_n^{H L V}\left(-1,v \right) \right\rangle.
\end{equation}
\end{conjecture}

Some specializations of this conjecture are already proved. The formula obtained after specialization to the $E$-polynomial is proved by Hausel, Letellier and Rodriguez-Villegas \cite{aha} for semisimple conjugacy classes  and by Letellier \cite{zclosure} for any type of conjugacy classes. The proof relies on counting points of character varieties over finite fields and representation theory of $\GL_n(\mathbb{F}_q)$. The formula obtained after specialization to the Poincaré polynomial is proved by Schiffmann \cite{Schiffmann16} for one central conjugacy class and by Mellit \cite{mellit_p} for any $k$-uple of semisimple conjugacy classes. The proof relies on counting point of moduli space of stable parabolic Higgs bundles over finite field. \newline

\subsubsection{Poincaré polynomial of character varieties with semisimple conjugacy classes at punctures}\label{subsect_poincare_ss}

Let $\boldsymbol{\mathcal{S}}=(\mathcal{S}_1,\dots,\mathcal{S}_k)$ a generic $k$-uple of semisimple conjugacy classes. Then $\mathcal{S}_j$ has the form $\mathcal{C}_{\underline{\mu}^j,\sigma^j}$ with 
$\underline{\mu^j}=(1^{\nu^j_1},\dots,1^{\nu^j_{l_j}})$ and
\[
s_{\boldsymbol{\mu'}} = \prod_{j=1}^k\prod_{i=1}^{l_j} s_{(\nu^j_i)}[X_j]  = \prod_{j=1}^k h_{\nu^j}[X_j] = h_{\boldsymbol{\nu}}.
\]

The conjecture from Hausel, Letellier, Rodriguez-Villegas \cite{aha} for the mixed-Hodge structure of the character varieties with monodromies specified by $\boldsymbol{\mathcal{S}}$ reads
\[
I H_c \left(\mathcal{M}_{\boldsymbol{\mathcal{S}}};q ,v \right) = (v \sqrt{q})^{d_{\boldsymbol{\mathcal{S}}}}\left\langle h_{\boldsymbol{\nu}},  \mathbb{H}_n^{H L V}\left(\frac{-1}{\sqrt{q}},v\sqrt{q} \right)\right\rangle.
\]
Note that as the conjugacy classes are generic semisimple, the character variety is smooth and the intersection cohomology coincides with the usual cohomology. Then the specialization to Poincaré polynomial of the conjecture is
\begin{equation}\label{eq_poincare_ss}
    P_c\left(\mathcal{M}_{\boldsymbol{\mathcal{S}}}, v  \right) = \sum_{i} v^i \dim H_c^i(\mathcal{M}_{\boldsymbol{\mathcal{S}}},\kappa) = v^{d_{\boldsymbol{\mathcal{S}}}} \left\langle h_{\boldsymbol{\nu}}, \mathbb{H}^{H L V}_n\left(-1,v\right)\right\rangle.
\end{equation}

After a change of variable $v=\frac{-1}{\sqrt{u}}$ and applying Poincaré duality, this formula is equivalent to Mellit's result  \cite[Theorem 7.12]{mellit_p} and we have the following theorem.
\begin{theorem}
For $\mathcal{\boldsymbol{S}}=\left(\mathcal{S}_1,\dots,\mathcal{S}_k\right)$ a generic $k$-uple of semisimple conjugacy classes. If the multiplicities of the eigenvalues of $\mathcal{S}_j$ are given by a partition $\nu^j\in \Part_n$ for $1\le j \le k$. Then the Poincaré polynomial of the character variety $\mathcal{M}_{\boldsymbol{\mathcal{S}}}$ is
\begin{equation}\label{eq_poincare_mellit}
    P_c \left( \mathcal{M}_{\boldsymbol{\mathcal{S}}}; v \right) = v^{d_{\boldsymbol{\mathcal{S}}}} \left\langle h_{\boldsymbol{\nu}}, \mathbb{H}_n^{H L V}(-1,v) \right\rangle.
\end{equation}
\end{theorem}

\subsubsection{Weyl group actions on the cohomology}
In \ref{subsect_relative_weyl_char} a Weyl group action on the cohomology of resolutions of character varieties was introduced. The conjecture about the mixed-Hodge structure also concerns this Weyl group action. We present the implications in terms of Poincaré polynomial using Notations \ref{notations_h_tilde} and \ref{notations_k_h_tilde}.
\begin{definition}[$\boldsymbol{\eta}$-twisted Poincaré polynomial]\label{def_twisted_Poincare}
$\mathcal{C}_{\boldsymbol{\mu},\boldsymbol{\sigma}}$ is a generic $k$-uple of conjugacy classes and $\mlps$ is the resolution of $\mcmus$. For $\boldsymbol{\eta}$ indexing a conjugacy class in $W_{\boldsymbol{M}}(\boldsymbol{L})$, the $\boldsymbol{\eta}$-twisted Poincaré polynomial of $\mlps$ is
\[
P_c^{\boldsymbol{\eta}}\left(\mlps,v\right) := \sum_{i} \tr\left( \boldsymbol{\eta}, H_c^i\left( \mlps, \kappa\right)  \right)v^i.
\]
\end{definition}

Letellier proved that the Weyl group action on the cohomology of the resolution $\mlps$ preserves the weight filtration. Therefore similarly to the $\boldsymbol{\eta}$-twisted Poincaré polynomial one can define the $\boldsymbol{\eta}$-twisted mixed-Hodge polynomial $IH_c^{\boldsymbol{\eta}}\left(\mlps,q,v\right)$. 
\begin{conjecture}[Letellier \cite{zclosure} Conjecture 1.8]\label{conj_twisted_mixed_resol}
Let $\cmus$ a generic $k$-uple of conjugacy classes. For $\mlps$ the resolution of a character variety $\mcmus$ and $\boldsymbol{\eta}$ a conjugacy class in $W_{\boldsymbol{M}}(\boldsymbol{L})$, the $\boldsymbol{\eta}$-twisted Poincaré polynomial is
\[
IH_c^{\boldsymbol{\eta}}\left(\mlps,q,v\right) = (-1)^{r(\boldsymbol{\eta})}(v\sqrt{q})^{d_{\boldsymbol{\mu}}} \left\langle \widetilde{h}_{\boldsymbol{\eta}}, \mathbb{H}^{H L V}_n \left( \frac{-1}{\sqrt{q}},v \sqrt{q}\right) \right\rangle. 
\]

\end{conjecture}

\section{Diffeomorphism between a resolution of a character variety and a semisimple character variety}\label{sect_diffeo_res}
In this section the diffeomorphism between a resolution $\mlps$ and a character variety with semisimple monodromies $\mathcal{M}_{\boldsymbol{\mathcal{S}}}$ is constructed thus proving the following theorem:

\begin{theorem}\label{th_diffeo_resol_ss}
$\mathcal{C}_{\boldsymbol{\mu},\boldsymbol{\sigma}}$ is a generic $k$-uple of conjugacy classes and $\mlps$ is the resolution of $\mcmus$. Then $\mlps$ is diffeomorphic to a character variety $\mathcal{M}_{\boldsymbol{\mathcal{S}}}$. With $\boldsymbol{\mathcal{S}}=\left(\mathcal{S}_1,\dots,\mathcal{S}_k \right)$ and $\mathcal{S}_j$ is the class of an element with centralizer in $\GL_n$ equal to $L^j\cong \GL_{{\underline{\mu}^j}'}$.
\end{theorem}


First the example of the sphere with four punctures and rank $n=2$ is studied in \ref{sect_example_sphere}. There, we can obtain the expected diffeomorphism using only tools from algebraic geometry. This example has been studied for a long time by Vogt \cite{Vogt} and Fricke-Klein \cite{Fricke}. The character varieties are affine cubic surfaces satisfying Fricke-Klein relation. Cubic surfaces and line over them have been extensively studied. They are classified for instance by Cayley \cite{Cayley}, see also Bruce-Wall \cite{Bruce_cubic}, Manin \cite{manin1986cubic} and Hunt \cite{hunt_b}. This rich theory proves that the minimal resolution is diffeomorphic to a character variety with semisimple monodromies. Both appear to be diffeomorphic to the projective plane blown up in six points minus three lines.

In general the construction of the diffeomorphism is performed in few steps and relies on analytic technics.

The first step is the Riemann-Hilbert correspondence, it gives a diffeomorphism between the resolution $\mathcal{M}_{\overline{{\mathcal{C}}}_{\boldsymbol{\mu},\boldsymbol{\sigma}}}$ and a de Rham moduli space of parabolic connections. Riemann-Hilbert correspondence was developed by Deligne \cite{deligne_equadif_rh}, and Simpson  for the filtered case \cite{Simpson}. Yamakawa proved that this correspondence induces a complex analytic isomorphism between moduli spaces \cite{Yamakawa}.

The second step is the non-Abelian Hodge theory, a diffeomorphism between de Rham moduli space and Dolbeault moduli space. It was established by Hitchin \cite{Hitchin87} and Donaldson \cite{Donaldson} for compact curves. Corlette \cite{Corlette} and Simpson \cite{Simpson_88} generalized it for higher dimensions. The parabolic version over non-compact curves was proved by Simpson \cite{Simpson}. This is the one needed here. It was generalized for higher dimension by Biquard \cite{Biquard97}. The relevant moduli spaces to obtain this correspondence as a diffeomorphism were introduced by Konno \cite{Konno} and Nakajima \cite{Nakajima_1996}. Biquard-Boalch \cite{biquard_boalch_2004} generalized further to  wild non-Abelian Hodge theory and constructed the associated hyperkähler moduli spaces. We use their construction of the moduli spaces. Biquard, García-Prada and Mundet i Riera \cite{Biquard_Garcia_Prada_Riera} established a parabolic non-Abelian Hodge correspondence for real groups, generalizing Simpson construction for $\GL_n$.

After the diffeomorphism from non-Abelian Hodge theory we use the method from Nakajima \cite{Nakajima_1996} for $\GL_2$ and Biquard, García-Prada, Mundet i Riera \cite{Biquard_Garcia_Prada_Riera} for real groups. The weights defining the moduli space of parabolic Higgs bundles are changed. This is done before going back to another de Rham moduli space thanks to non-Abelian Hodge theory in the other direction. The change of stability on the Dolbeault side induces a change of eigenvalues of the residue on the de Rham side.

Finally Riemann-Hilbert correspondence is applied in the other direction. It gives a diffeomorphism to a character variety where the eigenvalues $\sigma$ have been perturbed, the monodromies are now semisimple.

\subsection{Example of the sphere with four punctures and rank two}\label{sect_example_sphere}
We study the particular case $n=2$, $k=4$. Then the character varieties are affine cubic surfaces. The defining equation was known by Vogt \cite{Vogt} and Fricke-Klein \cite{Fricke}. The theory of cubic surfaces allows to obtain the expected diffeomorphism. Cubic surfaces and lines over them have been extensively studied. They are classified for instance by Cayley \cite{Cayley}, see also Bruce-Wall \cite{Bruce_cubic}, Manin \cite{manin1986cubic} and Hunt \cite{hunt_b}. This particular example of character varieties also appear in the theory of Painlevé VI differential equation. In this context resolution of cubic surfaces were studied by Inaba-Iwasaki-Saito \cite{Inaba_Iwasaki_Saito,Inaba_Iwasaki_Saito_I,Inaba_Iwasaki_Saito_II} with Riemann-Hilbert correspondence. It was also studied on the Dolbeault side by Hausel \cite{Hausel_compact}.

\subsubsection{Fricke relation}

We consider representations of the fundamental group of the sphere with four punctures $\mathbb{P}^1\setminus \left\lbrace p_1,...,p_4 \right\rbrace$. First we prescribe no particular condition on the monodromies around the puncture
\begin{equation*}
    \mathcal{R}:=\left\lbrace (X_1,\dots,X_4)\in\SL_2^4\left|X_1\dots X_4=\Id\right.\right\rbrace
\end{equation*}
The group $\GL_2$ acts by conjugation on $\mathcal R$, its center acts trivially, hence the action factors through an action of $\PGL_2$. Points of the following GIT quotient represent closed orbits for this action.
\begin{equation*}
    \mathcal{M}:=\mathcal{R}//\PGL_2:=\Spec \mathbb{C} \left[\mathcal R\right]^{\PGL_2}
\end{equation*}
where  $\C \left[\mathcal R\right]^{\PGL_2}$ are the invariants under the $\GL_2$ action in the algebra of functions of the affine variety $\mathcal{R}$. There is an explicit description of the variety $\mathcal{M}$ known by Vogt \cite{Vogt} and Fricke-Klein \cite{Fricke}, see also Goldman \cite{Goldman2009} for a detailed discussion and Boalch-Paluba \cite{Boalch_Paluba} for applications to $G_2$ character varieties. $\mathcal{M}$ is the affine variety given by the Fricke relation
\begin{equation}\label{eq_fricke_rel}
x y z + x^2 +y^2 + z^2 + A x + B y + C z + D =0
\end{equation}
with
\[
    x=\tr(X_2 X_3), \   \  y=\tr(X_1 X_3),\  \   z=\tr(X_1 X_2)
\]
and
\begin{eqnarray*}
A &=& -\tr(X_1)\tr(X_1 X_2 X_3) -\tr(X_2)\tr(X_3)  \\
B &=& -\tr(X_2)\tr(X_1 X_2 X_3)-\tr(X_1)\tr(X_3)  \\
C &=& -\tr(X_3)\tr(X_1 X_2 X_3) -\tr(X_1)\tr(X_2) \\
D &=& \tr(X_1)\tr(X_2)\tr(X_3)\tr(X_1 X_2 X_3) + \tr(X_1)^2 \\ 
& &+ \tr(X_2)^2 + \tr(X_3)^2 + \tr(X_1 X_2 X_3)^2 -4.
\end{eqnarray*}

The character varieties we are interested in are obtained by specifying the Zariski closure of the conjugacy class of each $X_i$. First we assume that they are all semi-simple regular with determinant $1$. For $i=1,\dots 4$;  $\mathcal{S}_i$ is the conjugacy class of  
\begin{equation}\label{conj2}
     \left(\begin{array}{cc}
        \lambda_i & 0  \\
        0 & \lambda_i\inv 
    \end{array}\right).
\end{equation}
$\boldsymbol{\mathcal{S}}=\left( \mathcal{S}_1,\dots, \mathcal{S}_4\right)$ is assumed to be generic.
In terms of invariant functions, $X_i\in\mathcal{S}_i$ for all $i$, if and only if
\begin{eqnarray*}
    \tr(X_i) &=& \lambda_i+\lambda_i \inv \; \text{ for }1 \le i \le 3  \\
    \tr(X_1 X_2 X_3) &=& \lambda_4+\lambda_4\inv .
\end{eqnarray*}
Then Fricke relation translates in next proposition.

\begin{proposition}
The character variety $\mathcal{M}_{\boldsymbol{\mathcal{S}}}$ is a smooth cubic surface in $\mathbb{A}^3$ given by Fricke relation \eqref{eq_fricke_rel} with coordinates $x,y$ and $z$ and constants $A, B, C$ and $D$.
\end{proposition}

Now consider non-semisimple conjugacy classes $\boldsymbol{\mathcal{C}}=\left(\mathcal{C}_1,\mathcal{C}_2,\mathcal{C}_3,\mathcal{C}_4 \right)$. With $\mathcal{C}_1$  the conjugacy class of
\[
\begin{pmatrix}
-1 & 1 \\
0 & -1
\end{pmatrix}
\]
and $\mathcal{C}_2=\mathcal{C}_3 = \mathcal{C}_4$ are the conjugacy classes of
\[
\begin{pmatrix}
1 & 1 \\
0 & 1
\end{pmatrix}.
\]
Note that this $4$-uple of conjugacy classes is generic. The $(X_i)_{1\le i \le 4}$ are already assumed to have determinant $1$, then $X_1$ belongs to the closure $\overline{\mathcal{C}}_1$ if and only if
\[
\tr X_1 = -2.
\]
Similarly the condition $(X_2, X_3, X_4)\in \overline{\mathcal{C}}_2\times \overline{\mathcal{C}}_3 \times
\overline{\mathcal{C}}_4$ is equivalent to
\[
\tr X_2 = \tr X_3 = \tr (X_1 X_2 X_3)= 2.
\]
Substituting this parameters in Fricke relation, the character variety is again a cubic surface in $\mathbb{A}^3$ with equation:
\begin{equation}\label{eq_cayley_char}
x y z + x^2 +y^2 +z^2 - 4 =0.
\end{equation}
This cubic surface has exactly four singularities at $(-2,-2,-2)$, $(-2,2,2)$, $(2,-2,2)$ and $(2,2,-2)$. 
The classification of cubic surfaces (see for instance Bruce-Wall \cite{Bruce_cubic}) gives the following theorem:
\begin{theorem}\label{th_compactification_cayley}
After compactification in $\mathbb{P}^3$, the character variety $\mathcal{M}_{\overline{\boldsymbol{\mathcal{C}}}}$ is Cayley's nodal cubic, the only cubic surface with four singularities.
\end{theorem}
This particular character variety was studied by Cantat-Loray \cite{cantat_loray} in the context of Painlevé VI.

In this example, using only elementary algebraic geometry, we can prove that the minimal resolution of $\mathcal{M}_{\overline{\boldsymbol{\mathcal{C}}}}$ is diffeomorphic to the character varieties with semisimple monodromies $\mathcal{M}_{\boldsymbol{\mathcal{S}}}$. We shall see that both varieties are obtained as the plane blown-up in six points minus three lines.

\subsubsection{Projective cubic surfaces}

Let us recall an important result in the classification of cubic surfaces. Smooth projective cubic surfaces in $\mathbb{P}^3$ can be constructed by a blow-up of $\mathbb{P}^2$ in six points.

Let $\boldsymbol{P}=(P_1,...,P_6)$ be six distinct points in the projective plane $\mathbb{P}^2$. The blow-up of $\mathbb{P}^2$ with respect to those six points is denoted $Y_{\boldsymbol{P}}\to \mathbb{P}^2$.

\begin{definition}[Generic configuration for six points in $\mathbb{P}^2$]
Such a configuration $\boldsymbol{P}$ of $6$ points in $\mathbb{P}^2$ is called generic if no three of them lie on a line and no five of them lie on a conic.
\end{definition}
The two following theorems are well-known results about cubic surfaces, see for instance Manin \cite{manin1986cubic} and Hunt \cite{hunt_b}.
\begin{theorem}\label{th_cubic_six}
Up to isomorphism, smooth projective cubic surfaces in $\mathbb{P}^3$ are obtained as $\mathbb{P}^2$ blown-up in six points in generic position.
\end{theorem}
\begin{theorem}\label{th_cayley_six}
If the six points $\boldsymbol{P}=(P_1,...,P_6)$ are the intersection of four lines $(L_1,...,L_4)$ in $\mathbb{P}^2$, then $Y_{\boldsymbol{P}}$ is isomorphic to a minimal resolution of singularities of Cayley's nodal cubic.
\end{theorem}

Up to diffeomorphism, the manifold obtained by $\mathbb{P}^2$ blown-up in six distinct points, does not depend on the position of the points. This implies next proposition.

\begin{proposition}
The minimal resolution of the projective Cayley's nodal cubic is diffeomorphic to a smooth projective cubic surface. Both are obtained as the projective plane $\mathbb{P}^2$ blown-up in six points.
\end{proposition}

\subsubsection{Lines on cubic surfaces}
So far we saw that the minimal resolution of the projective Cayley's nodal cubic is diffeomorphic to a smooth projective cubic surface. However the variety we are interested in are not projective, they are affine. By Theorem \ref{th_compactification_cayley} the variety $\mathcal{M}_{\overline{\boldsymbol{\mathcal{C}}}}$ is the projective Cayley's nodal cubic minus three lines at infinity. Those three lines are given by the equation $x y z =0$, they form a triangle. Similarly the variety $\mathcal{M}_{\boldsymbol{\mathcal{S}}}$ is a smooth projective cubic surface minus the triangle at infinity $x y z = 0$. This triangle at infinity is a particular case of a general situation studied by Simpson \cite{simpson16} for $n=2$ and any number of punctures $k$.

The theory of lines on cubic surfaces has been thoroughly studied. See for instance Cayley \cite{Cayley}, Bruce-Wall \cite{Bruce_cubic}, Manin \cite{manin1986cubic} and Hunt  \cite{hunt_b}. 
\begin{proposition}[$27$ lines on smooth projective cubic surface]\label{proposition_27}
There are $2 7$ lines on a smooth projective cubic surface. They all have a nice description in terms of $\mathbb{P}^2$ blown-up in six points $(P_1,\dots,P_6)$.
\begin{itemize}
    \item Six of them are exceptional divisors $E_i$ over $P_i$.
    \item Fifteen of them are the strict transform $\widetilde{L}_{i,j}$ of the line through $P_i$ and $P_j$.
    \item Six of them are the strict transform $\widetilde{C}_j$ of the conic through all $P_i$ except $P_j$.
\end{itemize}
Following picture is an example of six generic points in the plan, the line $L_{1,6}$ as well as the conic $C_6$ are drawn. 
\end{proposition}
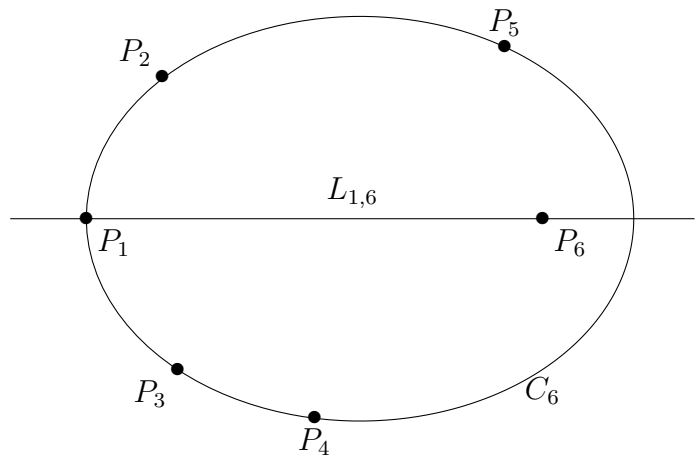
\[
\begin{tikzpicture}
\coordinate (P_1) at (-6,0);
\coordinate (P_6) at (0,0);

\draw (P_1) node[below right] {$P_1$}
node{$\bullet$} ;
\draw (-5, 1.88) node[above left] {$P_2$} node{$\bullet$};
\draw (-4.8, -2) node[below left] {$P_3$} node{$\bullet$} ;
\draw (-3, -2.64) node[below] {$P_4$} node{$\bullet$};
\draw (-0.5, 2.28) node[above] {$P_5$} node{$\bullet$};
\draw (P_6) node[below right] {$P_6$}
node{$\bullet$} ;

\draw (-7,0) -- (2,0) node[midway,above] {$L_{1,6}$};

\draw [domain= -180 : 180 , samples = 200] plot (\x : {3/(1.5 +  cos(\x))} );
\draw (0,-2.28) node {$C_6$};
\end{tikzpicture}
\]
Now consider six points not in generic position. Take four lines $(L_1,\dots,L_4)$ in $\mathbb{P}_2$ with exactly six intersections $(P_1,\dots,P_6)$, those lines are black in next figure. Consider the three lines $L_{1,2}, L_{3,4}$ and $L_{5,6}$ with $L_{i,j}$ containing $P_i$ and $P_j$, those lines are blue in next figure. Up to relabelling we may assume $L_{i,j}\ne L_k$ for all $i,j,k$. Cayley's nodal cubic is obtained by blowing up the six points and then blowing down the strict transform of the four lines $(L_1,\dots,L_4)$. The four points image of this four lines under the blow-down are exactly the four singular points. See  Hunt \cite[Chapter 4]{hunt_b} for more pictures.

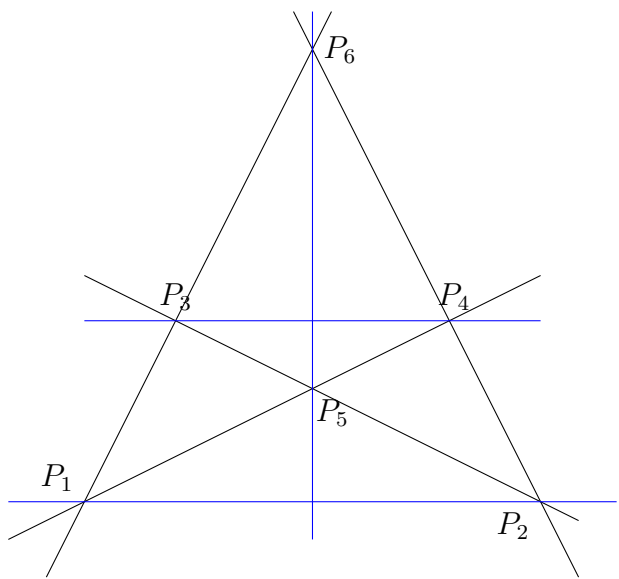
\[
\begin{tikzpicture}
\draw (-0.5,-1)--(3.25,6.5) ;
\draw (-1,-0.5)--(6,3);
\draw (2.75,6.5)--(6.5,-1);
\draw (0,3)--(6.5,-0.25);

\draw [blue] (-1,0)--(7,0);
\draw [blue] (3,-0.5)--(3,6.5);
\draw [blue] (0,2.4)--(6,2.4);

\draw (0,0) node[above left] {$P_1$};
\draw (1.2,2.4) node[above] {$P_3$};
\draw (4.5,2.4) node[above right] {$P_4$};
\draw (2.9,1.5) node[below right] {$P_5$};
\draw (3,6) node[right] {$P_6$};
\draw (6,0) node[below left] {$P_2$};
\end{tikzpicture}
\]

\begin{proposition}[lines on Cayley's nodal cubic]
There are $9$ lines on Cayley's nodal cubic.
\begin{itemize}
    \item Six of them are the exceptional divisors $E_i$ over $P_i$.
    \item Three of them are the strict transform of $L_{1,2}, L_{3,4}$ and $L_{5,6}$.
\end{itemize}
\end{proposition}
\begin{proposition}
The variety $\mathcal{M}_{\overline{\boldsymbol{\mathcal{C}}}}$ is Cayley's nodal cubic minus the images of $L_{1,2}, L_{3,4}$ and $L_{5,6}$.
\end{proposition}
\begin{proof}
We saw that $\mathcal{M}_{\overline{\boldsymbol{\mathcal{C}}}}$ is Cayley's nodal cubic minus the three lines at infinity $xyz=0$. Those three lines does not contain any of the four singularities. Therefore they are not the image of the exceptional divisors. Then they must be the three remaining lines, the blue lines on the picture.
\end{proof}
\begin{theorem}
The character variety with generic semisimple conjugacy classes at punctures $\mathcal{M}_{\boldsymbol{\mathcal{S}}}$ is diffeomorphic to the minimal resolution of singularities of the character variety $\mathcal{M}_{\overline{\boldsymbol{\mathcal{C}}}}$. Both are obtained as the projective plane $\mathbb{P}^2$ blown up in six points $(P_1,\dots,P_6)$ minus three lines $\widetilde{L}_{1,2}, \widetilde{L}_{3,4}, \widetilde{L}_{5,6}$.
\end{theorem}
\begin{proof}
The statement about the minimal resolution of $\mathcal{M}_{\overline{\boldsymbol{\mathcal{C}}}}$ follows from previous proposition. $\mathcal{M}_{\boldsymbol{\mathcal{S}}}$ is a smooth projective cubic surface minus three lines forming a triangle. As those three lines intersect each other they cannot be any triple among the $27$ lines over the surface, there are some restriction:
\begin{itemize}
    \item Exceptional divisor $E_i$ do not intersect each other.
    \item Strict transform $\widetilde{C}_j$ do not intersect each other.
    \item Strict transforms of two distinct line containing a same point $P_i$ do not intersect.
\end{itemize}
Therefore the only possible triples of lines forming a triangle on a smooth cubic surface have the following form:
\begin{enumerate}
    \item $(\widetilde{L}_{1,2},\widetilde{L}_{3,4}, \widetilde{L}_{5,6})$ 
    \item $(E_1, \widetilde{L}_{1,6},\widetilde{C}_6) $.
\end{enumerate}
The first case is exactly the expected result. The second case is illustrated by the picture below Proposition \ref{proposition_27}, the conic $C_6$ and the line $L_{1,6}$ are drawn. To relate the second case to the first, proceed in two steps. First $\mathbb{P}^2$ is blown-up in the three points $P_1, P_2$ and $P_3$. The resulting variety is blown-down along $\widetilde{L}_{1,2},\widetilde{L}_{1,3}$ and $\widetilde{L}_{2,3}$ (three lines with self-intersection $-1$). The variety obtained is again isomorphic to $\mathbb{P}_2$. We consider this copy of the projective plane as the starting point. This plane is blown up in six points $(P_1',\dots,P_6')$ with
\begin{itemize}
    \item $P_1'$ the blow-down of $\widetilde{L}_{2,3}$
    \item $P_2'$ the blow-down of $\widetilde{L}_{1,3}$
    \item $P_3'$ the blow-down of $\widetilde{L}_{1,2}$
    \item $P_j'$ the image of $P_j$ for $j=4,5,6$.
\end{itemize}
The construction obtained from the new copy of $\mathbb{P}^2$ and the points $(P_1',\dots,P_6')$ are labelled with a prime. Then the triple $\left(E_1, \widetilde{L}_{1,6},\widetilde{C}_6\right) $ becomes $\left(\widetilde{L}'_{2,3},\widetilde{L}'_{1,6},\widetilde{L}'_{4,5}\right)$. In any cases the triangle of lined removed at infinity has the expected form.
\end{proof}
\begin{remark}
There is an action of the Weyl group of $\mathbb{E}_6$ on the configuration of the $27$ lines on a smooth cubic surface. The Dynkin diagramm of $\mathbb{E}_6$ is
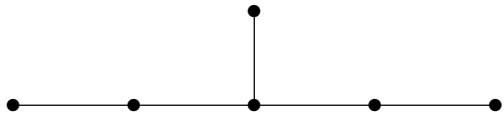
\[\begin{tikzcd}
	&& \bullet \\
	\bullet & \bullet & \bullet & \bullet & \bullet
	\arrow[no head, from=2-1, to=2-2, start anchor=center, end anchor = center ]
	\arrow[no head, from=2-3, to=2-4,start anchor=center, end anchor = center ]
	\arrow[no head, from=2-4, to=2-5,start anchor=center, end anchor = center ]
	\arrow[no head, from=1-3, to=2-3,start anchor=center, end anchor = center ]
	\arrow[no head, from=2-2, to=2-3,start anchor=center, end anchor = center ]
\end{tikzcd}\]
The generator of the upper vertex corresponds to the transformation previously described sending $\left(E_1, \widetilde{L}_{1,6},\widetilde{C}_6\right) $ to  $\left(\widetilde{L}'_{2,3},\widetilde{L}'_{1,6},\widetilde{L}'_{4,5}\right)$. See Hartshorne \cite[V-Exercise 4.11]{hartshorne}.
\end{remark}

\subsection{Moduli spaces}\label{sect_moduli_spaces}

In general the construction of the diffeomorphism relies on analytical techniques  and go through various moduli spaces. Let $\Sigma$ a compact Riemann surface endowed with a complex structure, then $\Sigma$ is seen as a smooth complex projective curve. Let $D$ the divisor $D=p_1+\dots +p_k$ for $k$ distincts points $p_1,\dots,p_k$. 

\subsubsection{de Rham moduli space}\label{subsect_dR_alg}
Parabolic holomorphic bundles were introduced by Mehta-Seshadri \cite{Mehta_Seshadri}, they generalized Narasimhan-Seshadri \cite{Narasimhan_Seshadri} result to the parabolic case. Parabolic bundles appear in various area in mathematics and physics, for instance Pauly \cite{Pauly} related those parabolic bundles with conformal field theory. In this section some definitions are recalled.

\begin{definition}[Filtered holomorphic bundles]
A filtered holomorphic bundle is the data of a holomorphic vector bundle $E$ together with filtrations of $E^j$ the fiber of $E$ at $p_j$ for $j=1,\dots,k$
\[
\left\lbrace 0\right \rbrace = E^j_0 \subset E^j_1\subset\dots\subset E^j_{m_j} = E^j.
\]
The type $\tau$ of the filtration is defined by
\[
\tau^j_i =\dim E^j_i/E^j_{i-1} 
\]
for $j=1,\dots,k$ and $i=1,\dots,m_j$.
\end{definition}
\begin{definition}[parabolic degree]
Let $E$ a filtered holomorphic bundle of type $\tau$. Let $\beta= \left(\beta^j_i\right)_{ \substack{1\le j \le k \\ 1 \le i \le m_j} }$ with $\beta_i^j\in\mathbb{R}$ a stability parameter. The parabolic degree of $E$ is
\[
\pdeg_{\beta} E =\deg E + \sum_{i,j} \beta^j_i \dim \left(E^j_i/E^j_{i-1}\right).
\]
\end{definition}

Let $E$ a holomorphic vector bundle on $\Sigma$. A logarithmic connection on $E$ is a map of sheaves $D: E\to E\otimes \Omega^1_{\Sigma}(\log D)$ satisfying Leibniz rule
\[
D(f s) = df\otimes s + f D(s)
\]
for all $f$ holomorphic function and $s$ section of $E$.

For $z$ a coordinate vanishing at a point $p_j$, in a trivialization of $E$ in a neighborhood of this point the connection reads
\[
D = d + A(z) \frac{d z}{z}.
\]
$A(0)$ is called the residue of $D$ at $p_j$ and denoted by $\Res_{p_j} D$

Fix some parabolic weights $\beta^j_i\in \left[ 0, 1\right[$ satisfying $\beta^j_i > \beta^j_{i-1}$. For $j=1,\dots, k$ and $i=1,\dots,m_j$ fix $A^j_i \in \C$ to specify a polar part. A logarithmic connection $(E,D)$ is compatible with the parabolic structure if the endomorphism 
\[
\Res_{p_j} D : E^j \to E^j
\]
satisfies $\left(\Res_{p_j} D \right)   E^j_i   \subset E^j_{i}$. A logarithmic connection compatible with the parabolic structure is called a parabolic connection.

It is compatible with the specified polar part if in addition the map induced by $\Res_{p_j} D$ on the graded spaces $   E^j_i/E^j_{i-1} $ is $A^j_i \Id$.
A logarithmic connection compatible with the parabolic structure is $\beta$-semistable if and only if, for all sub-bundle $F\subsetneq E$ preserved by $D$
\[
\frac{\pdeg_{\beta} F}{\rank F} \le \frac{\pdeg_{\beta} E}{\rank E}
\]
it is stable if the inequality is strict unless $F=0$. Two pairs of filtered holomorphic bundle and parabolic connections $(E,D)$ and $(E',D')$ are isomorphic if there is an isomorphism of holomorphic bundle $f:E\to E'$ compatible with the filtrations and such that $\left(f\otimes Id \right)\circ D = D'\circ f$. A connection is flat if its curvature vanishes.

\begin{notations}[de Rham moduli space]
The de Rham moduli space $\mathcal{M}^{d R}_{A,\beta}$ classifies isomorphism classes of flat $\beta$-stable parabolic connections with prescribed polar part $A$ and parabolic degree $0$.
\end{notations}

\subsubsection{Filtered local systems and resolutions of character varieties}\label{subsect_resol_filtered_local}
\begin{definition}[Filtered local system]
A filtered local system is a local system $\mathcal{L}$ over $\Sigma\setminus\left\lbrace p_1,\dots,p_k\right\rbrace$ together with a filtration of the restrictions $\mathcal{L}_{\left|U_j\right.}$ to $U_j$ some punctured neighborhood of $p_j$. Namely for all $j=1,\dots,k$ there are local systems $\mathcal{L}^j_i$ such at
\[
0 = \mathcal{L}^j_0 \subsetneq \mathcal{L}^j_1 \subsetneq \dots\subsetneq \mathcal{L}^j_{m_j} = \mathcal{L}_{\left|U_j\right.}.
\]
The type $\tau$ of the filtered local system is defined by
\[
\tau^j_i := \rank \mathcal{L}^j_i/ \mathcal{L}^j_{i-1}.
\]
\end{definition}

\begin{definition}[Parabolic degree of a filtered local system]
Let $\gamma = \left(\gamma^j_i\right)_{ \substack{1\le j \le k \\ 1 \le i \le m_j} }$ a stability parameter. The parabolic degree of the filtered local system is defined by
\[
\pdeg_{\gamma} \mathcal{L} = \sum_{i,j} \gamma^j_i \rank \mathcal{L}^j_i / \mathcal{L}^j_{i-1}
\]
\end{definition}

A filtered local system $\mathcal{L}$ is $\gamma$-semistable if and only if for all sub local system $ 0 \subsetneq \mathcal{L}' \subsetneq \mathcal{L}$
\[
\frac{\pdeg_{\gamma} \mathcal{L}'}{\rank \mathcal{L}'} \le \frac{\pdeg_{\gamma} \mathcal{L}}{\rank \mathcal{L}}
\]
it is $\gamma$-stable if the inequality is strict.

Consider a character variety $\mcmus$ with a resolution of singularities $\widetilde{\mathcal{M}}_{\boldsymbol{L},\boldsymbol{P}, \boldsymbol{\sigma} }$. By the usual equivalence of category between local systems and representations of the fundamental group, the character variety $\mcmus$ is the moduli space of local system with monodromy around $p_{j}$ in $\overline{\mathcal{C}}_{\underline{\mu}^j,\sigma^j}$. This correspondence extends to the resolution $\mlps$ and the moduli space of filtered local system.
\begin{proposition}\label{proposition_filtered_moduli_resolu}
$\mlps$ is the moduli space of filtered local system with filtration around $p_j$ of type ${\underline{\mu}^j}'$ and such that the endomorphism induced by the monodromy on $\mathcal{L}^j_i/\mathcal{L}^j_{i-1}$ is $\sigma^j_i \Id$.
\end{proposition}
\begin{proof}
An element $g_j P^j \in \GL_n /P^j$ identifies with a partial flag of type ${\underline{\mu}^j}'$.
The condition $g_j\inv X_j g_j \in \sigma^j U_{P^j}$
is exactly that the partial flag is preserved by $X_j$
and that the induced endomorphism on the graded spaces are $\sigma_i^j \Id$.
Note that we study only character varieties for generic choices of conjugacy classes at punctures. For such a generic choice, the stability parameter is irrelevant as the local system does not admit any sub local system.
\end{proof}

\subsubsection{Dolbeault moduli space}\label{subsect_dol_alg}
A parabolic Higgs bundle is a pair $(E,\phi)$ with $E$ a filtered holomorphic vector bundle on $\Sigma$ and a Higgs field $\phi : E \to E\otimes \Omega^1(\log D) $ such that $\Res_{p_j}\phi (E^j_i )\subset E^j_i$. Let $\alpha = \left(\alpha^j_i\right)_{ \substack{1\le j \le k \\ 1 \le i \le n_j} }$ a stability parameter. A parabolic Higgs bundle $(E,\phi)$ is $\alpha$-semistable if and only if for all $0\subsetneq F\subsetneq E$ a sub bundle preserved by $\phi$
\[
\frac{\pdeg_{\alpha} F}{\rank F} \le \frac{\pdeg_{\alpha} E}{\rank E}.
\]
it is $\alpha$-stable if the inequality is strict.
As for the parabolic connections, it is interesting to specify the residue of the Higgs field. For all $i,j$ fix a semisimple adjoint orbit $B^j_i$ in $\gl_{\nu^j_i}$. The parabolic Higgs bundle has the prescribed residue if, in an holomorphic trivialization, the map induced on $E^j_i / E^j_{i-1}$ by the residue lies in the adjoint orbit $B^j_i$. Note that contrarily to the parabolic connections, the prescribed adjoint orbits on the graded spaces are not necessarily central. In fact much more general polar parts are considered by Biquard-Boalch \cite{biquard_boalch_2004}, we restrict here to what is necessary for our purpose.

\begin{notations}[Dolbeault moduli space]
The Dolbeault moduli space $\mathcal{M}^{D o l}_{B,\alpha}$ classifies isomorphism classes of
$\alpha$-stable parabolic Higgs bundles with prescribed residue $B$ and parabolic degree $0$.
\end{notations}

\subsubsection{Various steps of the diffeomorphism}\label{sect_analytic_moduli}
In the remaining of the section, analytic construction of the moduli spaces are recalled. Those spaces are endowed with a manifold structure. They will be used to obtain a diffeomorphism from a resolution $\widetilde{\mathcal{M}}_{\boldsymbol{L},\boldsymbol{P}, \boldsymbol{\sigma} }$ to a character variety $\mathcal{M}_{\mathcal{\boldsymbol{S}}}$ with semisimple conjugacy classes at punctures. The picture is is the following:
\begin{equation}\label{diagram_path_ch4}
\begin{tikzcd}
\widetilde{\mathcal{M}}_{\boldsymbol{L},\boldsymbol{P}, \boldsymbol{\sigma} } \arrow[r,"\text{R.H}"] & \mathcal{M}^{d R}_{A,\beta} \arrow[r,"\text{N.A.H}"] & \mathcal{M}^{Dol}_{B,\alpha} \arrow[d,"\alpha\mapsto \widetilde{\alpha}"] \\
\mathcal{M}_{\boldsymbol{\mathcal{S}}} & \arrow[l,"\text{R.H}"] \mathcal{M}^{d R}_{\widetilde{A},\widetilde{\beta}} &\arrow[l,"\text{N.A.H}"] \mathcal{M}^{Dol}_{B,\widetilde{\alpha}}
\end{tikzcd}
\end{equation}
All the arrows are diffeomorphisms, R.H stands for Riemann-Hilbert correspondence and N.A.H for non-Abelian Hodge theory. The vertical arrow accounts for a change of stability parameter $\alpha\to\widetilde{\alpha}$. This is similar to Biquard, García-Prada and Mundet i Riera \cite[Theorem 7.10]{Biquard_Garcia_Prada_Riera}. It is detailed in the remaining of the section for this particular application.

\subsection{Local model}\label{subsect_local_model}
In this section the local model used by Biquard-Boalch \cite{biquard_boalch_2004} to construct moduli spaces is recalled.

\subsubsection{Local model for Riemann-Hilbert correspondence}\label{subsub_local_RH}
Before constructing the moduli space, let us present what happens locally, near a puncture, and how the parameters of the moduli spaces are related. Consider $\mathcal{L}$ a rank $n$ filtered local system on a punctured disk $\mathbb{D}^0=\left\lbrace z\in \mathbb{C} \left| 0 < |z|< 1 \right.\right\rbrace$. It is assumed that the monodromy induces a central endomorphism on the successive quotients of the filtration. The monodromy $X$ has eigenvalues $\sigma_i$ with multiplicity $\nu_i$ for $1\le i \le l$. We assume the filtration of the local system is finer than a filtration spanned by generalized eigenspaces of $M$. Then in a trivialization $(l_j)_{1\le j \le n}$ compatible with the filtration, the monodromy reads
\[
X =
\begin{pmatrix}
X_{\sigma_1} &  *& \\
0 & X_{\sigma_2} & *\\
\vdots  &  0 & \ddots 
\end{pmatrix}
\]
with $X_{\sigma_i}$ a block of size $\nu_i$ with further decomposition
\[
X_{\sigma_i} =
\begin{pmatrix}
\sigma_i \Id_{{\mu^i_1}'} &  *& \\
0 & \sigma_i \Id_{{\mu^i_2}'} & *\\
\vdots  &  0 & \ddots
\end{pmatrix}
\]
The type of the filtration is $\underline{\mu}'=({\mu^1_1}',{\mu^1_2}',\dots,{\mu^2_1}',{\mu^2_2}',\dots)$. Let $A_i \in \C$ such that 
\[
\exp(- 2 i \pi A_i)=\sigma_i
\]
and $0 \le \Ree A_i<1$. Then $A$ is the diagonal matrix with diagonal coefficients
\[
\left(\underbrace{A_1,\dots,A_1}_{\nu_1},\dots,\underbrace{A_l,\dots,A_l}_{\nu_l}\right).
\]
Let $a$ a block strictly upper triangular matrix such that $\exp\left(-2i\pi (A+a)\right)=X$. Define $E$ a rank $n$ holomorphic bundle on the disk $\mathbb{D}=\left\lbrace z\in\mathbb{C} \left| |z|<1\right. \right\rbrace$ spanned by $\tau_j=e^{(A+a) \log z} l_j$ for $1\le j\le n$. Let $D$ the parabolic connection on $E$ defined in the holomorphic trivialization $(\tau_j)_{1\le j\le n}$ by
\begin{eqnarray*}
D &=& d + \frac{A+a}{z} dz \\
&=& D_0 + \frac{a}{z} dz
\end{eqnarray*}
Then the filtered local system $\mathcal{L}$ is nothing but the local system of flat sections of the parabolic connection $(E,D)$. This describes locally the Riemann-Hilbert correspondence between a resolution of a character variety and a de Rham moduli space.

\subsubsection{Metric and parabolic structure}
The connection $D_0$ will be the local model for parabolic connections:
\[
D_0 = d + \frac{A}{z} dz
\]
with $A$ diagonal. In order to continue the path presented in Diagram \eqref{diagram_path_ch4}, we need to introduce an Hermitian metric. It will be related to a choice of stability parameters. Chose some stability parameters $\beta_{r,s}\in [0,1[$ for each graded spaces of the filtration of type $\underline{\mu}$. Introduce a diagonal matrix $\beta$ with diagonal coefficients
\[
\left(\beta_1,\beta_2,\dots,\beta_n\right) := \left( \overbrace{\underbrace{\beta_{1,1},\dots,\beta_{1,1}}_{{\mu^1_1}'},\underbrace{\beta_{1,2},\dots,\beta_{1,2}}_{{\mu^1_2}'},\dots}^{\nu_1},\dots,\overbrace{  \underbrace{\beta_{l,1},\dots,\beta_{l,1}}_{{\mu^l_1}'},\underbrace{\beta_{l,2},\dots,\beta_{l,2}}_{{\mu^l_2}'},\dots   }^{\nu_l} \right)
\]
so that the $\beta_i$ are the $\beta_{r,s}$ repeated according to the multiplicities ${\mu^r_s}'$. Moreover assume that $\beta_i \ge \beta_{i+1}$ and $\beta_{r,s}\ne \beta_{u,v}$ if $(r,s)\ne(u,v)$.
\begin{remark}\label{remark_various_beta}
In this local model there is just one puncture $p_1$, therefore the stability parameter intorduced in \ref{subsect_dR_alg} are $\left(\beta^1_i\right)_{1\le i \le m_1}$. They are related to the stability parameters introduced in this section by
\[
\left(\beta^1_1,\beta^1_2,\dots,\beta^1_{m_1}\right)=\left(\beta_{1,1},\beta_{1,2},\dots,\beta_{2,1},\beta_{2,2},\dots\right)
\]
We apologize for the multiplication of similar notations. $(\beta^1_i)_{1\le i \le m_1}$ are adapted to the algebraic definition of stability whereas $(\beta_{r,s})_{ \substack{1\le r \le l \\ 1 \le s \le \mu^r_1  }  }$ are adapted to the description of the connections and $(\beta_1,\beta_2,\dots,\beta_n)$ to explicit construction of trivializations. 
\end{remark}
Define a Hermitian metric $h$ on $E$ such that $|\tau_j| = |z|^{\beta_j}$. This metric determines the filtration of $E$:
\[
E_i = \left\lbrace s\in E \left| | s(z) |_h = \mathcal{O}\left( \left| z \right|^{\beta^1_i} \right) \right. \right\rbrace.
\]
with $|\dots|_h$ the norm with respect to the metric $h$. We obtained an Hermitian vector bundle $\mathbb{E}$ on the disk $\mathbb{D}$  with an orthonormal trivialization $\left(\frac{\tau_j}{|z|^{\beta_j}}\right)_{1\le j \le n}$.
\begin{notations}
The symbol $\mathbb{E}$ represents a vector bundle in the sense of differential geometry, with smooth transition functions; whereas the symbol $E$ represents a holomorphic bundle. 
\end{notations}
The parabolic connection $D_0$ on the holomorphic bundle $E$ induces a connection on $\mathbb{E}$, in the orthonormal trivialization $\left(\frac{\tau_j}{|z|^{\beta_j}}\right)_{1\le j \le n}$ it reads
\[
D_0 = d + \left( A -\frac{\beta}{2} \right) \frac{d z}{z} -\frac{\beta}{2} \frac{d \overline{z}}{\overline{z}}.
\]


\subsubsection{Local behaviour for non-Abelian Hodge theory}\label{subsection_local_naht}
$D_0$ decomposes as unitary connection plus a self-adjoint term
\[
D_0= D_0^h + \Phi_0.
\]
In the orthonormal trivialization $\left(\frac{\tau_j}{|z|^{\beta_j }}\right)_{1\le j \le n}$
\[
D_0^h = d + \frac{A}{2}  \frac{d z}{z} - \frac{A^{\dagger}}{2} \frac{d\overline{z}}{\overline{z}}
\]
and
\[
\Phi_0 = \frac{1}{2} \left( A\frac{d z}{z} + A^{\dagger}\frac{d\overline{z}}{\overline{z}} - \beta \frac{d z}{z} - \beta \frac{d \overline{z}}{\overline{z}}\right)
\]

Consider the basis $(e_j)_{1\le j \le n} $ defined by
\[
e_j := \frac{\tau_j}{\left|z\right|^{\beta_j- i \Imm A_j}}
\]
with $\Imm A_j$ the imaginary part of the $j$-th diagonal term of the matrix $A$.
\begin{notations}[Canonical form]\label{notations_canonical_form}
The expression of $D_0$ in the orthonormal trivialization $(e_j)_{1\le j\le n}$ is
\begin{eqnarray*}
D_0 &=& D_0^h + \Phi_0 \\
D_0^h &=& d +\frac{1}{2} \Ree (A) \left(\frac{d z}{z}- \frac{d\overline{z}}{\overline{z}} \right) \\ 
\Phi_0 &=& \frac{1}{2} \left( A\frac{d z}{z} + A^{\dagger}\frac{d\overline{z}}{\overline{z}} - \beta \frac{d z}{z} - \beta \frac{d \overline{z}}{\overline{z}}\right).
\end{eqnarray*}
Such expressions will be referred to as canonical forms.
\end{notations}

Let $\overline{\partial}^F$ be the $(0,1)$-part of $D_0^h$ and $\theta_0$ the $(1,0)$-part of $\Phi_0$. In the basis $(e_j)_{1\le j \le n} $ one has
\[
\overline{\partial}^F = \drond - \frac{1}{2} \Ree ( A) \frac{d\overline{z}}{\overline{z}}.
\]
This operator defines an holomorphic bundle over the punctured disk with holomorphic sections killed by $\drond^F$. This holomorphic bundle can be extended over the puncture to an holomorphic bundle $F$, taking as a basis of holomorphic sections $(f_j)_{1\le j \le n}$ defined by
\[
f_j = |z|^{\alpha_j} e_j.
\]
with $\alpha_j$ the real part of the $j$-th diagonal term of the matrix $A$. Then
\[
|f_j|_h = \left|z\right|^{\alpha_j}
\]
Similarly to the correspondence \ref{remark_various_beta} between $(\beta_1,\dots,\beta_n)$ and $(\beta^1_i)_{1\le i\le m_1}$, a stability parameter $(\alpha^1_i)_{1\le i \le n_1}$ is associated to $(\alpha_1,\dots,\alpha_n)$.
This stability parameter provides a parabolic structure
\[
F_i =\left\lbrace s\in F\left|\; |s|_h=\mathcal{O}\left(\left|s(z)\right|^{\alpha_i^1}\right) \right.  \right\rbrace 
\]
Note that the holomorphic bundle $F$ is different from the holomorphic bundle $E$. Even the type of the parabolic structure differ, $E$ is of type $\underline{\mu}'$ whereas $F$ is of type $\nu$.

Note that $\theta_0$, the $(1,0)$ part of $\Phi_0$, provides an Higgs field:
\[
\theta_0 = \frac{1}{2} \left(A  - \beta \right)\frac{d z}{z}.
\]
This is the local behaviour of the non-Abelian Hodge theory for the model connection.
To summarize, starting from a flat parabolic connection $D_0$ with polar part $A$, a metric $h$ and a parabolic structure $\beta$ we obtain a parabolic Higgs bundle with residue of the Higgs field $B$ and parabolic structure $\alpha$. The relation between those parameters are as described by Simpson \cite{Simpson}
\begin{equation}\label{eq_nah_parameters}
\begin{array}{ccc}
B &=& \frac{1}{2}\left( A - \beta \right) \\
\alpha &=& \Ree A.
\end{array}
\end{equation}

\subsubsection{Local description of weighted Sobolev spaces}\label{subsect_local_weighted}
\begin{definition}[Weighted $L^2$ spaces]
$r=|z|$ is the radial coordinate on the disk, for $\delta$ real, $L^2_{\delta}$ is the space of function $f$ on the disk such that $\frac{f}{r^{\delta+1}}$ is $L^2$.
\end{definition}
The hermitian metric $h$ on the vector bundle $\mathbb{E}$ induces a metric on $\End (\mathbb E)$ and $\End (\mathbb E)\otimes \Omega^1$. The definition of the spaces $L^2_{\delta}$ extends to section of such bundles using the induced metric. There is an orthogonal decomposition
\begin{equation}\label{eq_end_E}
\End(\mathbb{E}) = \End(\mathbb{E})_0 \oplus \End(\mathbb{E})_1
\end{equation}
with $\End( \mathbb{E})_0$ the space of endomorphism commuting with $A$. It induces an orthogonal decomposition
\[
\Omega^1 \otimes \End( \mathbb{E}) = \left( \Omega^1 \otimes \End (\mathbb{E})_0 \right) \oplus \left( \Omega^1 \otimes \End (\mathbb{E})_1 \right)
\]
For $f\in \Omega^1\otimes \End( \mathbb{E}) $ this orthogonal decomposition reads
\[
f=f_0 +f_1.
\]

\begin{definition}[Sobolev spaces $L^{k,2}_{\delta}$]
\[
L^{k,2}_{\delta}\left(\Omega^1\otimes\End (\mathbb{E}) \right)  :=\left\lbrace f \in L^2_{\delta} \left| \nabla^j f_0, \;  \frac{\nabla^j f_1}{r^{k-j}} \in L^{2}_{\delta} \; \text{ for }0\le j \le k\right. \right\rbrace 
\]
with $\nabla$ the covariant derivative with respect to the unitary connection $D_0^h$. 
\end{definition}

\begin{definition}[Space of admissible connections]\label{def_admissible_connections_loc}
The space of admissible connections is 
\[
\mathcal{A} =\left\lbrace D_0+a \left| a \in L^{1,2}_{-2+\delta} \left(\Omega^1\otimes \End (\mathbb{E})\right)\right.\right\rbrace.
\]
\end{definition}
\begin{remark}
Note that the space of admissible connections is chosen so that the connection $D=D_0+a$ introduced at the beginning of this section is admissible. Indeed, in the orthonormal trivialization  $\left(e_j\right)_{1\le j \le n}$, the matrix $a$ is strictly block upper triangular. The non zero coefficients strictly above the diagonal have the following form
\[
|z|^{\beta_i-\beta_j} \frac{a_{i,j}}{z}
\]
with $\beta_i > \beta_j$ and $a_{i,j}$ constant. Thus $a\in L^{1,2}_{-2+\delta}$ for small enough parameter: \[
0<\delta < \beta_i-\beta_j.
\]
\end{remark}

\subsubsection{Variation of the stability parameters and the metric}\label{subsect_local_variation_metric}
In order to pursue the path announced in Diagram \eqref{diagram_path_ch4}, slightly modify the stability parameter $\alpha$ to a parameter $\widetilde{\alpha}$, a diagonal matrix with coefficients
\[
\left(\widetilde{\alpha}_1,\widetilde{\alpha}_2,\dots\right)=\left( \underbrace{\widetilde{\alpha}_{1,1},\dots, \widetilde{\alpha}_{1,1}}_{{\mu^1_1}' }, \underbrace{\widetilde{\alpha}_{1,2},\dots, \widetilde{\alpha}_{1,2}}_{ {\mu^1_2}'}, \dots \right)
\]
The associated metric $\widetilde{h}$ is defined such that the holomorphic trivialization $(f_j)_{1\le j\le n}$ of the holomorphic bundle $F$ is orthogonal and
\[
|f_j|_{\widetilde{h}} =|z|^{\widetilde{\alpha}_j}.
\]
This provide an hermitian bundle with orthonormal trivialization $(\widetilde{e}_j)_{1\le j\le n}$ defined by
\[
\widetilde{e}_j = \frac{f_j}{|z|^{\widetilde{\alpha}_j}}.
\]
We follow the same process as before in the opposite direction. $D_0^{\widetilde{h}}$ is the $\widetilde{h}$-unitary connection with $(0,1)$-part $\overline{\partial}^F$. And 
\[
\widetilde{\Phi}_0 := \theta_0 + \theta_0^{\widetilde{\dagger}}
\]
the adjoint is taken with respect to the metric $\widetilde{h}$. Then
\[
\widetilde{D}_0 := D_0^{\widetilde{h}} + \widetilde{\Phi}_0.
\]
In the trivialization $(\widetilde{e}_j)_{1\le j\le n}$ it reads
\begin{eqnarray*}
\widetilde{\Phi}_0 &=& \frac{1}{2}\left(A-\beta\right)\frac{d z}{z} + \frac{1}{2}\left(A^{\widetilde{\dagger}}-\beta\right)\frac{d \overline{z}}{\overline{z}} \\
D_0^{\widetilde{h}} &=& d + \frac{1}{2}\widetilde{\alpha} \left(\frac{d z}{z}-\frac{d\overline{z}  }{\overline{z}  }\right)
\end{eqnarray*}
Setting $\widetilde{A}= \widetilde{\alpha}+i \Imm A$ and $\widetilde{\beta}= \beta + \widetilde{\alpha} - \alpha$ we obtain a canonical form like in Notations \ref{notations_canonical_form}
\begin{eqnarray*}
D_0^{\widetilde{h}}= d + \frac{1}{2}\Ree(\widetilde{A})\left(\frac{d z}{z}-\frac{d \overline{z}}{\overline{z}} \right) \\
\widetilde{\Phi}_0 = \frac{1}{2}\left( \widetilde{A}\frac{d z}{z} + \widetilde{A}^{\widetilde{\dagger}}\frac{d\overline{z}}{\overline{z}} - \widetilde{\beta} \frac{d z}{z} - \widetilde{\beta} \frac{d \overline{z}}{\overline{z}}    \right) 
\end{eqnarray*}
Continuing in the opposite direction, the $(0,1)$-part of $\widetilde{D}_0$ defines an holomorphic bundle $\widetilde{E}$ with holomorphic trivialization $(\widetilde{\tau}_j)_{1\le j\le n}$
\[
\widetilde{\tau}_j:=|z|^{\widetilde{\beta}_j- i \Imm \widetilde{A}_j} \widetilde{e}_j.
\]
$\widetilde{D}_0$ defines a logarithmic connection on $\widetilde{E}$, in the trivialization $(\widetilde{\tau}_j)_{1\le j\le n}$ it reads
\[
\widetilde{D}_0 = d + \widetilde{A}\frac{d z}{z}
\]
and $\widetilde{A}$ has distinct eigenvalues on each graded of the filtration of type $\nu$ and so does the monodromy of the local system of flat sections. 

Let us summarize the local behaviour for Diagram \eqref{diagram_path_ch4} in terms of residue. We look at a particular block of size $\nu_j$. The stability parameter associated to the graded of the filtration is specified with over-brace. N.A.H stands for non-Abelian Hodge theory.
\begin{equation*}
\begin{tikzcd}
 \begin{pmatrix}
\overbrace{A_j \Id_{ {\mu^j_{1}}' }   }^{\beta_{j,1}} &  *& \\
0 & \overbrace{A_j  \Id_{ {\mu^j_{2}}' }  }^{\beta_{j,2}} & *\\
\vdots  &  0 & \ddots 
\end{pmatrix}  \arrow[r,"N.A.H"] & 
\left(
\overbrace{
\begin{array}{ccccc}
   (A_j-\beta_{j,1})  \Id_{  {\mu^j_{1}}' }  & * &   \\
   0  &  (A_j-\beta_{j,2} )  \Id_{ {\mu^j_{2}}'  } & *  \\
 \vdots & 0 & \ddots 
\end{array}
}^{\alpha^1_j }
\right) \arrow[d,"\alpha \mapsto \widetilde{\alpha}"] \\
\begin{pmatrix}
\overbrace{\widetilde{A}_{j,1} \Id_{  {\mu^j_{1}}'  }}^{\widetilde{\beta}_{j,1}} &  *& \\
0 & \overbrace{\widetilde{A}_{j,2} \Id_{ {\mu^j_{2}}'  }}^{\widetilde{\beta}_{j,2}} & *\\
\vdots  &  0 & \ddots 
\end{pmatrix} & \arrow[l,"N.A.H"]
 \begin{pmatrix}
 \overbrace{\left(A_j-\beta_{j,1}\right) \Id_{ {\mu^j_{1}}' }}^{\widetilde{\alpha}_{j,1}} &  *& \\
0 & \overbrace{\left(A_j-\beta_{j,2}\right)  \Id_{{\mu^j_{2}}'  }}^{\widetilde{\alpha}_{j,2}} & *\\
\vdots  &  0 & \ddots 
\end{pmatrix} 
\end{tikzcd}
\end{equation*}
With $\widetilde{A}_{j,i}= \widetilde{\alpha}_{j,i}+ i \Imm A_j$ and $\widetilde{\beta}_{j,i}=\beta_{j,i}+\widetilde{\alpha}_{j,i}-\alpha^1_j$.

\subsection{Diffeomorphism between moduli spaces}\label{sect_diffeo_ms}
\subsubsection{Analytic construction of the moduli spaces}\label{sect_global_moduli}
The local study on the disk actually extends to global moduli spaces for objects defined over punctured Riemann surfaces. Analytic construction of moduli spaces relies on methods from Kuranishi \cite{Kuranishi},  Atiyah-Hitchin-Singer \cite{AHSinger} and Atiyah-Bott \cite{Atiyah_Bott}. In this section we recall the analytic construction of the moduli spaces involved in the parabolic version of non-Abelian Hodge theory. Some particular cases of those moduli spaces were constructed by Konno \cite{Konno} and Nakajima \cite{Nakajima_1996}. However we need more general construction in order to allow not necessarily central action of the residues of the Higgs fields on the graded of the filtration. The construction we follow is the one from Biquard-Boalch \cite{biquard_boalch_2004}. Note that a larger family of groups was considered by Biquard, García-Prada, Mundet i Riera \cite{Biquard_Garcia_Prada_Riera}.

The local canonical model introduced in \ref{notations_canonical_form} is used to represent behaviour of connections near the punctures $p_j$. Let $\mathbb{E}$ a vector bundle on $\Sigma$ endowed with an hermitian metric $h$. Notation $\mathbb{E}$ refers to a vector bundle from differential geometry point of view whereas $E$ refers to holomorphic bundle. Let $D_0$ a model connection such that on the neighborhood of the punctures it coincides with the local model connection of previous subsection. The connection decomposes as
\[
D_0 = D_0^h + \Phi
\]
with $D_0^h$ unitary and $\Phi$ self-adjoint with respect to the metric $h$. We assume for this model connection that in an orthonormal trivialization $(e_i)_{1\le i \le n}$ of $\mathbb{E}$ near the puncture $p_j$:
\[
D^h_0 = d + \frac{1}{2} \Ree(A^j) \left(\frac{d z}{z} - \frac{d \zbar}{\zbar}\right)
\]
and
\[
\Phi = \frac{1}{2} \left( A^j\frac{d z}{z} + (A^j)^{\dagger}\frac{d\overline{z}}{\overline{z}} - \beta^j \frac{d z}{z} - \beta^j \frac{d \overline{z}}{\overline{z}}\right)
\]
with $A^j$ and $\beta^j$ the residue and the stability parameter for the de Rham moduli space at the puncture $p_j$. They correspond to the local parameter $A$ and $\beta$ from Section \ref{subsect_local_model}, they are constant diagonal matrices. The parameters of the de Rham moduli space are chosen so that it corresponds under Riemann-Hilbert correspondence to a resolution of a character varieties with generic monodromies $\mlps$. Therefore connections with such polar parts are necessarily irreducible. 

Take $r$ a function strictly positive on the punctured Riemann surface $\Sigma^0$ such that it coincides with the radial coordinate near each punctures. The global weighted Sobolev space is defined as the local one from \ref{subsect_local_weighted} with this positive function $r$. It is still denoted by $L^{k,2}_{\delta}\left(\Omega^1\otimes \End{(\mathbb{E})}\right)$. The space of admissible connections is 
\[
\mathcal{A} = \left\lbrace D_0 + a \left| a \in L^{1,2}_{-2+\delta} \left( \Omega^1\otimes \End(\mathbb{E}) \right) \right.\right\rbrace .
\]
This affine space is actually endowed with various complex structures. Decomposing according to $(1,0)$-part and $(0,1)$-part $a=a^{1,0}+a^{0,1}$
\[
I.a = i a
\]
and
\[
J.a = i(a^{0,1})^{\dagger}-i(a^{1,0})^{\dagger}
\]
The curvature of an admissible connection $D=D_0 + a$  is denoted by $F_D$. Consider the complex gauge group
\[
\mathcal{G}^I=\left\lbrace g\in\Aut(\mathbb{E}) \left| (D_0^h g) g\inv , \; g \Phi_0 g\inv \in L^{1,2}_{-2+\delta} \right. \right\rbrace 
\]
It acts on $ \mathcal{A}$ by
\[
g.D := g D g\inv = D - \left(D g \right)g\inv .
\]
Next theorem gives an analytic construction of the set of isomorphism classes of parabolic flat connections with prescribed polar part. Later on, this set will be endowed with a manifold structure.
\begin{theorem}[Biquard-Boalch \cite{biquard_boalch_2004} Section 8]
The de Rham moduli space of stable flat connection with prescribed polar part on the graded part of the filtration introduced in \ref{subsect_dR_alg} is the following set
\[
\mathcal{M}^{d R}_{A, \beta} = \left\lbrace D_0 + a \in \mathcal{A} \left| F_D =0 \right. \right\rbrace/ \mathcal{G}^{I}.
\]
With $F_D$ the curvature of $D=D_0+a$. The stability condition does not appear as it is imposed by the generic choice of eigenvalues of the residue of $D_0$.
\end{theorem}

Now starting from $ D=D_0 +a \in \mathcal{A}$ there is a natural candidate to produce a parabolic Higgs bundle, like in the local model. First decompose $D$ in a unitary part and a self-adjoint part
\begin{eqnarray*}
D &=&D^h + \Phi \\
&=& D_0^h +\frac{a-a^{\dagger} }{2} + \Phi_0 + \frac{ a + a^{\dagger}  }{2}
\end{eqnarray*}
The natural candidate for the underlying holomorphic structure of the parabolic Higgs bundle is, in the orthonormal trivialization $(e_j)_{1\le j \le k}$
\[
\overline{\partial}^E = \overline{\partial} - \frac{1}{2} \Ree(A) \frac{d\overline{z}}{\overline{z}} + \frac{a^{0,1}-(a^{1,0})^{\dagger}}{2}.
\]
and the Higgs field
\[
\theta = \theta_0 + \frac{a^{1,0} + (a^{0,1})^{\dagger}}{2}.
\]
This data provides a Higgs bundle if $\overline{\partial}^E \theta =0$, equivalently if the pseudo curvature $G_D$ vanishes. Note that the complex structure $J$ is compatible with the Higgs bundles point of view. Indeed if $\theta$ is the Higgs field associated to $D$ then $i \theta$ is the Higgs field associated to $J.D$. The complex gauge group acts on the Higgs bundles structures by
\[
g.(\overline{\partial}^E,\theta ) :=(g \overline{\partial}^E g\inv ,g \theta g\inv ).
\]

Next theorem gives an analytic construction of the set of isomorphism classes of parabolic Higgs bundles with prescribed residue. Later on, this set will be endowed with a manifold structure.
\begin{theorem}[Biquard-Boalch \cite{biquard_boalch_2004} Section 7]
The Dolbeault moduli space of stable parabolic Higgs bundles with prescribed polar part on the graded part of the filtration introduced in \ref{subsect_dol_alg} is the following set
\[
\mathcal{M}^{Dol}_{B, \alpha} = \left\lbrace D_0 + a \in \mathcal{A} \left| \overline{\partial}^E \theta =0 \right. \right\rbrace/ \mathcal{G}^{J}.
\]
The stability condition does not appear as it is imposed by the generic choice of eigenvalues of the residue. As a group $\mathcal{G}^J$ is just $\mathcal{G}^I$, we change the upper index to precise which action is considered, the $I$-linear action or the $J$-linear action.
\end{theorem}

The non-Abelian Hodge theory gives a correspondence between Dolbeault and de Rham moduli spaces. The parameters are intertwined as in the local model. A nice way to state this correspondence is with hyperkähler geometry. Introduce the unitary gauge group
\[
\mathcal{G}=\left\lbrace g\in U(\mathbb{E}) \left| \left( D_0 g\right) g\inv \in L^{1,2}_{-2+\delta} \right. \right\rbrace .
\]
Consider the moduli space
\[
\mathcal{M}=\left\lbrace D\in\mathcal{A} \left|\; \overline{\partial}^E \theta=0, \; F_D=0  \right. \right\rbrace/\mathcal{G}.
\]

The equations defining $\mathcal{M}$ can be interpreted as vanishing of an hyperkähler moment map. Then the moduli space $\mathcal{M}$ is an  hyperkähler reduction as in \cite{hklr}.
\begin{theorem}[Biquard-Boalch \cite{biquard_boalch_2004} Theorem 5.4]
The moduli space $\mathcal{M}$ carries an hyperkähler manifold structure.
\end{theorem}
\begin{proof}
The deformation theory for the moduli space $\mathcal{M}$ at a point $[D]$ is encoded in the following complex
\[
\begin{tikzcd}
L^{2,2}_{-2+\delta} \left(\mathfrak{u}(\mathbb{E}) \right) \arrow[r,"D"] & L^{1,2}_{-2+\delta} \left(\Omega^1\otimes\End \mathbb{E}\right) \arrow[r,"D+D^*"] & L^{2}_{-2 +\delta} \left(\left(\Omega^2\otimes \End \mathbb{E}  \right)\oplus i\mathfrak{u}(\mathbb{E})\right)
\end{tikzcd}.
\]
$D^*$ is the formal adjoint of $D$ with respect to the $L^2$ inner product and the metric $h$. The analytic study of this complex is detailed in \cite{biquard_boalch_2004}. Its first cohomology group is represented by the harmonic space $\boldsymbol{H}^1 \subset L^{1,2}_{-2+\delta} \left(\Omega^1\otimes\End \mathbb{E}\right)  $. The Kuranishi slice at $[D]$ is defined by
\begin{equation}\label{eq_kuranishi_slice}
\mathcal{S}_D :=\left\lbrace D+a \left| \Imm( D^* a )= 0, \; G_{D+a}=0, \; F_{D+a}=0\right. \right\rbrace.
\end{equation}
Taking a small enough neighborhood of $D$ in the Kuranishi slice, one obtains a finite dimensional manifold transverse to the $\mathcal{G}$-orbits. The Kuranishi map provides an isomorphism between a neighborhood of $0$ in $\boldsymbol{H}^1$ and a neighborhood of $D$ in the Kuranishi slice, see Konno \cite[Lemma 3.8, Theorem 3.9]{Konno}. This provides an hyperkähler manifold structure on the moduli space. 

\end{proof}

Now the non-Abelian Hodge theory can be described the following way.
\begin{theorem}[Biquard-Boalch \cite{biquard_boalch_2004} Theorem 6.1]
The manifold $\mathcal{M}$ endowed with the complex structure $I$ is the moduli space $\mathcal{M}^{d R}_{A,\beta}$.

The manifold $\mathcal{M}$ endowed with the complex structure $J$ is the moduli space $\mathcal{M}^{Dol}_{B,\alpha}$.
\end{theorem}

\subsubsection{Construction of the diffeomorphisms}
\begin{theorem}[Riemann-Hilbert correspondence]
The moduli space $\mathcal{M}^{d R}_{A,\beta}$ is complex analytically isomorphic to a resolution of character varieties $\mlps$.
\end{theorem}
\begin{proof}
As explained in \ref{proposition_filtered_moduli_resolu}, $\mlps$ is nothing but the moduli space of filtered local systems with prescribed graded part of the monodromy around the punctures. Filtered version of the Riemann-Hilbert correspondence is established as an equivalence of category by Simpson \cite{Simpson}. Yamakawa \cite{Yamakawa} proved that it is a diffeomorphism using a particular construction of the de Rham moduli space from Inaba \cite{Inaba_}. The same argument apply with the de Rham moduli space endowed with the manifold structure from $\mathcal{M}$. Starting from a flat connection, the associated local system is obtained by taking flat sections i.e. solving a differential equation. When the parameters of the equation vary complex analytically, so does the solution.
\end{proof}

Then $\mathcal{M}^{dR}_{A,\beta}$ and $\mathcal{M}^{Dol}_{B,\alpha}$ are diffeomorphic as both are $\mathcal{M}$ with a particular complex structure. The first line in the path announced in Diagram \ref{diagram_path_ch4} is now constructed. The second line is obtained exactly like the first, but in the other direction. It remains to describe the vertical arrow between two Dolbeault moduli spaces $\mathcal{M}^{Dol}_{B,\alpha}$ and $\mathcal{M}^{Dol}_{B,\widetilde{\alpha}}$. This is given by Biquard, García-Prada, Mundet i Riera \cite[Theorem 7.10]{Biquard_Garcia_Prada_Riera}. The construction of the diffeomorphism is detailed in the remaining of the section.

Because of genericity of the eigenvalues of the residue, the stability parameter $\alpha$ is irrelevant. The parameter $\alpha$ can be changed to a stability parameter $\widetilde{\alpha}$ with different values for each graded of the filtration. Namely one can chose $\widetilde{\alpha}$ such that the associated matrix satisfies $Z_{\GL_n}(\widetilde{\alpha}^{i})=Z_{\GL_n}(B^i)$ and such that the parabolic degree remains $0$. The local behaviour near each puncture is described by the right hand side of the diagram at the end of \ref{subsect_local_variation_metric}.  

We introduce the following notation
\[
\epsilon_i := \widetilde{\alpha}_i-\alpha_i.
\]
For the construction of the diffeomorphism in Theorem \ref{th_diffeo_var_alpha}, it will be conveniant to assume
\[
\max_{i,j}|\epsilon_i-\epsilon_j| < \delta 
\]
with $\delta$ the parameter appearing in the weighted Sobolev space $L^{1,2}_{-2+\delta}$. 
\begin{proposition}\label{prop_bijection_dol}
For such choice of parameter there is a natural bijection between $\mathcal{M}^{Dol}_{B,\alpha}$ and $\mathcal{M}^{Dol}_{B,\widetilde{\alpha}}$.
\end{proposition}
\begin{proof}
$\mathcal{M}^{Dol}_{B,\alpha}$ classifies isomorphism classes of parabolic Higgs bundles with parabolic structure at $p_j$
\[
0= F^j_0 \subsetneq F^j_1 \subsetneq \dots \subsetneq F^j_{n_j}=F^j
\]
and with the residue of the Higgs fields preserving this filtration and acting as a semisimple endomorphism $B^j_i$ on the graded spaces
\[
F^j_i / F^j_{i-1}.
\]
Such spaces decomposes as direct sum of eigenspaces for $B^j_i$. After ordering the eigenvalues, we obtain a uniquely determined refinement of the initial parabolic structure:
\[
0 = \widetilde{F}^j_0 \subsetneq \widetilde{F}^j_1 \subsetneq \dots \subsetneq \widetilde{F}^j_{m_j} = F^j.
\]
Then the residue of the Higgs field acts as a central endomorphism on the graded $\widetilde{F}^j_i / \widetilde{F}^j_{i-1}$.
This gives a map $f:\mathcal{M}^{Dol}_{B,\alpha}\to \mathcal{M}^{Dol}_{B,\widetilde{\alpha}}$. Stability is not an issue as the polar part of the residue is generic. The map forgetting part of the filtration is an inverse so that there is a natural bijection between both moduli spaces.
\end{proof}

Before proving that this bijection is a diffeomorphism the manifold structure on $\mathcal{M}_{B,\widetilde{\alpha}}^{Dol}$ is detailed. It is constructed just like $\mathcal{M}^{Dol}_{B,\alpha}$ but with different parameters.

Similarly to $\mathcal{M}$, construct a moduli space $\mathcal{M}_{\widetilde{h}}$. Instead of the initial metric $h$, we use a metric $\widetilde{h}$, similar to the local model from \ref{subsect_local_variation_metric}. Namely it is chosen so that near each puncture it admits as an orthonormal trivialization $(\widetilde{e}_i)_{1\le i \le n}$ with
\[
\widetilde{e}_i = r^{\epsilon_i} e_i.
\]
Where $(e_i)_{1\le i \le n}$ is the orthnormal trivialization with respect to $h$ near the puncture and $\epsilon_i=\widetilde{\alpha}_i-\alpha_i$.

First we construct $\widetilde{D}_0$, a starting point to construct an affine space of admissible connections. Recall that
\[
D_0 = D_0^h +\Phi_0
\]
with $D_0^h$ a $h$-unitary connection and $\Phi_0$ self-adjoint with respect to $h$. Take ${D^h_0}''$ the $(0,1)$-component of $D_0^h$ and $\Phi_0^{1,0}$ the $(1,0)$-component of $\Phi_0$. There exists a unique ${D^{\widetilde{h}}_0}'$ such that ${D^{\widetilde{h}}_0}'+{D^h_0}''$ is $\widetilde{h}$-unitary. Let ${\Phi_0^{1,0}}^{\widetilde{\dagger}}$ the adjoint of $\Phi_0^{1,0}$ with respect to the metric $\widetilde{h}$. Then $\widetilde{D}_0$ is defined by
\[
\widetilde{D}_0 := {D^{\widetilde{h}}_0}'+{D^h_0}''+ \Phi_0^{1,0} + {\Phi_0^{1,0}}^{\widetilde{\dagger}}.
\]
Near the puncture, in the trivialization $(\widetilde{e}_i)_{1\le i \le n}$, the connection $\widetilde{D}_0$ behaves exactly like the local model with the same name introduced in \ref{subsect_local_variation_metric}. Define the affine space of admissible connections with respect to $\widetilde{D}_0$ and the metric $\widetilde{h}$.
\[
\mathcal{A}_{\widetilde{h}}:=\left\lbrace \widetilde{D}_0  + \widetilde{a} \left| \right. \widetilde{a} \in L^{1,2}_{-2+\widetilde{\delta}}\left( \Omega^1 \otimes \End(\mathbb{E})\right)  \right\rbrace
\]
The weighted Sobolev space $L^{1,2}_{-2+\widetilde{\delta}} \left(\Omega^1\otimes \End(\mathbb{E})\right)$ is also defined using the metric $\widetilde{h}$.

Moreover notice that we do not chose the same parameter $\delta$ for $\mathcal{A}$ and for $\mathcal{A}_{\widetilde{h}}$. It will be convenient to chose $\widetilde{\delta}$ such that
\begin{equation}\label{eq_delta_t_small}
0 < \widetilde{\delta} < \delta - \max_{i,j} |\epsilon_i-\epsilon_j|.
\end{equation}
With this set up, we are ready to prove that the bijection from previous proposition is a diffeomorphism. 
\begin{theorem}\label{th_diffeo_var_alpha}
The natural bijection between $\mathcal{M}^{Dol}_{B,\alpha}$ and $\mathcal{M}^{Dol}_{B,\widetilde{\alpha}}$ is a diffeomorphism.
\end{theorem}
\begin{proof}
$\mathcal{M}^{Dol}_{B,\alpha}$ is identified with the manifold $\mathcal{M}$ with the complex structure $J$.

Take an element in $\mathcal{M}^{Dol}_{B,\alpha}$ identified with an element $[D]\in\mathcal{M}$. $[D]$ is the class of $D=D_0+a$ an admissible connection with vanishing curvature and pseudo-curvature. By construction of the manifold structure, a neighborhood of $[D]$ in $\mathcal{M}$ is diffeomorphic with a neighborhood of $D$ in the Kuranishi slice $\mathcal{S}_{D}$ defined in \eqref{eq_kuranishi_slice}. We shall prove that the bijection from Proposition \ref{prop_bijection_dol} induces a smooth map from a neighborhood of $D$ in $\mathcal{S}_{D}$ to $\mathcal{A}_{\widetilde{h}}$.

First we describe the image of the connection $D$, it is obtained exactly the same way $\widetilde{D}_0$ is obtained from $D_0$. It decomposes as a connection $h$-unitary plus a hermitian part
\[
D= D_{0}^h+ \frac{a-a^{\dagger}}{2}  + \Phi^0 +\frac{a + a^{\dagger}}{2}.
\]
It can be decomposed further in components of type $(1,0)$ and $(0,1)$. Then the $(0,1)$-component of the $h$-unitary part is
\[
\overline{\partial}^F={D_{0}^h}''+ \frac{a^{0,1}-{a^{1,0}}^{\dagger} }{2}
\]
and the $(1,0)$-component of the self-adjoint part is
\[
\theta = \Phi^{1,0} + \frac{a^{1,0}+{a^{0,1}}^{\dagger} }{2}.
\]
The parabolic Higgs bundle associated to $D$ is $(\overline{\partial}^F,\theta)$. Now we switch to the metric $\widetilde{h}$. Near each puncture, in the $\widetilde{h}$-orthonormal trivialization $(\widetilde{e}_i)_{1\le i\le n}$
\begin{eqnarray*}
\overline{\partial}^F= {D_0^h}'' +\left(\frac{\widetilde{\alpha}-\alpha}{2}\right)\frac{d\overline{z}}{\overline{z}}  + \widetilde{H} \frac{a^{0,1}-{a^{1,0}}^{\dagger} }{2} \widetilde{H} \inv
\end{eqnarray*}
and
\[
\theta = \phi^{1,0} + \widetilde{H}\frac{a^{1,0}+{a^{0,1}}^{\dagger} }{2} \widetilde{H}\inv.
\]
with $\widetilde{H}$ a diagonal matrix with coefficients $r^{\epsilon_i}$. Using the metric $\widetilde{h}$ we construct $D_{\widetilde{h}}'$ such that $D_{\widetilde{h}}'+\overline{\partial}^F$ is $\widetilde{h}$-unitary. And $\theta^{\widetilde{\dagger}}$ the adjoint of $\theta$ with respect to $\widetilde{h}$. We want to prove that 
\[
D_{\widetilde{h}}'+\overline{\partial}^F+\theta + \theta^{\widetilde{\dagger}}
\]
belongs to the space of admissible connections $\mathcal{A}_{\widetilde{h}}$. Let
\[
\widetilde{a}:=D_{\widetilde{h}}'+\overline{\partial}^F+\theta + \theta^{\widetilde{\dagger}} - \widetilde{D}_0.
\]
Components of $\widetilde{a}$ are obtained from components of $a$ by multiplication by $r^{\epsilon_i-\epsilon_j}$. Thus for $\widetilde{\delta}$ small enough \eqref{eq_delta_t_small}, $\widetilde{a}$ belongs to $L^{1,2}_{-2+\widetilde{\delta}}$. Therefore the bijection from $\mathcal{M}^{Dol}_{B,\alpha}$ to $\mathcal{M}^{Dol}_{B,\widetilde{\alpha}}$ comes from a map
\[
\begin{array}{ccc}
    \left\lbrace D_0 + a \in \mathcal{A} \left|\; F_{D_0+a} = G_{D_0+a}= 0 \right. \right\rbrace & \to &  \left\lbrace \widetilde{D}_0 + \widetilde{a} \in \mathcal{A}_{\widetilde{h}} \left| \; G_{\widetilde{D}_0 + \widetilde{a}}= 0 \right. \right\rbrace \\
     D_0 + a & \mapsto & \widetilde{D}_0 + \widetilde{a}. 
\end{array}
\]
This restricts  to a diffeomorphism from a neighborhood of $D$ in the Kuranishi slice $\mathcal{S}_D$ to a  manifold transverse to the $\mathcal{G}^J$-orbits in a neighborhood of $D$. Therefore the map $\mathcal{M}_{B,\alpha}^{Dol}\to \mathcal{M}_{B,\widetilde{\alpha}}^{Dol}$ is a diffeomorphism.
\end{proof}

To finish, let us detail the last step at the bottom left corner of Diagram \eqref{diagram_path_ch4}. Applying successively non-Abelian Hodge theory and Riemann-Hilbert correspondence, the moduli space $\mathcal{M}_{B,\widetilde{\alpha}}^{Dol}$ is diffeomorphic to a moduli space of filtered local system $\widetilde{\mathcal{M}}_{\boldsymbol{L},\boldsymbol{P},\boldsymbol{\widetilde{\sigma}}}$. The parameters are such that  $Z_{\GL_n}(\widetilde{\sigma}^j) = L^j$ for $1\le j \le k$. The map $p^{\boldsymbol{\widetilde{\sigma}}}: \widetilde{\mathcal{M}}_{\boldsymbol{L},\boldsymbol{P},\boldsymbol{\widetilde{\sigma}}} \to \mathcal{M}_{\boldsymbol{\mathcal{S}}}$ from \ref{def_resol_char} is an isomorphism. $\mathcal{M}_{\boldsymbol{\mathcal{S}}}$ is the character variety with monodromy at the puncture $p_j$ in $\mathcal{S}_j$ the conjugacy class of $\widetilde{\sigma}^j$. Theorem \ref{th_diffeo_resol_ss} is proven.

\section{Computation of the intersection cohomology of character varieties}\label{sect_computation_pc}

In this section the Poincaré polynomial for intersection cohomology of character varieties with the closure of conjugacy classes of any Jordan type at each puncture:
\[
P_c \left(\mathcal{M}_{\overline{{\mathcal{C}}}_{\boldsymbol{\mu},\boldsymbol{\sigma}}},v\right) = v^{d_{\boldsymbol{\mu}}}\left\langle s_{\boldsymbol{\mu'}}, \mathbb{H}_n^{H L V}(-1,v) \right\rangle .
\]

This proves the Poincaré polynomial specialization of a conjecture from Letellier \cite{zclosure}. The idea is to express the intersection cohomology of character varieties in terms of usual cohomology of resolutions of character varieties. We proved in previous section that such resolution are diffeomorphic to semisimple character varieties. We conclude as the cohomology of semisimple character varieties was computed by Mellit \cite{mellit_p}. 

\subsection{Computation of the Poincaré polynomial}\label{sect_poincare_comput}

Consider a generic $k$-uple of conjugacy classes $\mathcal{C}_{\boldsymbol{\mu},\boldsymbol{\sigma}} =  \left(\mathcal{C}_{\underline{\mu}^1,\sigma^1}, \dots, \mathcal{C}_{\underline{\mu}^k,\sigma^k} \right)$. As usual, the class $\mathcal{C}_{\underline{\mu}^j,\sigma^j}$ is characterized by its eigenvalues 
\[
\underbrace{\sigma^j_1,\dots,\sigma^j_1}_{\nu^j_1},\dots,\underbrace{\sigma^j_{l_j},\dots, \sigma^j_{l_j}}_{\nu^j_{l_j}}  
\]
and by $\mu^{j,i}\in \Part_{\nu^j_i}$ the Jordan type of the eigenvalue $\sigma^{j}_i$. Denote by ${\mu^{j,i}}'$ the transposed partition. For each of this conjugacy classes consider the resolution of the closure (see \ref{subsect_resolution_conjugacy})
\[
\mathbb{\widetilde{X}}_{L^j, P^j, \sigma^j}\to \overline{\mathcal{C}}_{\underline{\mu}^j,\sigma^j} .
\]
The group $L^j$ used to construct the resolution is
\[
L^j \cong \underbrace{\GL_{ {\mu^{j,1}_1}'}\times \GL_{{\mu^{j,1}_2}'} \times \dots }_{\subset \GL_{\nu^j_1}}\times \dots \times \underbrace{\GL_{{\mu^{j,l_j}_1}'}\times \GL_{{\mu^{j,l_j}_2}'} \times \dots}_{\subset \GL_{\nu^j_{l_j}}}.
\]
As explained in \ref{subsect_resolution_character}, resolution of closure of conjugacy classes fit together in  $\mlps$ a resolution of the character variety $\mathcal{M}_{\overline{\mathcal{C}}_{\boldsymbol{\mu},\boldsymbol{\sigma}}}$.

Springer theory provides a combinatoric relation between the cohomology of $\mlps$ and intersection cohomology of character varieties $\mathcal{M}_{\overline{\mathcal{C}}_{\boldsymbol{\rho},\boldsymbol{\sigma}}}$ (see Theorem \ref{th_coho_resol_char})
\begin{equation}\label{eq_relation_cohom_resol}
H^{i+d_{\boldsymbol{\mu}}}_c
\left(\mlps,\kappa \right) 
\cong 
\bigoplus_{\boldsymbol{\rho}\in \Part_{\nu^1}\times\dots\times \Part_{\nu^k} }
A_{\boldsymbol{\mu'},\boldsymbol{\rho} }\otimes
I H_c^{i+d_{\boldsymbol{\rho}} } \left(\mathcal{M}_{\overline{\mathcal{C}}_{\boldsymbol{\sigma},\boldsymbol{\rho}}  } ,\kappa \right).
\end{equation}
This relation is the main tool allowing to go from usual cohomology of smooth varieties to intersection cohomology of singular varieties. In previous section (Theorem \ref{th_diffeo_resol_ss}) we saw that the resolution $\mlps$ is diffeomorphic to a character variety $\mathcal{M}_{\boldsymbol{\mathcal{S}}}$ with generic semisimple conjugacy classes at punctures. Precisely $\boldsymbol{\mathcal{S}}=\left(\mathcal{S}_1,\dots,\mathcal{S}_k \right)$ and $\mathcal{S}_j$ is the class of an element with centralizer in $\GL_n$ equal to $L^j\cong \GL_{{\underline{\mu}^j}'}$.

As the Poincaré polynomial is a topological invariant
\[
P_c \left( \mlps ,t \right)  = P_c\left(\mathcal{M}_{\boldsymbol{\mathcal{S}}},t \right).
\]
Let us translate \eqref{eq_relation_cohom_resol} in terms of Poincaré polynomial.
\begin{equation}\label{eq_combi_ss_resol}
t^{-d_{\boldsymbol{\mu}}} P_c\left(\mathcal{M}_{\boldsymbol{\mathcal{S}}},t \right) = \sum_{\boldsymbol{\rho\preceq\mu}} \left(\dim A_{\boldsymbol{\mu',\rho}}\right) t^{-d_{\boldsymbol{\rho}} } P_c\left( \mathcal{M}_{\overline{\mathcal{C}}_{\boldsymbol{\rho},\boldsymbol{\sigma}}} , t \right).
\end{equation}
The idea is now to invert this relation. First we compute the dimension of the multiplicity spaces  $\dim A_{\boldsymbol{\mu',\rho}}$.
\begin{lemma}\label{lemma_dim_Amurho}
The dimension of the multiplicity space is given by
\[
\dim A_{\boldsymbol{\mu',\rho}} = \prod_{  \substack{  1\le j \le k     \\ 1 \le i \le l_j   }  } \left\langle h_{{\mu^{j,i}}'}, s_{ {\rho^{j,i}}'  }   \right\rangle 
\]
\end{lemma}
\begin{proof}
By definition
\begin{eqnarray*}
A_{\boldsymbol{\mu',\rho}} &=& \Hom_{W_{\boldsymbol{M}}} 
\left(\Ind_{W_{\boldsymbol{L}}}^{W_{\boldsymbol{M}} } \epsilon_{\boldsymbol{\mu'} }, V_{\boldsymbol{\rho} }  \right) \\
& = & \bigotimes_{1 \le j \le k} \left(\bigotimes_{1 \le i \le l_j}\Hom_{\Sym_{\nu^j_i}}\left(\epsilon_{{\mu^{ j,i } }'  },V_{\rho^{ j,i }} \right)  \right).
\end{eqnarray*}
We conclude with Lemma \ref{lemma_ind_sign}.
\end{proof}

\begin{theorem}\label{th_sect_poincare}
For a generic $k$-uple of conjugacy classes ${\mathcal{C}}_{\boldsymbol{\mu},\boldsymbol{\sigma}}$, the Poincaré polynomial for compactly supported intersection cohomology of the character variety $\mathcal{M}_{\overline{\mathcal{C}}_{\boldsymbol{\mu},\boldsymbol{\sigma}}}$ is
\[
P_c \left(\mathcal{M}_{\overline{\mathcal{C}}_{\boldsymbol{\mu},\boldsymbol{\sigma}}},v\right) = v^{d_{\boldsymbol{\mu}}}\left\langle s_{\boldsymbol{\mu'}}, \mathbb{H}_n^{H L V}(-1,v) \right\rangle .
\]
\end{theorem}
\begin{proof}
The complete symmetric functions $(h_{\mu})_{\mu\in \Part_m}$ and the Schur functions $(s_{\rho})_{\rho\in\Part_m}$ are two basis of the space of degree $m$ symmetric functions. Let $(M_{\mu,\rho})_{\mu,\rho\in\Part_m }$ the transition matrix between those basis then
\[
h_{\mu} = \sum_{\rho \in\Part_m} M_{\mu,\rho} s_{\rho}.
\]
As the Schur functions form an orthonormal basis, the transition matrix is given explicitly by
\[
M_{\mu,\rho} = \left\langle h_{\mu},s_{\rho} \right\rangle.
\]
It is invertible and denote by $(N_{\mu,\rho})_{\mu,\rho \in \Part_m }$ its inverse. Combining Equation \eqref{eq_combi_ss_resol}, Lemma \ref{lemma_dim_Amurho} and the formula for Poincaré polynomial of character varieties with semisimple conjugacy classes:
\[
\left\langle \prod_{j=1}^k\prod_{i=1}^{l_j} h_{ {\mu^{j,i}}'}[X_j] ,\mathbb{H}_n^{H L V}(-1,v) \right\rangle =  \sum_{\boldsymbol{\rho\preceq\mu}} \prod_{j = 1}^k \prod_{i= 1}^{l_j} \left\langle h_{ {\mu^{ j,i} }'  }, s_{ {\rho^{ j,i} }' }     \right\rangle    v^{-d_{\boldsymbol{\rho}} } P_c\left( \mathcal{M}_{\overline{\mathcal{C}}_{\boldsymbol{\rho},\boldsymbol{\sigma}}} , v\right). 
\]
This relation can now be inverted. Fix $\boldsymbol{\lambda}\in \Part_{\nu^1}\times\dots\times\Part_{\nu^k}$. Multiply previous equation by $N_{{\lambda^{1,1}}',{\mu^{1,1} }'}$ and sum over ${\mu^{1,1} }'\in \Part_{\nu^1_1}$. Repeating this process gives the expected result:
\[
\left\langle s_{\boldsymbol{\lambda}'}, \mathbb{H}_n^{H L V}(-1,v) \right\rangle = v^{-d_{\boldsymbol{\lambda}}} P_c\left( \mathcal{M}_{\overline{\mathcal{C}}_{\boldsymbol{\lambda},\boldsymbol{\sigma}}} , v \right).
\]
\end{proof}

\subsection{Weyl group action and twisted Poincaré polynomial}\label{subsect_twisted_poincare}

As in \cite[Proposition 1.9]{zclosure},  twisted Poincaré polynomial can be computed thanks to previous theorem. Using notations from \ref{subsect_relative_weyl_char} and Definition \ref{def_twisted_Poincare} for $\boldsymbol{\eta}$-twisted Poincaré polynomial we have the following theorem
\begin{theorem}\label{cor_twiste_Poincare_resol}
$\mathcal{C}_{\boldsymbol{\mu},\boldsymbol{\sigma}}$ is a generic $k$-uple of conjugacy classes and $\mlps$ is the resolution of $\mcmus$. For $\boldsymbol{\eta}$ indexing a conjugacy class in $W_{\boldsymbol{M}}(\boldsymbol{L})$, the $\boldsymbol{\eta}$-twisted Poincaré polynomial of $\mlps$ is
\[
 P_c^{\boldsymbol{\eta}}\left(\mlps,v\right) = (-1)^{r(\boldsymbol{\eta})}v^{d_{\boldsymbol{\mu}}}\left\langle \widetilde{h}_{\boldsymbol{\eta}},\mathbb{H}_n^{H L V}\left(-1,v \right)  \right\rangle 
\]
\end{theorem}
\begin{proof}
Theorem \ref{th_tr_resol} implies
\[
v^{-d_{\boldsymbol{\mu}}} P_c^{\boldsymbol{\eta}}\left(\mlps,v\right) = \sum_{\boldsymbol{\rho\preceq\mu}}  \left(\prod_{j=1}^k\prod_{i=1}^{l_j} c^{\rho^{j,i}}_{\eta^{j,i}}  \right) v^{-d_{\boldsymbol{\rho}}} P_c \left( \mcmur,t\right). 
\]
Apply Theorem \ref{th_sect_poincare}:
\begin{eqnarray*}
v^{-d_{\boldsymbol{\mu}}} P_c^{\boldsymbol{\eta}}\left(\mlps,v\right) = \sum_{\boldsymbol{\rho\preceq\mu}}  \left(\prod_{j=1}^k\prod_{i=1}^{l_j} c^{\rho^{j,i}}_{\eta^{j,i}}  \right) \left\langle s_{\boldsymbol{\rho'}}, \mathbb{H}_n^{H L V}\left(-1,v \right) \right\rangle. 
\end{eqnarray*}
Then using the relation $c_{\omega}^{\mu}=(-1)^{r(\omega)} c_{\omega'}^{\mu'}$ (see Lemma \ref{lemma_transpose_r_omega}) and Notations \ref{notations_k_h_tilde}
\[
v^{-d_{\boldsymbol{\mu}}} P_c^{\boldsymbol{\eta}}\left(\mlps,v\right) = (-1)^{r(\boldsymbol{\eta})}\left\langle \widetilde{h}_{\boldsymbol{\eta}},\mathbb{H}_n^{H L V}\left(-1,v \right)  \right\rangle 
\]
\end{proof}

Theorem \ref{th_diffeo_resol_ss}  gives a diffeomorphism between $\mlps$ and a character variety with semisimple monodromies $\mathcal{M}_{\boldsymbol{\mathcal{S}}}$. The diffeomorphism transports the action on the cohomology of $\mlps$ to an action on the cohomology of $\mathcal{M}_{\boldsymbol{\mathcal{S}}}$ and we have the following corollary.
\begin{corollary}
$W_{\boldsymbol{M}}(\boldsymbol{L})$ acts on the cohomology of $\mathcal{M}_{\boldsymbol{\mathcal{S}}}$ and the $\boldsymbol{\eta}$-twisted Poincaré polynomial is
\[
 P_c^{\boldsymbol{\eta}}\left( \mathcal{M}_{\boldsymbol{\mathcal{S}}} ,v\right) = (-1)^{r(\boldsymbol{\eta})}v^{d_{\boldsymbol{\mu}}}\left\langle \widetilde{h}_{\boldsymbol{\eta}},\mathbb{H}_n^{H L V}\left(-1,v \right)  \right\rangle.
\]
\end{corollary}

\bibliographystyle{alpha}
\bibliography{biblio_these}

\end{document}